%% file: main.tex
\documentclass[11pt]{article}
\input{preamble_report}

\def\R{\mathbb{R}}
\def\tr{\mathrm{tr}}

\def\eps{\epsilon}
\def\epss{\varepsilon}

\usepackage{subfigure}
\usepackage{caption}
\usepackage{enumitem}
\newboolean{isjournal}
\setboolean{isjournal}{true}

\newboolean{nocolors}
\setboolean{nocolors}{false}

\newboolean{usemtpro}
\setboolean{usemtpro}{false}

\ifthenelse{\boolean{nocolors}}{\hypersetup{colorlinks=false}}{}

\makeatletter
\newcommand\footnoteref[1]{\protected@xdef\@thefnmark{\ref{#1}}\@footnotemark}
\makeatother

\setauthA{{\bfseries Daniel Hill}$^{1*}$}

\setauthB{{\bfseries Martin Slawski}$^{2}$}

\begin{document}
\thispagestyle{firststyle}

\makemytitle{{\Large {\bfseries {Linear Regression from 1-bit Quantized Data}}}}
\vskip 3ex
%
{\large\begin{center}
\printA
\printB
\vskip1.5ex
{\scriptsize $^{1}$Department of Statistics, George Mason University, Fairfax, VA 22030, USA $\; \; \;$}  \\  {\scriptsize $^{2}$Department of Statistics, University of Virginia, Charlottesville, VA 22903, USA $\; \; \;$} \\
{\small \mymail{dhill29@gmu.edu}   
$\quad$ \mymail{ebh3ep@virginia.edu} } 
\end{center}}
\vskip 2ex

\begin{abstract}\vspace*{-.5ex}  Motivated by the prevalence of environments in which data is abundant while resources for storage and/or transmission might be scarce, we study linear regression when predictors, their squares, and responses 
are subject to single-bit dithered quantization. An estimator relying on plug-in estimation of the quadratic and linear terms in the quadratic program formulation of the least squares problem is proposed. We provide a non-asymptotic bound on the $\ell_2$-estimation error of this estimator and obtain its asymptotic distribution when the number of predictors is fixed, which can be used for inference and an investigation of the mean-square error efficiency relative to the ordinary least squares estimator. It is shown that for the quantization protocol under consideration, substantial improvements over the proposed estimator cannot be expected. A compression pipeline in which the underlying data is first subject to sketching and subsequently
quantization can be studied within our framework as well. We also present an extension to address high-dimensional predictors. Numerical experiments with synthetic data complement our theoretical findings. 
\end{abstract}

\section{Introduction}\label{sec:intro}
\vspace*{-1.3ex}
The information age has led to the generation of vast and unquantifiable amounts of data, but technology has struggled to keep pace with the growing demand for efficient storage and transmission. For example, Internet of Things (IoT) devices are limited in terms of power, bandwidth, and computational resources; as a result, they must reduce the data they store, transmit, and compute. In edge computing \cite[e.g.,][]{Shi2016, Abbas2017}, the bulk of incoming data needs to be processed directly by the local device receiving them since transmission to the cloud would yield unacceptable latencies, or is infeasible because of limited communication bandwidths. Increasing
energy demands of modern data centers have triggered an intense public debate \cite{turn1search18}, re-emphasizing the need for resource-efficient computing architectures. To this end, data and model compression techniques play an important role, reducing the footprint on the available computational resources while approximately preserving task-specific performance. Quantization, a traditional method of lossy data compression \cite[e.g.,][]{GershoGray1991, NeuhoffGray}, continues to be a widely used ingredient in data processing pipelines. In short, quantization maps high-precision measurements to values in a much coarser alphabet, thereby achieving savings in storage and transmission. In the past two decades, quantization has been a popular tool in compressed sensing \cite{LaskaBaraniuk2011, PlanVershynin2013a, BoufounosBaraniuk2008, Jacques2013a, Gopi2013, Baraniuk2017, Zymnis2009}, information retrieval \cite{Charikar2002, Li2014, Rane2013, SlawskiMitzenmacherLi2016, SlawskiLi2017}, distributed as well as federated optimization and machine learning \cite{Rabbat2005, alistarh2017qsgd, Wu2018, Beznosikov2023}, deep learning \cite{Gholami2022, Rokh2023, zhang2017zipml}, and edge computing \cite{EdgeBuilding, Lyu2018, Liu2020}. 

Despite the large body of work on the subject including various papers on theoretical aspects, there are only very few studies on the use of (coarse) quantization in the statistics literature. A notable exception is the paper \cite{Dirksen2022covariance} that considers covariance estimation -- a classical subject in multivariate and high-dimensional statistics -- under one-bit scalar quantization of the samples. The paper \cite{Dirksen2022covariance} develops covariance estimators under two different quantization protocols and demonstrates that in a broad range of cases the associated non-asymptotic estimation errors in spectral norm agree with those when full precision (i.e., unquantized) data is used, modulo logarithmic factors. This pioneering work has stimulated a series of follow-up studies that build on and further develop the approach in \cite{Dirksen2022covariance}, particularly the insight that pairs of one-bit samples obtained 
through dithered quantization give rise to a simple unbiased estimator of the covariance matrix provided that the underlying random vector has bounded coordinates. The paper \cite{Chen2023high} adopts this approach for sparse linear regression and matrix completion. The results on the former problem are refined and generalized to heavy-tailed covariate distributions in \cite{Chen2023heavy}. Improvements on the dithering-based covariance estimator in \cite{Dirksen2022covariance} are developed in \cite{Chen2025cov, Dirksen2024}. Extensions to various other estimation problems are considered in \cite{Chen2023LR, Hou2025, Xu2024bit, Dirksen2025subspace, Abdalla2026robust}.         
\vskip1.5ex
\noindent {\bfseries Contributions}. In this paper, we build on the approach in \cite{Dirksen2022covariance} to develop a linear regression estimator that operates
on 1-bit quantized predictor variables, their squares, and responses. The proposed 
quantization protocol roughly aligns with the dithering-based scheme in \cite{Dirksen2022covariance}, with the difference that the approach considered herein 
quantizes the squared predictors separately rather than generating quantized pairs per predictor. This modification is shown to achieve a reduction in variance when estimating 
the diagonal entries of the covariance matrix of the predictors, and requires corresponding changes in
the statistical analysis compared with \cite{Dirksen2022covariance, Chen2023high}. We derive 
a non-asymptotic bound on the $\ell_2$-estimation error of a quadratic programming-based estimator that is obtained by substituting the covariance matrix of the predictors and the cross-covariance of predictors and response by moment-based estimates. We also obtain
central limit theorems for this estimator under fixed and random designs with independent 
sub-Gaussian rows when the number of predictors is fixed. By analyzing the expression of 
the resulting asymptotic covariance matrix, we shed light on the asymptotic (MSE)-efficiency
relative to the ordinary least squares estimator based on the unquantized (i.e., full precision) data. It is shown that the leading term in this comparison is unavoidable since it arises from the quantization protocol rather than the specific estimator under consideration.  Sketching followed by quantization is identified as a specific random design within our analysis framework. Finally, we consider an extension of the proposed approach that includes an $\ell_1$-penalty to address scenarios with high-dimensional predictors under sparsity, and derive bounds on the estimation error (in the $\ell_1$- and $\ell_2$-norms). The use of the debiasing method in \cite{Zhang2014, Javanmard2014, vandegeer2014} is proposed to perform asymptotic inference for low-dimensional parameters such as individual regression coefficients. Empirical studies complement and corroborate our analysis, and illustrate potential practical gains concerning data transmission.      
\vskip1.5ex
\noindent {\em Related Work}. In addition to the references around the work \cite{Dirksen2022covariance} cited above, the literature on quantization touching 
on aspects important to the present paper is voluminous, and thus a review cannot
be comprehensive. Quantization with dithering, sometimes referred to as ``stochastic
quantization" has been proposed for compressed sensing \cite{Xu2020}, information retrieval \cite{Jacques2015}, and distributed
optimization and learning \cite{Suresh2017, Konecny2018}, to name a few recent applications. Several recent papers study the use of quantization in regression contexts. Building on the method of compressive least squares \cite{Maillard2009} in which a random projection is applied to lower the dimension of the predictors, the paper \cite{LiLi2019} studies the generalization error when quantization is applied subsequently. Sparse linear regression based on quantized predictors and samples is the subject of \cite{Cerone2019}. Nonparametric curve estimation via smoothing splines based on quantized responses is studied in \cite{Li2024}. It is worth distinguishing a separate line 
of work concerned with model compression, in which the goal is to compress the model 
parameters as opposed to inferring parameters from compressed data. A recent study of model compression for linear models is provided in \cite{Saha2022}.   
\vskip3ex
\noindent {\bfseries Notation}. 

\begin{table}[h!!!]
{\footnotesize
\begin{center}
\begin{tabular}{|ll|ll|}
\hline & & &\\[-1.5ex]
$n$, $d$ & sample size, dimension   & $Q$ & quantizer \\[1ex]
$\Sigma, \Sigma_n$ & (population) covariance matrix    & $R_n, L_n$ & quantizer ranges  \\[1ex]
$\Sigma_{Xy}$  & population cross-covariance matrix    & $\wt{[\ldots]}$ & quantized data $[\ldots]$  \\[1ex]
$\M{X}$, $\M{y}$ & design matrix, response vector  & $\sigma$ & noise level  \\[1ex]
$\M{S} $ & Sketching matrix  & $\beta_*$  & population regression parameter \\[1ex]
$\nnorm{\cdot}_p$ & $p$-norm   & $\E[\cdot]$, $\p(\ldots)$ & expectation, probability \\[1ex]
$\nnorm{\cdot}_{\psi_2}$  & sub-Gaussian norm     & $\{ e_j \}_{j = 1}^{d}$ & canonical basis of $\R^d$  \\[1ex]
$\nnorm{\cdot}_{\text{op}}$  & operator norm   & $a \wedge b, a \vee b$  & short for $\min\{a,b\}$, $\max\{a,b\}$  \\[1ex]
$\lambda_{\max}(\cdot)$   & maximum eigenvalue   & $A^{\textsf{c}}$  & complement of event $A$  \\[1ex]
$\lambda_{\min}(\cdot)$  & minimum eigenvalue   & $\overset{\text{D}}{\rightarrow}$ & convergence in distribution \\[1ex]
$\nnorm{\cdot}_{\textsf{F}}$    &  Frobenius norm & $\text{diag}(M)$ & short for a diagonal matrix whose \\[1ex]
$\nnorm{\cdot}_{\infty}$ & entry-wise absolute maximum  & & diagonal  entries are the same as   \\[1ex]
$\mathbb{B}_p(r)$ &  $\ell_p$-ball around $0$ of radius $r$  & & those of some matrix $M$  \\
\hline
\end{tabular}
\end{center}}
\vspace*{-3ex}
\caption{Summary of notation used repeatedly in this paper.}\label{tab:notation}
\end{table} 

\noindent Let $X$ be a zero-mean random variable. For $p \in \{1,2\}$, we define $\nnorm{X}_{\psi_p} =  \inf\{t > 0: \, \E[\exp(|X|^p / t^p)] \leq 2\}$. If $\nnorm{X}_{\psi_1}$ and 
$\nnorm{X}_{\psi_2}$ are finite, we say that $X$ is {\em sub-exponential} and {\em sub-Gaussian}, respectively. Let now $X = (X_j)_{j = 1}^d$ be a zero-mean random vector of dimension $d$. We let 
$\nnorm{X}_{\psi_2} = \sup_{x \in \R^d: \, \nnorm{x}_2 = 1} \nnorm{\nscp{X}{x}}_{\psi_2}$ and call
$X$ a sub-Gaussian random vector if $\nnorm{X}_{\psi_2}$ is finite.

We use $O(\cdot)$, $\Omega(\cdot)$ etc.~for the usual Landau symbols, and $O_{\p}(\cdot)$ and $o_{\p}(\cdot)$ for the corresponding stochastic order symbols \cite[][$\S$2.2]{vanderVaart1998}. Furthermore, 
we write $\wt{O}(\cdot)$ and $\wt{\Omega}(\cdot)$ to hide logarithmic factors. We use $\lesssim$, $\gtrsim$ in alignment with $O(\cdot)$ and $\Omega()$, and write $\asymp$ if both directions $\lesssim$ and $\gtrsim$ hold true. We write 
$C, C', \overline{C}, \wt{C}, c$ etc.~to denote positive constants; their values may change from instance to instance. A summary of frequently used notation is provided in Table \ref{tab:notation}. 

\section{Preliminaries}
{\bfseries Setup}. We consider a linear regression setup defined in terms of $(X,Y)$-pairs $\{(X_i, Y_i) \}_{i = 1}^n$ satisfying the additive error linear regression model 
\begin{equation}\label{eq:linear_model}
Y_i = X_i^{\T} \beta_* + \sigma \eps_i, \quad 1 \leq i \leq n, 
\end{equation}
where the $\{ \eps_i \}_{i = 1}^n$ are i.i.d.~zero mean, unit variance random variables independent of the $\{ X_i \}_{i = 1}^n$, which are of dimension $d$. 
In a fixed design scenario, the $\{ X_i \}_{i = 1}^n$ are considered fixed, and $\Sigma = \Sigma_n = \frac{1}{n} \su X_i X_i^{\T}$. In a random design scenario, the $\{ (X_i, Y_i) \}_{i = 1}^n$ are jointly random with 
$\Sigma = \E[X X^{\T}]$. To avoid cumbersome case distinctions, in the fixed design case we frequently identify $X$ with a random variable whose distribution consists of point masses $1/n$ at each of the $\{ X_i \}_{i = 1}^n$, so that operations such as $\E[\cdot]$ yield compatible results. Unless specified otherwise, we assume that 
$X$ is centered, i.e., $\E[X] = 0$ so that $\Sigma$ indeed constitutes a covariance matrix. Accordingly, there is no need for an intercept term in \eqref{eq:linear_model}. We note that the target parameter $\beta_*$ 
satisfies the relationship $\Sigma \beta_* =  \Sigma_{Xy}$ with $\Sigma_{Xy} = \E[XY]$. We assume that 
$\lambda_{\min}(\Sigma)$, i.e., the minimum eigenvalue of $\Sigma$, is bounded away from zero, hence $\beta_* = \Sigma^{-1} \Sigma_{Xy}$. 

Unlike standard regression, we assume that the $\{ (X_i, Y_i) \}_{i = 1}^n$ are not observed. Instead, 
we only observe triplets $\{ (\wt{X}_i, \wt{X_i^2}, \wt{Y}_i) \}_{i = 1}^n$ obtained through 1-bit dithered
quantization (as explained in detail in $\S$\ref{subsec:quantization} below) of $\{ (X_i, X_i^2, Y_i) \}_{i = 1}^n$. The collection and quantization of the squared $X$'s is performed to construct unbiased estimators 
of the covariance matrix. Quantization reduces storage and transmission requirements to three bits per $(X,Y)$-pair. In the sequel, we provide details on the specifics of the quantization protocol, the proposed estimation procedure, its statistical properties, and extensions including $\ell_1$-penalization for a ``large $d$" setting. 

We emphasize that in the setting under consideration, {\em both} $X$'s and $Y$'s are subject to (1-bit) quantization. This renders inference much more challenging than if only the $Y$'s were quantized, and indeed many approaches have been devised for the latter case \cite[e.g.,][]{BoufounosBaraniuk2008, Jacques2013a,  PlanVershynin2013a, Gopi2013, Baraniuk2017, Risuleo2019}. For example, if the distribution of the error terms is known, the 
resulting maximum likelihood estimator (MLE), is able to recover $\beta_* / \sigma$ \cite{Zymnis2009}. It is well-known that depending on the specific distribution assumed for the $\eps$'s, this MLE coincides with the MLE of popular binary regression models (probit, logistic, $\ldots$) \cite{Mcc1989}. 

If both $X$'s and $Y$'s are quantized and a specific distribution is assumed for both the $X$'s and the $\eps$'s, computing the MLE for general $d$ is not straightforward. Furthermore, it is restrictive to assume that the distribution of the $X$'s is known. Therefore, we propose a different route hinging on (unbiased) estimators of $\Sigma$ and $\Sigma_{Xy}$. Given those, $\beta_*$ can be estimated simply by plugging in those estimators. 

\subsection{Quantization}\label{subsec:quantization}
\begin{figure}
\centering
\includegraphics[height = 0.15\textwidth]{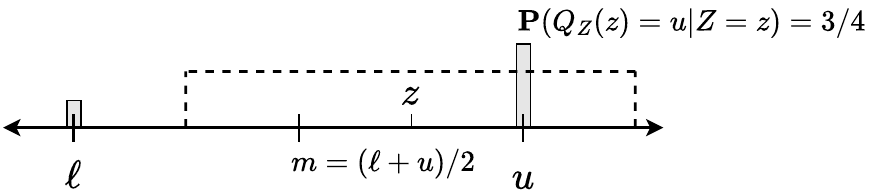}
\caption{Illustration of the quantization method for a random variable $Z$ supported on an interval $[\ell, u]$. Conditional on $\{Z = z = \frac{1}{4} \ell + \frac{3}{4} u \}$, the quantizer outputs $u$ with probability $3/4$ and $\ell$ with probability $1/4$ so that (conditional) unbiasedness holds, i.e., $\E[Q_Z(z)| Z = z] = z$.}\label{fig:quantization_dithering}
\end{figure}

This subsection is dedicated to the specific (scalar) quantization protocol under consideration, where ``scalar" means that all components of $X$'s are quantized separately. The approach taken is that of quantization with dithering, which randomly perturbs the quantizer input so that the quantizer output
is unbiased in the sense that in expectation (w.r.t.~the random perturbation), the latter equals the input. Let
$Z$ be a bounded random variable whose range is contained in $[\ell, u]$. Then its quantized version is 
obtained as $\wt{Z} = Q_Z(Z)$, where conditional on the event $\{Z = z\}$, $z \in [\ell, u]$, $Q_Z$ is defined as the following $\{\ell,u\}$-valued randomized map.
\begin{equation*}
Q_Z(z) = \begin{cases}
         u \; &\text{with probability} \;\, \dfrac{z - \ell}{u - \ell}, \\[2ex]
         \ell \; &\text{with probability} \;\, \dfrac{u - z}{u - \ell}. 
         \end{cases}
\end{equation*}
Let $\Delta = u - \ell$. Then equivalently, $Q_Z(z) \overset{\text{D}} = q(z + \xi)$, $z \in [\ell, u]$ with $\overset{\text{D}}{=}$ denoting equality in distribution, where $\xi \sim \mathrm{Unif}(-\frac{1}{2} \Delta, \frac{1}{2} \Delta)$ independent of $Z$ and $q(t) = \argmin_{v \in \{\ell, u \}}  |t - v|$, $t \in [\ell, u]$. Figure 
\ref{fig:quantization_dithering} illustrates this construction.

We note that $Q_Z$ depends on $Z$ only through its range, but we carry the random variable in the subscript to indicate the random variables that various quantizers act on. In the sequel, we shall consider quantizers $Q_{X}$, $Q_{X^2}$, and $Q_Y$. Even though the components of $X$ are quantized
separately and independently, we refrain from using varying quantizers $\{ Q_{X_j} \}_{j = 1}^d$; instead, we assume for simplicity that the ranges of the components are all equal to $[-R, R]$, permitting the 
use of a single quantizer (this assumption can be relaxed -- it suffices to assume lower
and upper bounds on component-specific interval boundaries $\{\ell_j\}_{j = 1}^d$ and $\{u_j\}_{j = 1}^d$). Similarly, $Q_{X^2}$ quantizes the squares $(X_j^2)_{j = 1}^d$ taking values in $[0, R^2]$. Finally, $Q_Y$ quantizes $Y$ whose range is assumed to be $[-L, L]$ with $L$ depending 
on the distributions of $X$ and $\eps$ (in particular $\sigma$) and $\beta_*$. 

\subsection{Approach}
Given quantized data $\wt{X}_i = \Big( Q_X(X_{ij}) \Big)_{j = 1}^d$, $\wt{X_i^2} = \Big( Q_{X^2}(X_{ij}^2) \Big)_{j = 1}^d$ and $\wt{Y}_i = Q_Y(Y_i)$, $1 \leq i \leq n$, we obtain 
\begin{equation}\label{eq:plugin_estimators}
\wh{\Sigma} \coloneq \frac{1}{n} \sum_{i = 1}^n \Big\{ \wt{X}_i \wt{X}_i^{\T} + \text{diag}(\wt{X_{i1}^2} - \underbrace{\wt{X}_{i1}^2}_{=R^2}, \ldots, \wt{X_{id}^2} - \underbrace{\wt{X}_{id}^2}_{=R^2}) \Big\}, \quad 
\wh{\Sigma}_{Xy} \coloneq \frac{1}{n} \su \wt{X}_i \wt{Y}_i. 
\end{equation}
We consider the estimator 
\begin{equation}\label{eq:estimator}
\wh{\beta} = \argmin_{\beta \in \R^d} \left\{ \frac{1}{2} \beta^{\T} \wh{\Sigma} \beta - \beta^{\T} \wh{\Sigma}_{Xy} \right\}.
\end{equation}
Note that the existence of a unique minimizer $\wh{\beta}$ requires that $\wh{\Sigma}$ is positive definite, which cannot be guaranteed for any (small) value of $n$, but can be established for $n = \wt{\Omega}(d)$ large enough. Existence provided, $\wh{\beta}$ is a solution of the estimating equation 
\begin{equation}\label{eq:est_eqn}
\Psi(\beta) = 0, \quad \Psi(\beta) \coloneq \wh{\Sigma} \beta - \wh{\Sigma}_{Xy}. 
\end{equation}
This estimating equation is {\em unbiased} in the sense that $\E[\Psi(\beta_*)] = 0$, where 
$\E[\cdot]$ is taken both w.r.t.~the randomness in $Y$ (and possibly $X$) and the randomness
at the quantization stage. This follows directly from the following observations:  
\begin{equation*}
\E[\wh{\Sigma}|\{ X_i \}_{i = 1}^n] = \frac{1}{n} \su X_i X_i^{\T}, \quad
\E[\wh{\Sigma}_{Xy}|\{ (X_i, Y_i) \}_{i = 1}^n] = \frac{1}{n} \su X_i Y_i, 
\end{equation*}
where the above expectations are w.r.t.~the quantization only. These relationships are seen to hold since for $1 \leq i \leq n$, we have
\begin{align*}
&\E[\wt{X}_{ij} \wt{X}_{ik} | X_{ij}, X_{ik}] = \E[\wt{X}_{ij}  | X_{ij}] \E[\wt{X}_{ik} | X_{ik}]
= X_{ij} X_{ik}, \;\;\;\, j,k=1,\ldots,d, \; j \neq k, \\
&\E[\wt{X_{ij}^2} | X_{ij}] = X_{ij}^2, \qquad \E[\wt{X}_{ij} \wt{Y}_i | X_{ij}, Y_i] = \E[\wt{X}_{ij} | X_{ij}] \E[\wt{Y}_i | Y_i] = X_{ij} Y_i, \;\;\;\, j=1,\ldots,d.
\end{align*}
Note that the diagonal modification of $\wh{\Sigma}$ in \eqref{eq:plugin_estimators} is necessary 
since $\p(\wt{X}_{ij}^2 = R^2) = 1$ by construction, hence $\wt{X}_{ij}^2$ in expectation 
(w.r.t.~$Q_X$) in general is not equal to $X_{ij}^2$, $1 \leq j \leq d$, $1 \leq i \leq n$. Note that another way of obtaining an unbiased estimator of the squares is to obtain a second 
independent quantized sets of the $X$'s \cite{Dirksen2022covariance, Chen2023high}, say, $\{ \wt{\wt{X}}_{ij} \}$, so that 
$\E[\wt{X}_{ij} \wt{\wt{X}}_{ij} | X_{ij}] = X_{ij}^2$, $1 \leq j \leq d$, $1 \leq i \leq n$. A shortcoming of the latter approach is a potentially serious increase in variance: it is easy to verify that 
\begin{equation*}
\var(\wt{X^2_{ij}} | X_{ij}) = R^2 X_{ij}^2 - X_{ij}^4, \qquad \var(\wt{X}_{ij} \wt{\wt{X}}_{ij} | X_{ij}) = R^4 - X_{ij}^4, \;\;\;\, j=1,\ldots,d, \; i=1,\ldots,n. 
\end{equation*}
We note that the former variance is never larger than the latter, and can be much smaller on
average if the distribution of the $\{ X_{ij} \}$ does not carry much mass at the boundary 
points $\{-R,R\}$ of its support. We discuss implications for the efficiency for estimating
the regression parameter in Appendix \ref{app:ARE_diagonalSigma}: we demonstrate that for a single Gaussian predictor,
the asymptotic relative efficiency is improved by a factor of $1.5$.

\subsection{Extension to random variables without fixed range}\label{subsec:nofixed_range}
The approach outlined so far requires that the data to be quantized have bounded range. We consider an extension to settings in which the quantizer input consists of random variables with sub-Gaussian tails. In this case, the framework remains applicable by having the quantization
ranges depend on the number of samples $n$, i.e., we vary $R = R_n$ and $L = L_n$ with $n$. Specifically, we shall consider the following situation: 
\begin{itemize}[leftmargin=10ex]
\item[(\textsf{SG}-$X$):] The $X_i = (X_{ij})_{j = 1}^d$, $1 \leq i \leq n$, are i.i.d.~zero-mean, unit variance sub-Gaussian random vectors with covariance matrix $\Sigma$.
\item[(\textsf{SG}-$\eps$):] The $\eps_i$, $1 \leq i \leq n$, are i.i.d.~zero-mean, unit variance sub-Gaussian random variables.  
\end{itemize}
Under these two assumptions made, the following can be shown.  
\begin{prop}\label{prop:quantizer_bound} Let {\em (\textsf{SG}-$X$)} and {\em (\textsf{SG}-$\eps$)} be satisfied. Let 
$K = \max_{1 \leq j \leq d} \nnorm{X_{1j}}_{\psi_2}$, $\overline{K} = \nnorm{\Sigma^{-1/2} X_1}_{\psi_2}$, and $K_{\eps} = \nnorm{\eps_1}_{\psi_2}$. Consider
the events 
\begin{equation*}
\mc{R} = \Bigg\{ \max_{\substack{1 \leq i \leq n \\ 1 \leq j \leq d}} |X_{ij}| \leq 2C_{K} \sqrt{\log(n \cdot d)} \Bigg\}, \quad \mc{L} = \big\{ \max_{1 \leq i \leq n} |Y_i| \leq 2(C_{\overline{K}} \nnorm{\Sigma^{1/2} \beta_*}_2 + \sigma K_{\eps}) \sqrt{\log n} \big \}. 
\end{equation*}
Then $\p(\mc{R}) \geq 1 - 2/(n \cdot d)$ and $\p(\mc{L}) \geq 1 - 4/n$.  
\end{prop}
\noindent For now, let us assume that $d \leq n$. In this case, according to Proposition \ref{prop:quantizer_bound}, we may choose 
$R_n \asymp \sqrt{\log n}$ and $L_n \asymp (\nnorm{\Sigma^{1/2} \beta_*}_2 + \sigma) \sqrt{\log n}$ so that $[-R_n, R_n]$ and $[-L_n, L_n]$ contain the $\{ X_{ij} \}$ and $\{ Y_i \}$, respectively, with high probability. Conditional on that event, quantization based on these ranges does not introduce any additional distortion compared to the fixed range case. Nevertheless, it is worth noting that $\E[\wh{\Sigma}|\mc{R} \cap \mc{L}] \neq \E[\wh{\Sigma}]$ and $\E[\wh{\Sigma}_{Xy}|\mc{R} \cap \mc{L}] \neq \E[\wh{\Sigma
}_{Xy}]$, where $\wh{\Sigma}$ and $\wh{\Sigma}_{Xy}$ are as in \eqref{eq:plugin_estimators}. The following statement shows that by choosing the constant in front of $\sqrt{\log n}$ in the quantization thresholds large enough, the discrepancy between the conditional and unconditional quantities becomes negligible.

\begin{prop}\label{prop:bias_control_truncation}
Let $q > 0$ and suppose that $n > 8^{1/q}$. Choose $R_n = \sqrt{2 C_K (q+1) \log (n \cdot d)}$ and $L_n =  \big(\sqrt{ 2(q+1) C_{\overline{K}}} \nnorm{\Sigma^{1/2} \beta_*}_2 + \sqrt{2 (q+1) C_{K_{\eps}}} \big) \sqrt{\log n}$, and consider the events $\mc{R} = \{\max_{i,j} |X_{ij}| \leq R_n \}$ and 
$\mc{L} = \{\max_i |Y_i| \leq L_n \}$. Then:
\begin{align*}
&\nnorm{\E[\wh{\Sigma} | \mc{R} \cap \mc{L}] - \Sigma}_{\text{{\em op}}} \leq d \, C \max\{ \nnorm{\Sigma^{1/2} \beta_*}_2 + \sigma, 1 \} \, n^{-q/2} \\  
&\nnorm{\E[\wh{\Sigma}_{Xy} | \mc{R} \cap \mc{L}] - \Sigma_{Xy}}_{\infty} \leq  C' \max\{ \nnorm{\Sigma^{1/2} \beta_*}_2 + \sigma, 1 \} \, n^{-q/2},
\end{align*}
for constants $C, C' > 0$ depending only on the sub-Gaussian norms $K$, $\overline{K}$, and $K_{\eps}$. 
\end{prop}
\noindent As will become clear in $\S$\ref{sec:main_result}, the sampling variation of 
$\wh{\Sigma}$ and $\wh{\Sigma}_{Xy}$ are $O(\sqrt{d/n})$
and $O(n^{-1/2})$ in operator and sup-norm, respectively, thus choosing $q \geq 3$ in Proposition \ref{prop:bias_control_truncation} yields bias terms that are negligible. 

\subsection{Sketching followed by Quantization}\label{subsec:Sketch+Quant}
Sketching is another important method of data compression \cite{Vempala2005, Achlioptas2003, Sarlos2006}. In sketched linear 
regression \cite{Mahoney2011, DrineasMahoney2016, PilanciWainwright2015, Yang2017}, a random matrix $\M{S}$ of dimension $m$-by-$n$ with (appropriately scaled) isotropic rows, $m \ll n$,  is applied to both $\M{X}$ and $\M{y}$, where
$\M{X}$ denotes the $n$-by-$d$ matrix whose rows are given by the $\{ X_i^{\T} \}_{i=1}^n$, and $\M{y} = (Y_i)_{i = 1}^n$. This yield the modified least squares problem 
\begin{equation*}
\min_{\beta} \, \left\{ \frac{1}{2n} \beta^{\T} \M{X}^{\T} \M{S}^{\T} \M{S} \M{X} \beta - \beta^{\T} \M{X}^{\T} \M{S}^{\T} \M{S} \M{y} \right\}
\end{equation*}
Since $\M{S}$ is isotropic, the expected value of the least squares objective w.r.t.~$\M{S}$ for any fixed $\beta$ equals that of the original least squares objective, which loosely justifies the use of this approach. Common choices for 
$\M{S}$ include i.i.d.~Gaussian entries, i.e., $S_{ij} \sim N(0,1/m)$, $1 \leq i \leq n$, $1 \leq j \leq m$, or the ternary distribution on $\{-1, 0, 1\}$ with probabilities $1/6, 2/3, 1/6$, respectively, scaled by $\sqrt{3/m}$ \cite{Achlioptas2003}. 

It has been observed that the compression achieved by sketching can be 
enhanced by subsequent quantization \cite{SlawskiMitzenmacherLi2016, SlawskiLi2017, Charikar2002, Rane2013, Jacques2015, Li2014}. This approach can be adopted here via 
a slight modification of the sketching step, assuming that the rows of $\M{S}$ 
are sub-Gaussian.

\begin{figure}
\begin{center}
\includegraphics[height = 0.1\textheight]{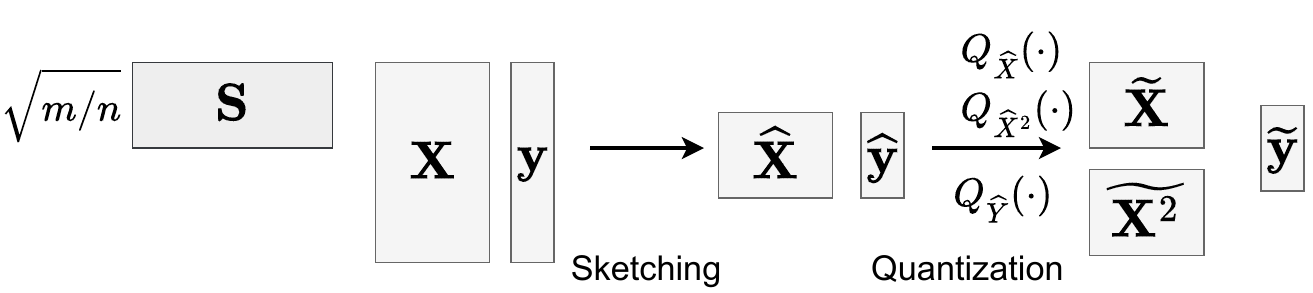}
\end{center}
\vspace*{-2.5ex}
\caption{Schematic summary of the sketching followed by quantization pipeline described in the text. In the first step, the $X$'s and $Y$'s are sketched via a rescaled sketching matrix $\M{S}$. In the second step, quantization is performed on the sketched data.}\label{fig:sketching_plus_quantization}
\end{figure}

\begin{itemize}
\item[1.] In the first step, we compute $\wh{X}_{ij} = (n/m)^{-1/2} \sum_{k = 1}^n S_{ik} X_{kj}$, $1 \leq i \leq m$, $1 \leq j \leq d$. Note that conditional on the $\{ X_{kj }\}_{k = 1}^n$, the $\wh{X}_{ij}$ is a zero-mean sub-Gaussian random variable with 
$\nnorm{\wh{X}_{ij}}_{\psi_2} \lesssim n^{-1/2} \left(\sum_{k = 1}^n X_{kj}^2 \right)^{1/2} = O(1)$, $1 \leq i \leq m$, $1 \leq j \leq d$. Here, the $O(1)$ is justified in scenarios in which the $X$'s are bounded or have been subject to centering and unit standard deviation scaling. The $O(1)$ also holds with high probability for random $X$ under suitable tail conditions. 
Similarly, we obtain $\wh{Y}_i = (n/m)^{-1/2} \sum_{k = 1}^n S_{ik} Y_k$
with $\nnorm{\wh{Y}_i}_{\psi_2} \lesssim n^{-1/2} \left( \sum_{k = 1}^n Y_k^2 \right)^{1/2} = O(1)$, $1 \leq i \leq m$, if the $Y$'s have fixed range or the $X$s and $\eps$s are sub-Gaussian. 
\item[2.] The $\{ \wh{X}_{ij} \}$ and $\{ \wh{Y}_{i} \}$ are quantized according to the protocol in  $\S$\ref{subsec:quantization} (if the entries of $\M{S}$, the $X$'s and the $Y$'s have fixed range) and $\S$\ref{subsec:nofixed_range}, respectively, yielding $\{ \wt{X}_{ij} \}$, $\{ \wt{X_{ij}^2} \}$ and $\{ \wt{Y}_i \}$. 
\end{itemize}
A diagram summarizing these two steps is provided in Figure \ref{fig:sketching_plus_quantization}. 

We observe that conditional on the $\{ X_i \}_{i = 1}^n$, the $\{ \wh{X}_i = (\wh{X}_{ij})_{j= 1}^d \}_{i = 1}^n$ are zero-mean sub-Gaussian random variables with covariance matrix $\Sigma = \M{X}^{\T} \M{X} / n$. Furthermore, 
\begin{equation*}
\wh{Y}_i = (n/m)^{-1/2} \sum_{k = 1}^n S_{ik} (X_k^{\T} \beta_* + \sigma \eps_k)
= \wh{X}_{i}^{\T} \beta_* + \wh{\eps}_i, \quad \wh{\eps}_i = \sigma (n/m)^{-1/2} \sum_{k = 1}^n S_{ik} \eps_k, \;\;\; 1 \leq i \leq m,
\end{equation*}
so that under the event $\{ \nnorm{\eps}_2 \leq C \sqrt{n} \}$ (which holds with high probability assuming \linebreak (\textsf{SG}-$\eps$)), the 
$\{ \wh{\eps}_i \}_{i = 1}^n$ are i.i.d.~zero mean sub-Gaussian random variables. Consequently, the resulting estimator $\wh{\beta}$ based on $\{ \wt{X}_{ij}, \wt{X_{ij}^2}$ and $\{ \wt{Y}_i \}$ obtained from $\{ \wh{X}_{ij} \}$ and 
$\{ \wh{Y}_i \}$, respectively, exhibits the same properties (cf.~the upcoming Section $\S$\ref{sec:main_result}) that hold under plain quantization in the case
of random design subject to condition (\textsf{SG}-$X$) and $n = m$. In view of this reduction and for the sake of brevity, we refrain from providing a separate analysis of the sketching-followed-by-quantization pipeline.

\section{Main results}\label{sec:main_result}
In this section, we analyze the estimator $\wh{\beta}$ proposed in \eqref{eq:estimator}. Our analysis consists of two parts: we first state asymptotic properties in the regime 
$n \rightarrow \infty$ and $d$ fixed, before stating a non-asymptotic bound on the 
$\ell_2$-estimation error in which $d$ may depend on $n$. We recover the usual
$\sqrt{d/n}$ error rate, apart from a log factor. Our analysis covers both fixed and random $X$. We also analyze the relative efficiency of our estimator compared with the (ordinary) least squares estimator applied to the full-precision data. That discussion is complemented with
a lower bound on the MSE of the MLE based on quantized data, which exhibits the same 
qualitative dependency on the quantizer ranges. We conclude that substantial improvements
concerning that dependency are unlikely, as long as the quantization protocol remains
unchanged. 

At the end of this section, we consider a high-dimensional setup for which we study
an $\ell_1$-penalized variant of the estimator \eqref{eq:estimator}. We present 
a bound on the estimation error in line with the sparse estimation literature, and 
propose a path towards inference for the individual regression coefficients via 
the debiasing approach in \cite{Zhang2014, Javanmard2014, vandegeer2014}. 

\subsection{Asymptotic analysis}

We start by showing asymptotic Normality of the proposed estimator. Note that the 
statement is sub-divided according to scenarios (I)--(III): random $X$ with fixed 
quantizer range, random $X$ with varying quantizer range, and fixed $X$. For the second
case, we implicitly adopt a triangular array setup with $\{ (\wt{X}_i^{(n)}, \wt{Y}_i^{(n)}) \}_{i = 1}^n$ since the distribution of the quantized data depends on $n$ via the quantizer
ranges $R_n$ and $L_n$ even if the underlying $(X,Y)$-pairs are i.i.d.; for ease of 
presentation, we refrain from re-iterating this technicality in the statement or proofs. 
\begin{theo}\label{theo:asymptotic}
Consider the linear model \eqref{eq:linear_model} and the estimator $\wh{\beta}$ defined in \eqref{eq:estimator} with $\wh{\Sigma}$ and $\wh{\Sigma}_{Xy}$ given by \eqref{eq:plugin_estimators}. Then:\\
(I) If $X$ is random with fixed range $[-R,R]$ and $Y$ has fixed range $[-L,L]$, respectively, we have 
\begin{equation*}
\sqrt{n} (\wh{\beta} - \beta_*) \overset{\text{\em D}}{\rightarrow} N_d(0, \Sigma^{-1}  \Gamma \Sigma^{-1}),
\end{equation*}
as $n \rightarrow \infty$ with $\Gamma = L^2 \{ \Sigma - \text{{\em diag}}(\Sigma) + R^2 I_d \} + \Gamma_a - \Gamma_b$, where
\begin{align*}
&\Gamma_a  = \E[(\wt{X}_1 \wt{X}_1^{\T} + \Delta_1) \beta_* \beta_*^{\T} (\wt{X}_1 \wt{X}_1^{\T} + \Delta_1)], \\
&\Gamma_b = \E[(\wt{X}_1 \wt{X}_1^{\T} + \Delta_1) \beta_* \wt{X}_1^{\T} \wt{Y}] + \E[\wt{X}_1 \wt{Y}  \beta_*^{\T} (\wt{X}_1 \wt{X}_1^{\T} + \Delta_1)],
\end{align*}
and $\Delta_i = \text{{\em diag}}(\wt{X_{i1}^2} - R^2, \ldots, \wt{X_{id}^2} - R^2)$, $1 \leq i \leq n$. Moreover,
$\nnorm{\Gamma_a - \Gamma_b}_{\infty} \leq 3 R^4 \nnorm{\beta_*}_1^2$.\\ 
\noindent (II) In the unbounded sub-Gaussian case, i.e., under conditions {\em (\textsf{SG}-$X$)} and {\em (\textsf{SG}-$\eps$)}, using the choices $R = R_n$ and $L = L_n$ 
according to Proposition \ref{prop:bias_control_truncation} with $q = 2$, for any 
sequence $\tau_n \asymp \log n$, we have 
\begin{equation*}
\tau_n^{-1} \sqrt{n}(\wh{\beta} - \beta_*) \overset{\text{\em D}}{\rightarrow} N_d(0, \Sigma^{1-}  \Gamma \Sigma^{-1}),
\end{equation*}
as $n \rightarrow \infty$, where $\Gamma = \lim_{n \rightarrow \infty} \{ \tau_n^{-2} [L_n^2 \{ \Sigma - \text{{\em diag}}(\Sigma) + R_n^2 I_d \} + \Gamma_{a,n} - \Gamma_{b,n}] \}$ is a finite,
symmetric positive definite matrix, and $\Gamma_{a,n}$, $\Gamma_{b,n}$ are defined as $\Gamma_a$ and 
$\Gamma_b$ above (but now depend on $n$).\\
\noindent (III) If the $X$'s are fixed elements of $\R^d$ such that $\lim_{n \rightarrow \infty} \frac{1}{n} \su X_i X_i^{\T} = \overline{\Sigma}$ for some finite, symmetric positive definite matrix 
$\overline{\Sigma}$, we have 
\begin{equation*}
\mathbb{V}_n^{-1/2} \overline{\Sigma} (\wh{\beta} - \beta_*) \overset{\text{\em D}}{\rightarrow} N_d(0, I_d)
\end{equation*}
as $n \rightarrow \infty$, where $\mathbb{V}_n = n^{-2} \su \cov(\wt{X}_i \wt{Y}_i - (\wt{X}_i \wt{X}_i^{\T} + \Delta_i) \beta_*) = O(n^{-1})$. 
\end{theo}
\noindent In informal terms, Theorem \ref{theo:asymptotic} states that for each of the cases (I)--(III)
\begin{equation}\label{eq:asymptotic_normal_informal}
\wh{\beta} \overset{\text{D}}{\approx} N\left(\beta_*, \Sigma^{-1} \cov\left( \frac{1}{n} \sum_{i = 1}^n \wt{X}_i \wt{Y}_i - (\wt{X}_i \wt{X}_i^{\T} + \Delta_i) \beta_* \right)  \Sigma^{-1} \right), 
\end{equation}
where $\overset{\text{D}}{\approx}$ means ``approximately distributed as" for large $n$. The form of covariance matrix is as would be anticipated from the characterization of $\wh{\beta}$ as root of \eqref{eq:est_eqn} and the theory of estimating equations (\cite[][$\S$5.3]{vanderVaart1998} and \cite{Stefanski2002}), noting that 
$\wh{\Sigma}$ equals the Jacobian of $\beta \mapsto \Psi(\beta)$ in \eqref{eq:est_eqn}. It is possible but rather tedious to evaluate the covariance matrix of the term inside $\cov(\ldots)$ in \eqref{eq:asymptotic_normal_informal}, which generally depends on the estimand $\beta_*$. It is more practical to use the consistent estimator 
\begin{equation}\label{eq:est_standard_error}
\left( \frac{1}{n} \su \{ \wt{X}_i \wt{Y}_i - (\wt{X}_i \wt{X}_i^{\T} + \Delta_i) \wh{\beta} \} \right) \Big / n, 
\end{equation}
which can then be used for asymptotic inference such as the construction of confidence intervals.  
\vskip1ex
\noindent {\em Asymptotic relative efficiency} (ARE). Besides asymptotic inference, the second important use of Theorem \ref{theo:asymptotic} is the evaluation of the ARE with respect to 
the ordinary least squares estimator $\wh{\beta}^{o}$  having access to the uncompressed data $\{ (X_i, Y_i) \}_{i = 1}^n$. Specifically, we consider the ARE w.r.t.~mean squared error (MSE). It is well-known that under \eqref{eq:linear_model}  
\begin{align}\label{eq:olsbound}
\text{MSE}(\wh{\beta}^{o}) = \E[\nnorm{\wh{\beta}^o - \beta_*}_2^2] = \tr(\cov(\wh{\beta}^o)) &= \frac{\sigma^2}{n}  \E \left[ \tr \left( \left\{ \frac{1}{n} \su X_i X_i^{\T} \right \}^{-1} \right) \right] \notag \\
&=  \frac{\sigma^2}{n} \tr(\Sigma^{-1}) + o(n^{-1})  \notag \\
\text{where}& \;\, 
\frac{\sigma^2}{\lambda_{\max}(\Sigma)} \frac{d}{n} \leq \frac{\sigma^2}{n} \tr(\Sigma^{-1}) \leq \frac{\sigma^2}{\lambda_{\min}(\Sigma)} \frac{d}{n},
\end{align}
as $n \rightarrow \infty$, for any sub-Gaussian $X$ as under consideration herein. According to Theorem \ref{theo:asymptotic} for random $X$, we have 
\begin{equation*}
\E[\nnorm{\wh{\beta} - \beta_*}_2^2] = \frac{1}{n} \tr(\Sigma^{-1}  [L_n^2 \{ \Sigma - \text{{diag}}(\Sigma) + R_n^2 I_d \} + \Gamma_{a,n} - \Gamma_{b,n}]  \Sigma^{-1}) + o(\tau_n^2/n) 
\end{equation*}
To make matters simple, suppose for the remainder of this subsection that $\Sigma = I_d$. Hence
\begin{align*}
\E[\nnorm{\wh{\beta} - \beta_*}_2^2] = \frac{d}{n} R_n^2 L_n^2 + \frac{1}{n} \tr(\Gamma_{a,n} - \Gamma_{b,n}) + o(\tau_n^2/n) \leq \frac{d}{n} \{ R_n^2 L_n^2 + 3 R_n^4 \nnorm{\beta_*}_1^2 \} + o(\tau_n^2/n),
\end{align*}
where the inequality uses the bound on $\nnorm{\Gamma_{a,n} - \Gamma_{b,n}}_{\infty}$ in Theorem \ref{theo:asymptotic}, which we hypothesize to be loose since it uses the bound $|X^{\T} \beta_*|^2 \leq R_n^2 \nnorm{\beta_*}_1^2$ effectively twice (cf.~Appendix \ref{app:theo:asymp}). However, the dependence on $R_n^2 L_n^2$ is likely inevitable in view of Theorem \ref{theo:lowerbound} below, which shows that for $d = 1$ even the variance of the MLE incurs this factor. The ARE of $\wh{\beta}$ relative to 
$\wh{\beta}^{o}$ is hence conjectured to be a multiple ($\geq 1$) of $R_n^2 L_n^2 / \sigma^2$. Observe that since 
$R_n^2 \geq \min_{1 \leq j \leq d} \E[X_j^2] = 1$, we have
\begin{equation}\label{eq:ARE_lower}
R_n^2 L_n^2 / \sigma^2  \geq L_n^2 / \sigma^2 \geq \frac{\nnorm{\beta_*}_2^2 + \sigma^2}{\sigma^2} \geq 1 + \frac{\nnorm{\beta_*}_2^2}{\sigma^2}.  
\end{equation}
assuming that $L_n \geq \nnorm{\beta_*}_2$, which is necessary for $[-L_n, L_n]$ to be a
valid quantizer range since $\E[(X^{\T} \beta_*)^2]^{1/2} = \nnorm{\beta_*}_2$. Eq.~\eqref{eq:ARE_lower} indicates that the ARE can be arbitrarily large depending on 
$\nnorm{\beta_*}_2^2 / \sigma^2$, which can be thought of as the signal-to-noise
ratio of the problem. 
Conversely, note that if the $X$'s are i.i.d.~Rademacher random variables, we have 
$R_n \equiv 1 = \E[X_j^2] = 1$, $1 \leq j \leq d$, i.e., and the lower bound \eqref{eq:ARE_lower} can be arbitrarily
close to $1$ depending on the SNR. Note that in this case, the $X$'s are already in 
1-bit form, i.e., quantization does not introduce any distortion. In general, 
$R_n^2 \geq 1$ represents the loss in efficiency introduced by quantization of the $X$'s 
alone -- e.g., if the $X$'s follow the uniform distribution on $[-3^{1/2}, 3^{1/2}]$
(scaled to unit variance), we have $R_n^2 \equiv 3$. 

Finally, if the $X$'s are unbounded, e.g., Gaussian random variables, then $\wh{\beta}$
 no longer enjoys the $O(1/n)$ rate for the MSE as $\wh{\beta}^{o}$, but the slower rate 
$O(\log^2(n) / n)$.

\subsection{Non-asymptotic analysis}
In this subsection, we provide a non-asymptotic bound on the $\ell_2$-estimation error for the regression parameter. As elaborated below, the stated bound yields similar findings 
as the asymptotic result presented in the preceding subsection. 

\begin{theo}\label{theo:nonasymptotic}
Consider the linear model \eqref{eq:linear_model} and the estimator $\wh{\beta}$ defined in \eqref{eq:estimator} with $\wh{\Sigma}$ and $\wh{\Sigma}_{Xy}$ given by \eqref{eq:plugin_estimators}. Then:\\
{\em (I)} In the fixed range case in which the $\{ |X_{ij}| \}$ and the $\{ |Y_i| \}$ are bounded 
by $R$ and $L$, respectively, if 
\begin{equation*}
n \geq  \left \{ \frac{512 d R^4 + 128 d R^2 \nnorm{\Sigma}_{\text{{\em op}}}}{\lambda_{\min}^2(\Sigma)} \vee  32 (2d+1)R^2 \right \} \{ \log(n) \vee \log(2d) \},
\end{equation*}
we have $\nnorm{\wh{\beta} - \beta_*}_2 \leq \text{{\em \textsf{err}}}(d,n,R,L)$ with probability at least $1 - 2/n$, where
{\small \begin{align*}
\text{{\em \textsf{err}}}(d,n,R,L) \coloneq& \frac{2}{\lambda_{\min}(\Sigma)} \Bigg( \nnorm{\beta_*}_2 \sqrt{\frac{\{ 32 d R^4 + 8dR^2 \nnorm{\Sigma}_{\text{\em op}} \}  \log(n) \vee \log(2d)}{n}} +  \\
&\nnorm{\beta_*}_2 \frac{8 (2d+1) R^2  \{ \log(n) \vee \log(2d) \}}{n} +  \sqrt{\frac{4d R^2 L^2  \{ \log(2d) \vee \log n \}}{n}} \Bigg) \Big.
\end{align*}}
\noindent {\em (II)} In the unbounded sub-Gaussian case, i.e., under conditions {\em (\textsf{SG}-$X$)} and {\em (\textsf{SG}-$\eps$)}, using the choices $R = R_n$ and $L = L_n$ 
according to Proposition \ref{prop:bias_control_truncation} with $q = 3$ and if additionally 
$n \geq C \{ \lambda_{\min}^{-1}(\Sigma) d \nnorm{\Sigma^{1/2} \beta_*}_2 + \sigma, 1 \}\}^{2/3}$, where $C > 0$ is a constant depending only on the sub-Gaussian norms $K, \overline{K}$ and $K_{\eps}$,  it holds that 
\begin{equation*}
\nnorm{\wh{\beta} - \beta_*}_2 \leq 2 \, \text{{\em \textsf{err}}}(d,n,R,L) + O\left(\frac{d}{n^{3/2}} \right)
\end{equation*}
with $\text{{\em \textsf{err}}}(d,n,R,L)$ as above, with probability at least $1 - 6/n - 2/(n\cdot d)$.  
\end{theo}
\noindent Theorem \ref{theo:nonasymptotic} yields a bound of the form $\nnorm{\wh{\beta} - \beta_*}_2 = \wt{O}(\sqrt{d/n})$ for $n \gtrsim d$ (modulo log factors) as expected. The additional (explicit) log factors in $n$ could be 
replaced by factors $\log(1/\delta)$ for any $\delta > 0$, for a failure probability 
of $1 - O(\delta)$. A logarithmic factor in $d$ is incurred from the use of the matrix
Bernstein inequality. Apart from numerical constants (which are unlikely to be optimal) and the pre-factor $1/\lambda_{\min}(\Sigma)$, which is 
expected (cf.~Eq.~\eqref{eq:olsbound}) and unrelated to quantization, the constants
preceding the $\sqrt{d/n}$ term are $\nnorm{\beta_*}_2 (R^2 + R) \nnorm{\Sigma}_{\text{op}}^{1/2} + R L$. Assuming that $\Sigma = I_d$, we need to have $L \geq \nnorm{\beta_*}_2$ as elaborated in the previous subsection. Depending on the setup, e.g., for Gaussian predictors, we have $L \geq R \nnorm{\beta_*}_2$. In this case, $RL$ becomes the leading ``constant" (where the quotes reflect that $R$ and $L$ actually scale as $\sqrt{\log n}$), in alignment with 
what is conjectured based on the asymptotic distribution of the estimator.

\subsection{Lower bound}\label{subsec:lowerbound}
According to the results stated in the previous sections, the MSE of the estimator $\wh{\beta}$ scales with $R^2 L^2$. It is worth investigating if the presence of this term can be avoided. In the sequel, we present a result that lower bounds the Cram\'{e}r-Rao bound for $d = 1$ and Gaussian $X$. In slightly simplified terms, this classical result represents a lower bound on the variance of asymptotically unbiased and Normal estimators \cite[e.g.,][Ch.~4.5]{Shao2003}. The lower bound presented below is proportional to $R^2 L^2$, indicating that this factor cannot be eliminated.    

We consider the following setting. 
\begin{align*}
X, \eps \overset{\text{i.i.d.}}\sim N(0, 1), \quad Y = X \beta_* + \sigma \eps, 
\end{align*}
Let $\{ (X_i, Y_i) \}_{i = 1}^n$ be i.i.d.~distributed as $(X,Y)$ in the above display
and let $\{ (\wt{X}_i, \wt{Y}_i) \}_{i = 1}^n$ be the quantized counterparts obtained 
according to the framework in $\S$\ref{subsec:quantization}, the maximum likelihood estimator (MLE) 
of $\beta_*$ for $\sigma$ known maximizes the likelihood function \cite{Li2016}
\begin{align}\label{eq:likelihood_quantdith}
\textsf{L}(\beta) = \prod_{i = 1}^n \{ \pi(\beta)^{c_i} (1-\pi(\beta))^{1-c_i} \},  \quad &c_i = I(\text{sign}(\wt{X}_i) = \text{sign}(\wt{Y}_i)),  \\[-2ex]
&\pi(\beta) = 2 \p(\wt{X}_1 > 0, \wt{Y}_1 > 0), \;\; 1 \leq i \leq n, \notag
\end{align}
where $\pi(\beta)$ is known as ``collision probability" \cite{Charikar2002}. It can be expressed as 
\begin{align}
 \pi(\beta) &= \frac{1}{2L}
  \frac{1}{2R} \int_{-L}^{L} \int_{-R}^R \int_{0}^{\infty}\int_{0}^{\infty}
  \frac{|V|^{-1/2}}{2\pi}  \exp\left(-\frac{1}{2}\begin{pmatrix}
      x - z_1 \\
      y - z_2 
      \end{pmatrix} V^{-1} \begin{pmatrix}
      x - z_1 \\
      y - z_2 
      \end{pmatrix} 
    \right) \; dx \, dy \, dz_1 \, dz_2, \notag \\[1ex]
 &\quad V \coloneq \begin{pmatrix}
1 & \beta_* \\
\beta_* &\beta_*^2 + \sigma^2. \label{eq:collisionprob}
 \end{pmatrix}
\end{align}
where we have used that quantization with dithering of $X$ and $Y$ is equivalent to
adding noise uniformly distributed on $[-R,R]$ and $[-L,L]$,
respectively. 

The following theorem provides an upper bound on the Fisher information associated
with the likelihood function \eqref{eq:likelihood_quantdith}. 
\begin{theo}\label{theo:lowerbound}
Consider the likelihood in \eqref{eq:likelihood_quantdith} based on $L = L_n \geq \sqrt{2(\sigma^2 + \beta_*^2) (\log \log n)^{1/2}}$ and $R = R_n = O(\sqrt{\log n})$. The following holds for the associated Fisher information: 
{\em \begin{equation*}
\frac{1}{n} \E\left[-\frac{d^2}{d \beta^2} \log \textsf{L}(\beta) \dev{\beta}{\beta_*} \right] = \{ \dot{\pi}(\beta_*) \}^2 \left(\frac{1}{\pi(\beta_*)} + \frac{1}{1 - \pi(\beta_*)}\right) \leq \frac{1}{L^2 \cdot R^2} C(\beta_*, \sigma),
\end{equation*}}
where $C(\beta_*, \sigma)$ is a constant depending only on $\beta_*$ and $\sigma$.  
\end{theo}
\noindent Note that Theorem \ref{theo:lowerbound} is valid for a fairly wide range of quantizer
ranges. The lower bound on $L$ involves the fourth root of an iterated logarithm, which behaves like a constant for all practical purposes.

\subsection{Lasso}
In this subsection, we consider high-dimensional predictors, i.e., $d$ is potentially
much larger than $n$. At the same time, we suppose that the target $\beta_*$ is sparse in the sense that it has few non-zero coordinates. Sparsity is incorporated by employing an $\ell_1$-penalty, also known as the Lasso \cite{Tib1996}. Specifically, we consider the following counterpart to \eqref{eq:estimator}:
\begin{equation}\label{eq:estimator_ell1}
\wh{\beta} = \argmin_{\beta \in \mathbb{B}_1^d(\sqrt{s}B)} \left\{ \frac{1}{2} \beta^{\T} \wh{\Sigma} \beta - \beta^{\T} \wh{\Sigma}_{Xy} + \lambda \nnorm{\beta}_1 \right \},
\end{equation}
where $\mathbb{B}_1^d(r) = \{\beta \in \R^d: \nnorm{\beta}_1 \leq r \}$ denotes the 
$\ell_1$-ball in $\R^d$ with radius $r$. The constraint set in the optimization problem \eqref{eq:estimator_ell1} reflects the assumptions that 
$\nnorm{\beta_*}_2 \leq B$ for some $B > 0$, and that $\beta_*$ has at most $s$ non-zero coordinates. The constraint is imposed since if $d/n$ is large, $\wh{\Sigma}$ may have a negative eigenvalue, in which case \eqref{eq:estimator_ell1} may not exist if the constraint is dropped. With the addition of the constraint, we follow prior literature on $\ell_1$-penalized estimation with non-convex objective (cf., e.g., \cite{Loh2012}).     

\begin{theo}\label{theo:lasso}
Consider the linear model \eqref{eq:linear_model} such that $\beta_*$  has
at most $s$ non-zero coefficients and $\nnorm{\beta_*}_2 \leq B$. 
Suppose that the $X$'s and $Y$'s are bounded by constants $R$ and $L$, respectively. Let the estimator $\wh{\beta}$ be as defined in \eqref{eq:estimator_ell1} with 
\begin{equation*}
\lambda \geq \gamma \lambda_0, \quad \gamma > 2,  \quad \lambda_0 \coloneq C (LR + R^2 ) \left( \sqrt{\frac{\log d}{n}} + \frac{\log d}{n} \right),
\end{equation*}
where $C > 0$ is a universal constant. If 
\begin{equation*}
n \geq \max \left\{ \frac{\lambda^{-1} \log(d) (4\sqrt{s} B) \wt{C}^2 (\nnorm{\Sigma}_{\text{{\em op}}}^{1/2} + R)  R^2}{\lambda_{\min}(\Sigma) \, (1/2 - 1/\gamma)},  \frac{16 R^4}{\lambda_{\min}^2(\Sigma)} \log d \right \} 
\end{equation*}
we have with probability at least $1 - 6/d$
\begin{alignat*}{2}
&\text{{\em (i)}} \;\; &&\nnorm{\wh{\beta} - \beta_*}_2 \leq \max \left\{\frac{24}{\lambda_{\min}(\Sigma)} \lambda \sqrt{s}, \, \frac{\wt{C} B  (\nnorm{\Sigma
}_{\text{{\em op}}}^{1/2} + R)}{\{ \lambda_{\min}(\Sigma) \}^{1/2}} \sqrt{\frac{s \log d}{n}}   \right\}, \\ 
&\text{{\em (ii)}} \;\; &&\nnorm{\wh{\beta} - \beta_*}_1 \leq  \max \left\{
\frac{96}{\lambda_{\min}(\Sigma)} \lambda s, \left( 4 + \frac{[\wt{C} \vee \wt{C}']^2 B^2  (\nnorm{\Sigma
}_{\text{{\em op}}}^{1/2} + R)^2}{\lambda_{\min}(\Sigma)} \right) s \sqrt{\frac{\log d}{n}}\right\}
\end{alignat*}
where $\wt{C},\wt{C}' > 0$ are constants depending only on the sub-Gaussian norm $\overline{K}$ of the $X$'s. 
\end{theo}
\noindent Note that as long as $\log d \leq n$, we may choose 
$\lambda \asymp \sqrt{\log(d)/n}$. Substituting this choice into the lower
bound required for $n$, we obtain the condition $n \gtrsim s \log d$, where 
the constant hidden in $\gtrsim$ is proportional to $B^2$, $\nnorm{\Sigma}_{\text{op}} / \lambda_{\min}^2(\Sigma)$, and $R^4$. The dependence
on these constants may not be optimal, but otherwise, the condition on $n$ and the bound on the estimation error in the $\ell_2$ and $\ell_1$ norms align
with those in the literature on sparse linear regression \cite[e.g.,][]{Raskutti2011, Wainwright2019, Negahban2012, Bickel2009}. Finally, observe that Theorem \ref{theo:lasso} requires bounded $X$ and $Y$ and fixed quantizer ranges. This requirement is likely not necessary; it is made for convenience here since it eliminates the need to deal with truncated random variables (cf.~Proposition \ref{prop:bias_control_truncation}), which cannot be addressed with the approach used in the low-dimensional case. 

\subsection{Debiasing and confidence intervals for the Lasso}
For the purpose of constructing asymptotic confidence intervals for the individual entries of $\beta_*$, we adopt the debiasing approach pioneered in \cite{Zhang2014} and further developed in \cite{Javanmard2014, vandegeer2014}. The following statement provides the underlying idea and justification.
\begin{prop}\label{prop:debias_lasso} Consider the setting in Theorem \ref{theo:lasso} with $\lambda \asymp \sqrt{\log(d)/n}$. Additionally, let $(M_n)$ be a fixed sequence of $d$-by-$d$ matrices such that 
\begin{equation*}
\nnorm{M_n \wh{\Sigma} - I_d}_{\infty} = O_{\p} \left(\sqrt{\frac{\log d}{n}} \right), \quad e_j^{\T} M_n^{\T} \mathbb{V} M_n e_j \rightarrow \vartheta_j, \; \; j=1,\ldots,d, 
\end{equation*}
as $n \rightarrow \infty$, where $\mathbb{V} = \cov(\wt{X}_1 \wt{Y}_1 - (\wt{X}_1 \wt{X}_1^{\T} + \Delta_1) \beta_*)$ with $\Delta_1$ as defined in Theorem \ref{theo:asymptotic}. Furthermore, suppose that $s = o(\sqrt{n}/\log d)$. It then holds that 
\begin{equation*}
\sqrt{n} (\wh{\beta}_j^{\text{{\em db}}} - \beta_{*j}) \overset{\text{{\em D}}}{\rightarrow} N(0, \vartheta_j), \;\; j=1,\ldots,d,      
\end{equation*}
where $\wh{\beta}^{\text{{\em db}}} = \wh{\beta} + M_n (\wh{\Sigma}_{Xy} - \wh{\Sigma} \wh{\beta})$.
\end{prop}
\noindent Note that the matrices $(M_n)$ in the above Proposition are assumed to be data-independent. In view of results in \cite{Javanmard2014}, choosing $M_n \equiv \Sigma^{-1}$ satisfies the conditions of the proposition. In practice, one computes an approximate inverse $\wh{M}_n$ of $\wh{\Sigma}$ satisfying $\nnorm{\wh{M}_n \wh{\Sigma}  - I}_{\infty} \lesssim \sqrt{\log(d)/n}$. Such an approximate inverse 
can be computed according to the optimization problem proposed in \cite{Javanmard2014} and implemented in \cite{sslasso}.  

\section{Empirical Results}\label{sec:empirical}
We provide empirical support for the theoretical results developed in the preceding sections. The first part of this section is dedicated to an exploration of (i) the 
relative efficiency of regression based on quantized data vs.~based on the full precision (i.e., unquantized) data, (ii) the statistical performance of the combined sketching/quantization pipeline presented in $\S$\ref{subsec:Sketch+Quant}, and (iii) savings in data transmission time when using this pipeline. In the second part of this section, we study asymptotic inference for the regression parameter including the case of $\ell_1$-penalization.
\begin{figure}[ht!]
\begin{center}
\begin{tabular}{cc}
\includegraphics[width = 0.42\textwidth]{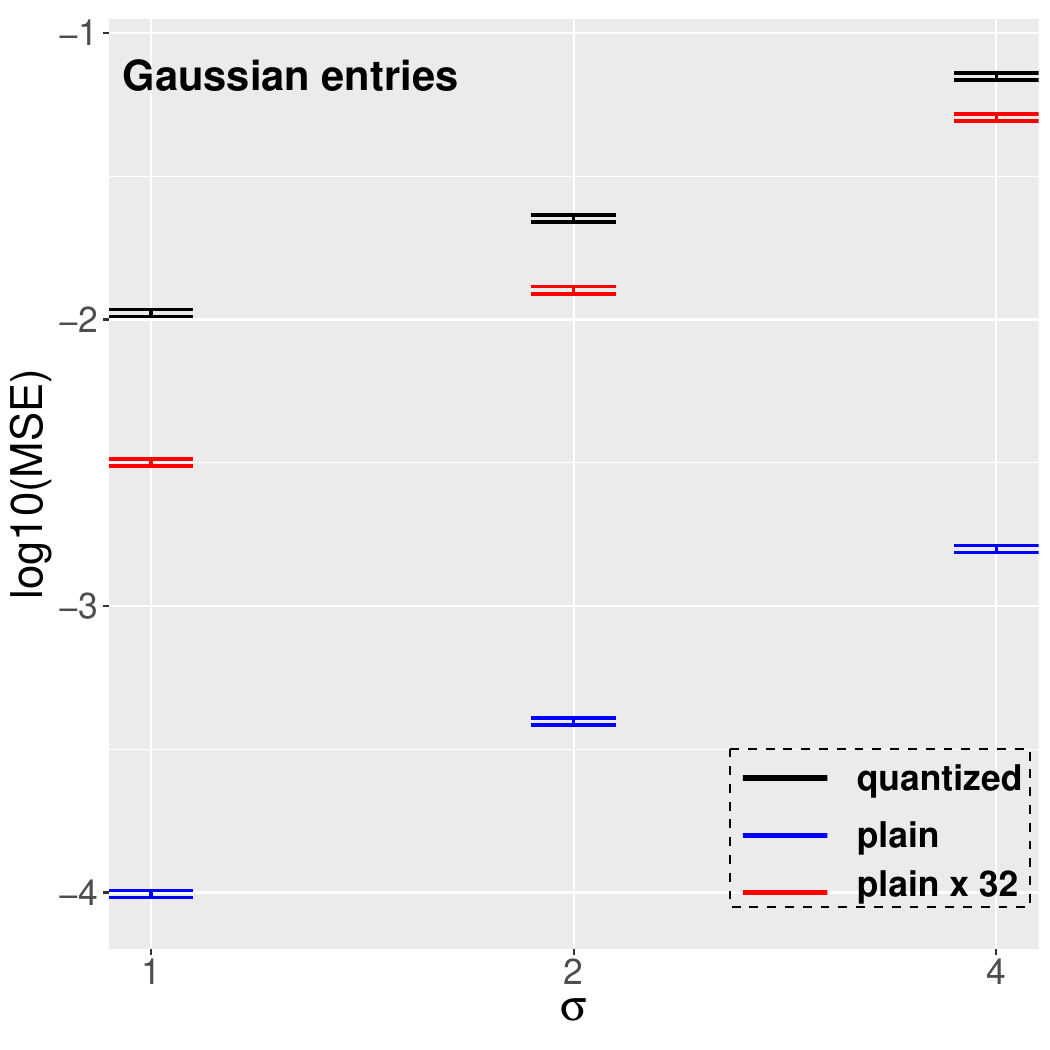} & \includegraphics[width = 0.42\textwidth]{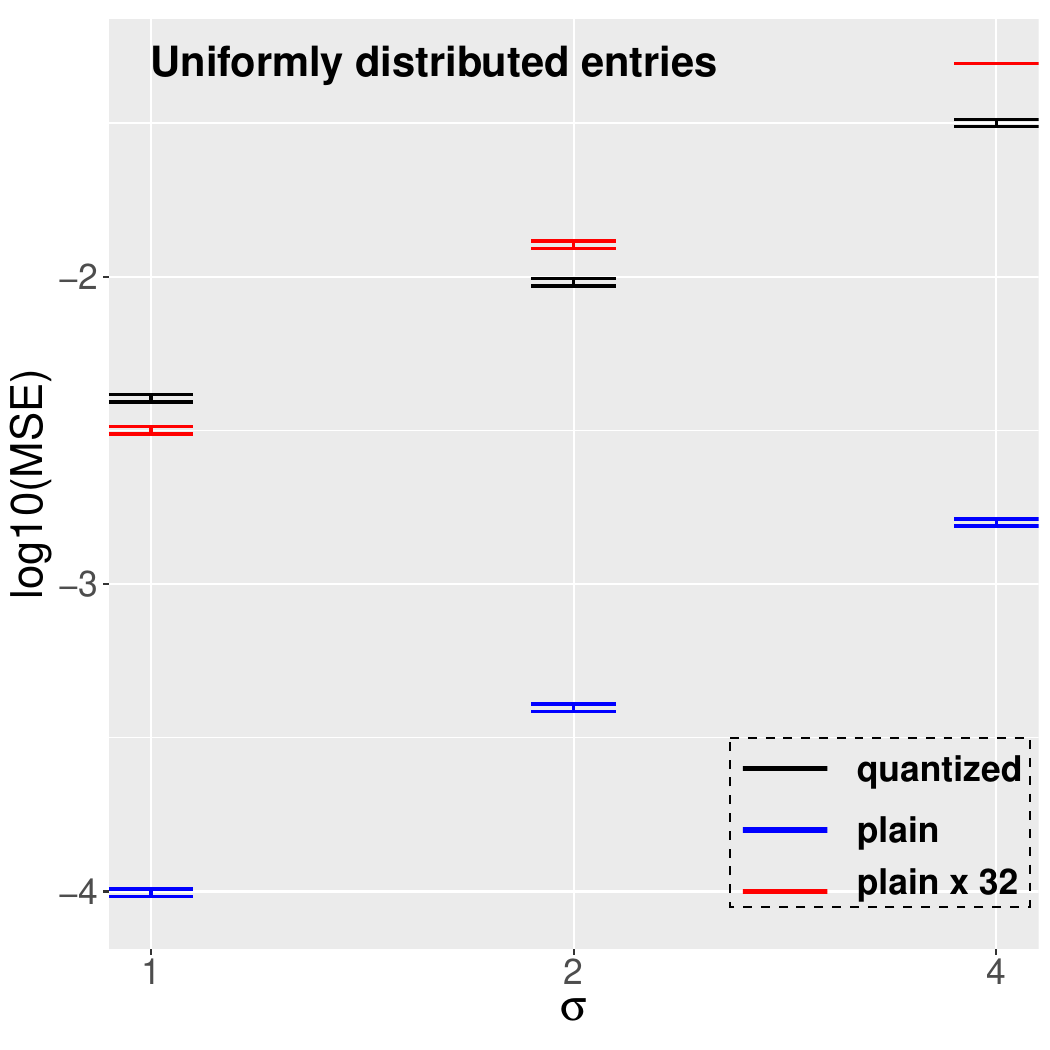}  
\end{tabular}
\begin{minipage}{\textwidth}
\vspace*{-7ex}
\begin{center}
\includegraphics[width = 0.49\textwidth]{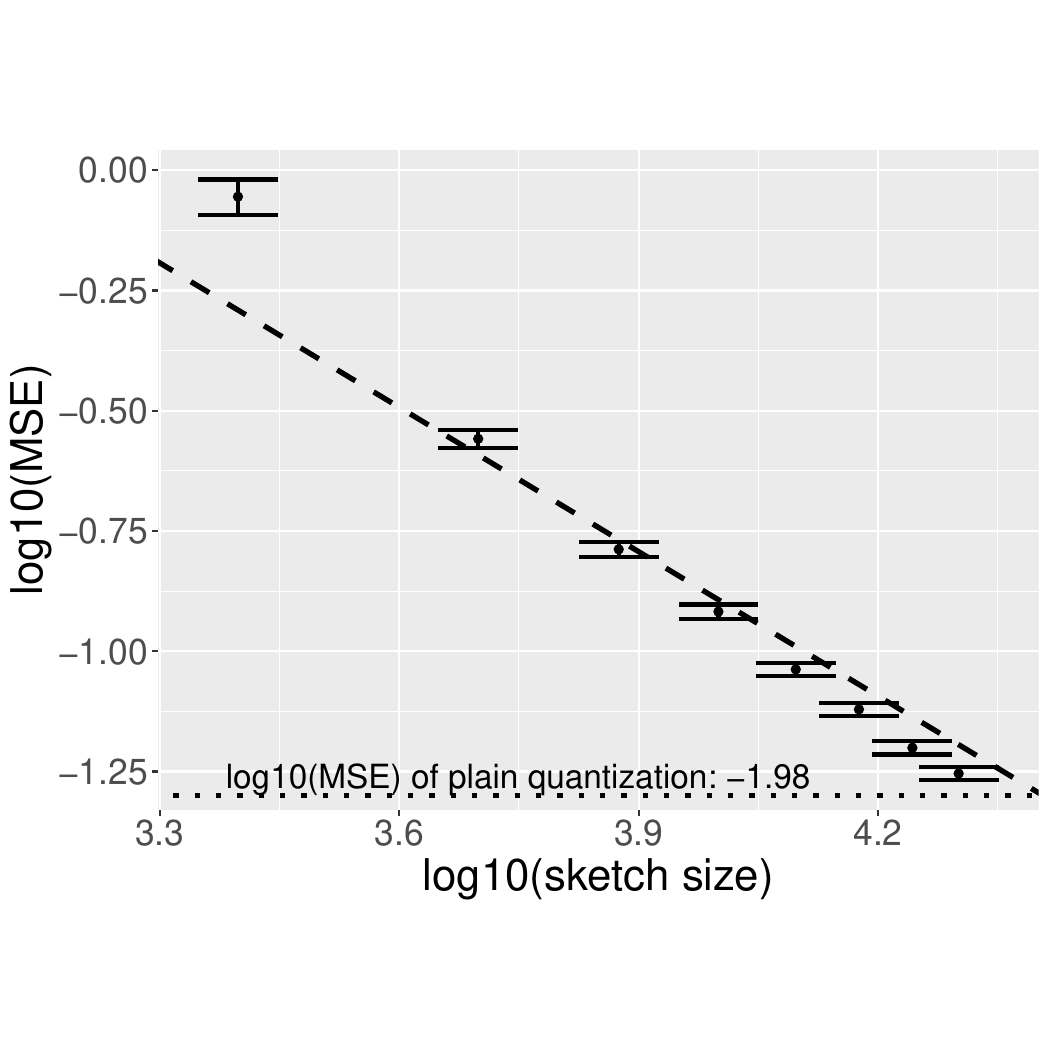}
\end{center}
\end{minipage}
\end{center}
\vspace*{-8.5ex}
\caption{Top: MSEs for estimating $\beta_*$ based on quantized and unquantized (``plain") data, including an adjustment (``plain $\times$ 32") to produce a comparison at the ``error per bit" level. The plots show means and $\pm 2\times$ standard error bars over 1k independent replications for each value of $\sigma$ (horizontal axis) and two scenarios -- (L) Gaussian $X$ and (R) uniformly distributed $X$. Bottom: MSE for estimating $\beta_*$ based on sketching + quantization under the first scenario (Gaussian $X$) and $\sigma = 1$ for varying sketch size $m$ (horizontal axis); the dashed straight line represents the least squares regression fit of $\log_{10} \text{MSE}$  on $\log m$ when fixing the slope to $m = -1$, while the dotted line indicates the MSE with quantization only (not to scale).}\label{fig:quant_sketch_mse}
\end{figure}

\subsection{Estimation error}
\begin{figure}[H]
    \centering
    \includegraphics[width=.8\linewidth]{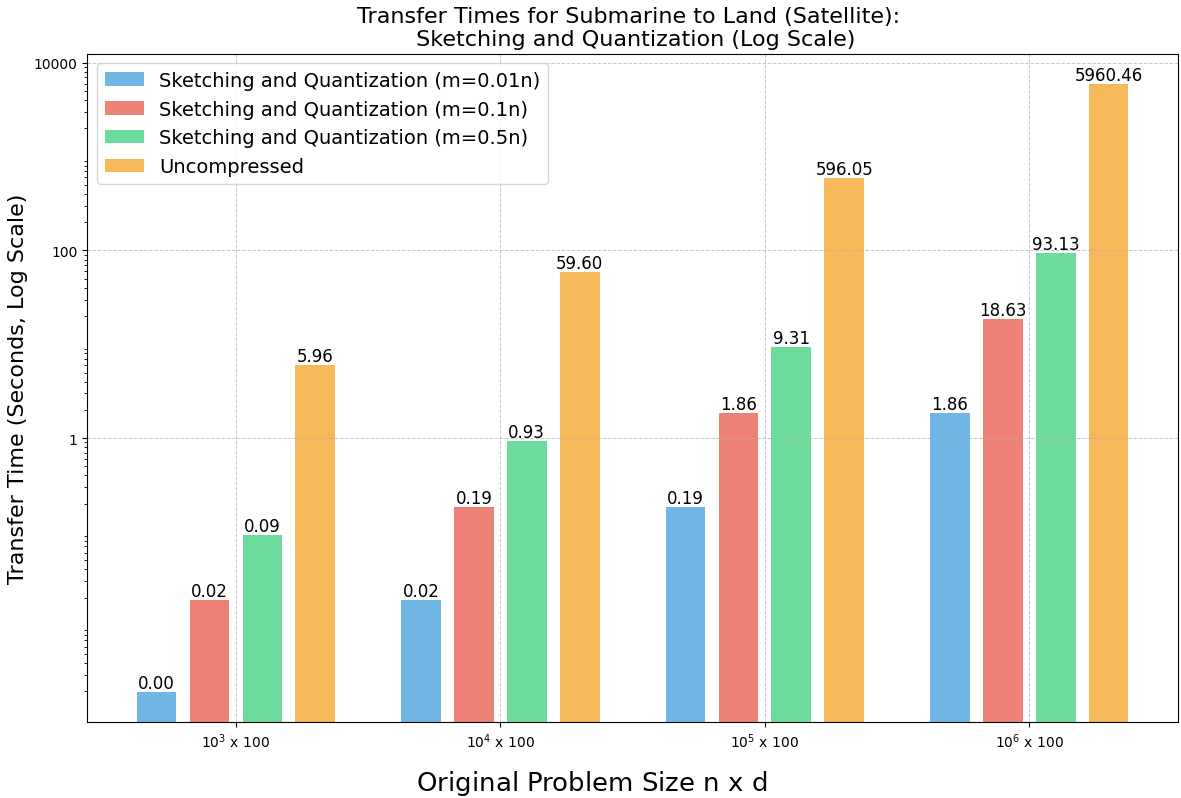}
    \caption{Example of low-bandwidth regression using compressed sufficient statistics.}
    \label{fig:submarine}
\vspace*{-1.5ex}
\end{figure}

We generate data from the linear model $Y = X^{\T} \beta_* + \sigma \eps$, where 
$\eps \sim N(0,1)$ is independent of $X$ and $\sigma \in \{1,2,4\}$ and $\beta_* \in \R^{10}$ is a unit vector with entries $\pm 0.15^{1/2}, \pm 0.15^{1/2}, \pm 0.1^{1/2}$, \\
$\pm 0.05^{1/2}, \pm 0.05^{1/2}$. Two scenarios
are considered for $X$: (i) the entries of are i.i.d.~$N(0,1)$ and (ii) the entries are i.i.d.~from the uniform distribution on $[-\sqrt{3}, \sqrt{3}]$. We use fixed quantizer ranges $R = 2.5$ for $X$ and $X^2$
in the first scenario ($R = \sqrt{3}$ in the second scenario), and $L = 2.5 \cdot \sqrt{\sigma^2 + \nnorm{\beta_*}_2^2}$ in both scenarios, thereby slightly departing from the choice of these ranges for unbounded sub-Gaussian data considered 
in our analysis. These less conservative choices yield improved
empirical performance by trading a rather small amount of bias against substantially reduced variance. The number of samples is set to 
$n=100$k. Under scenario (i) and $\sigma = 1$, we additionally consider sketching 
prior to quantization, adopting the protocol sketched in Figure \ref{fig:sketching_plus_quantization}: the entries of the $m$-by-$n$ sketching matrix $\M{S}$ are sampled from the $N(0,1/m)$ distribution, with $m$ varying between $2.5$k and $20$k in steps of $2.5$k; the sketched data is subsequently quantized according to the above choice
of the quantizer ranges. 1,000 independent replications are performed for each scenario and each value of $\sigma$. 
\vskip1.5ex
\noindent {\em Results}. Figure \ref{fig:quant_sketch_mse} provides a summary of the results. In addition to comparing the ``raw" MSEs $\nnorm{\cdot - \beta_*}_2^2$ that
result from using quantized and unquantized data, we also report the 
latter scaled up by $32$ to roughly account for the larger number of 
bits used in this case, assuming double precision, i.e., 64 bits relative to the 2 bits used for each entry of $X$ when quantizing. In the resulting comparison, the gap between the quantized and (adjusted)
unquantized MSE narrows as $\sigma$ increases, in agreement with the folklore that the smaller the ratio of signal ($\nnorm{\beta_*}$) to noise the more suitable it is to apply coarse quantization \cite{LaskaBaraniuk2011}. For Gaussian $X$, not
using quantization still performs better for all values of $\sigma$ after adjustment, whereas for uniformly distributed $X$, quantization
becomes slightly superior for $\sigma = 4$. The bottom panel of Figure \ref{fig:quant_sketch_mse} corroborates the rate $1/m$ for the MSE as function of the number of sketches as hypothesized based on the discussion in $\S$\ref{subsec:Sketch+Quant}. For $m = 2.5$k, there is 
an apparent departure from the straight line, which indicates that 
the overall amount of compression is too stringent to yield 
alignment. In general, the sample size required to achieve substantially
small MSEs even with coarse quantization tends to be large (10k or higher), which confirms that the use of coarse quantization is most appropriate when data is abundant. 
\vskip1.5ex
\noindent {\em Transmission times}. In the preceding discussion, we have focused on estimation error vs.~storage/transmission cost incurred. To give an indication of how transmission volumes can also impact transmission time, we consider an illustrative scenario in which data must be transmitted from a submarine to a remote server through a narrow-band satellite uplink operating at approximately 0.5~Mbps.  In this environment, minimizing transmission time is critical to reduce surface exposure and maintain stealth.  Figure~\ref{fig:submarine} shows the simulated upload time required to transmit design matrices of varying sizes under compressed and uncompressed schemes.  The results show that sketching and quantization can shorten transmission time by several orders of magnitude for large-scale data, enabling real-time or near-real-time communication that would otherwise be infeasible. When interpreted alongside the MSE results above, these findings demonstrate that compression can achieve favorable trade-offs: controlled degradation in estimator accuracy for a substantial reduction in communication latency. 



\subsection{Inference}
\begin{figure}
\begin{tabular}{|c|cc|}
\hline
{\small Low-dimensional Setup} & {\small Moderate-dimensional Setup}  &  {\small Moderate-dimensional Setup} \\
                               & {\small (non-zero coefficients)}     & {\small (first four zero coefficients)}  \\
\includegraphics[height = .32\textwidth]{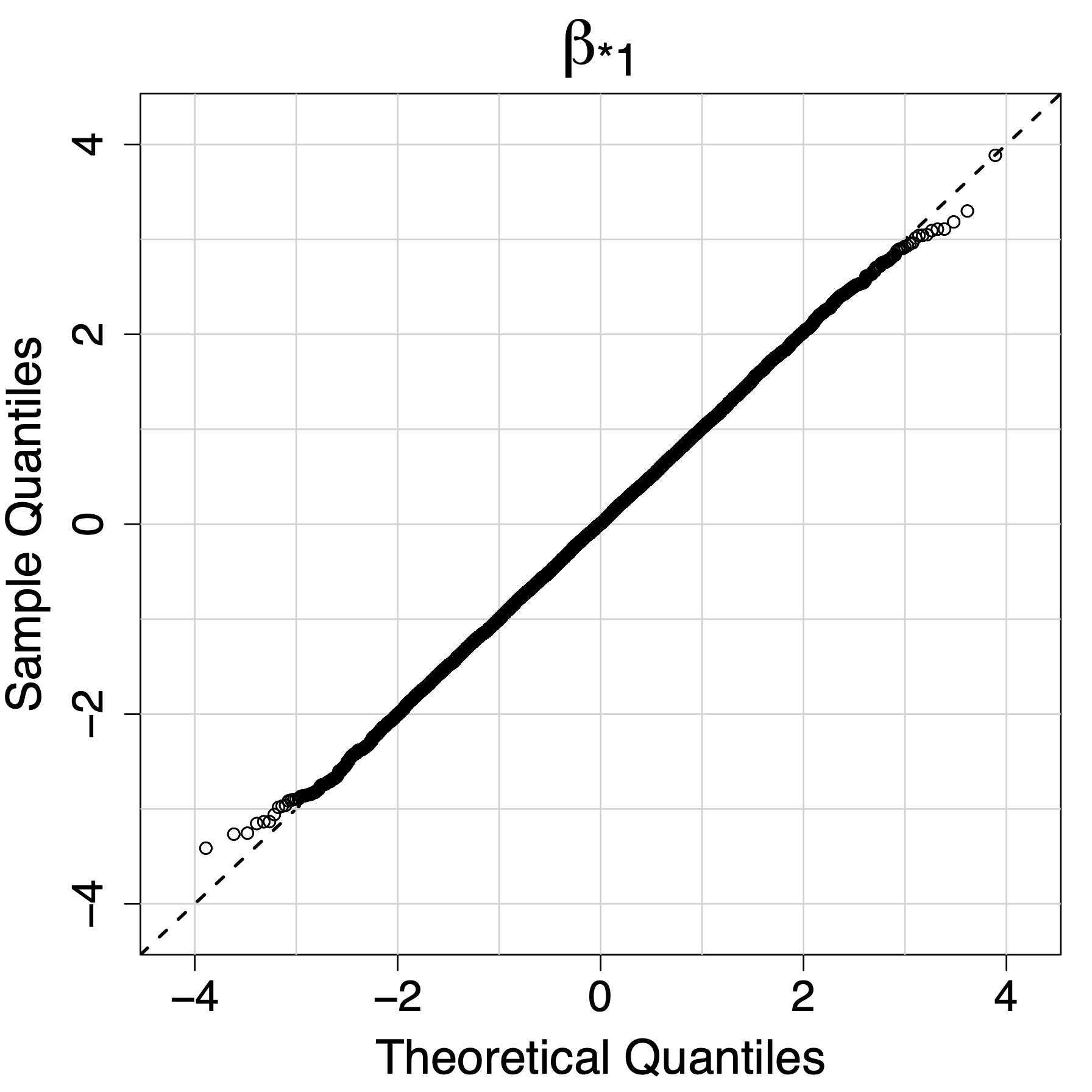} & \includegraphics[height = .32\textwidth]{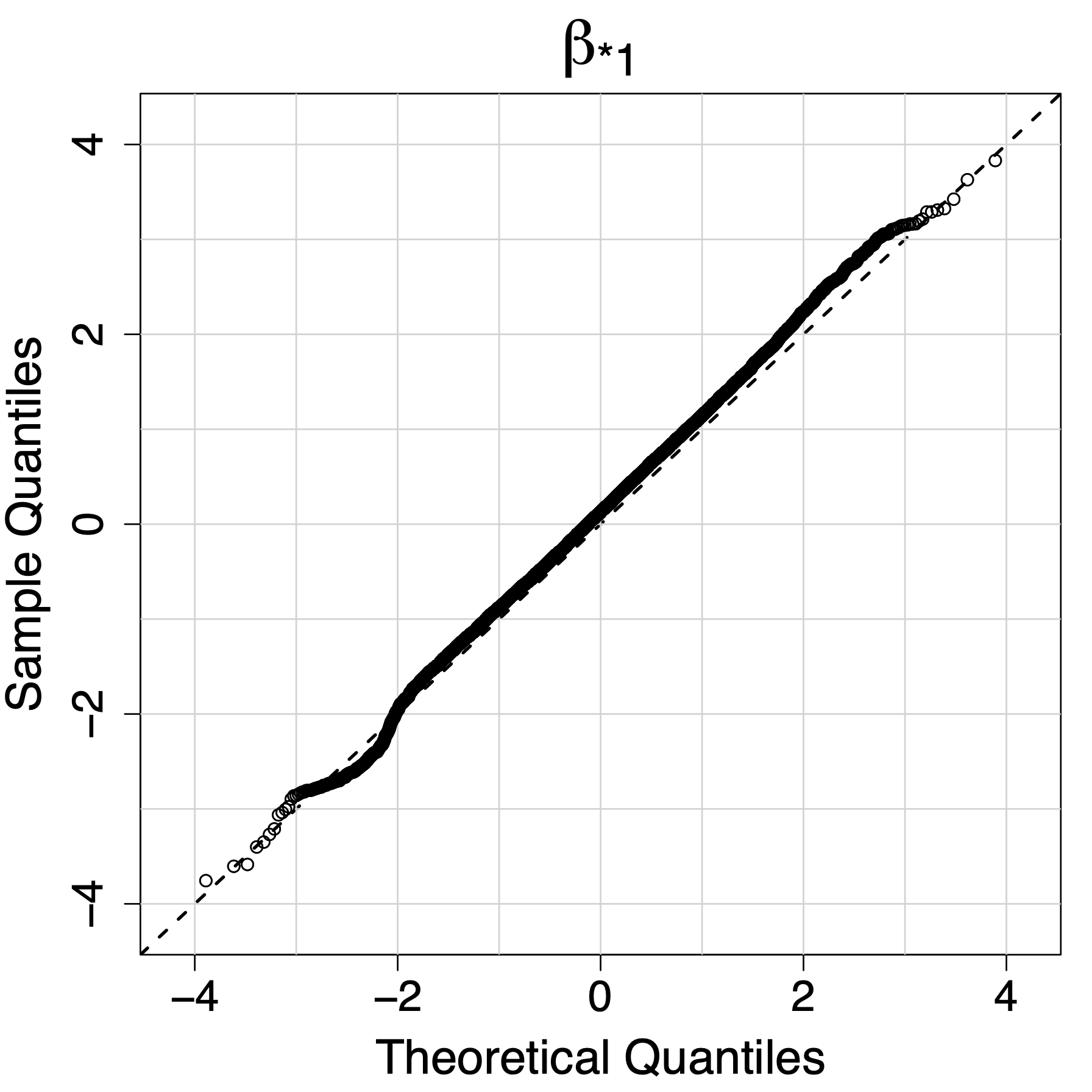}   & \includegraphics[height = .32\textwidth]{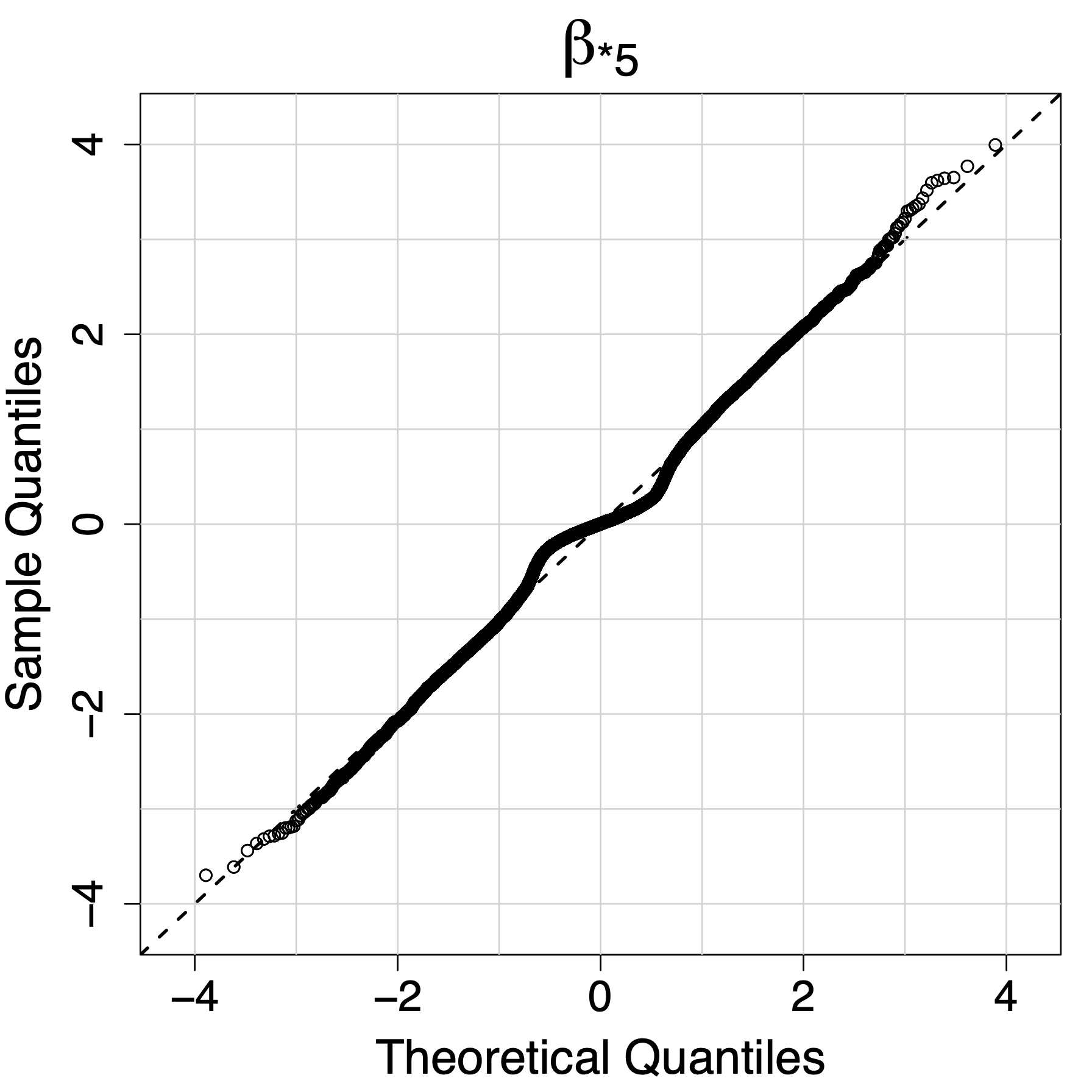} \\
\includegraphics[height = .32\textwidth]{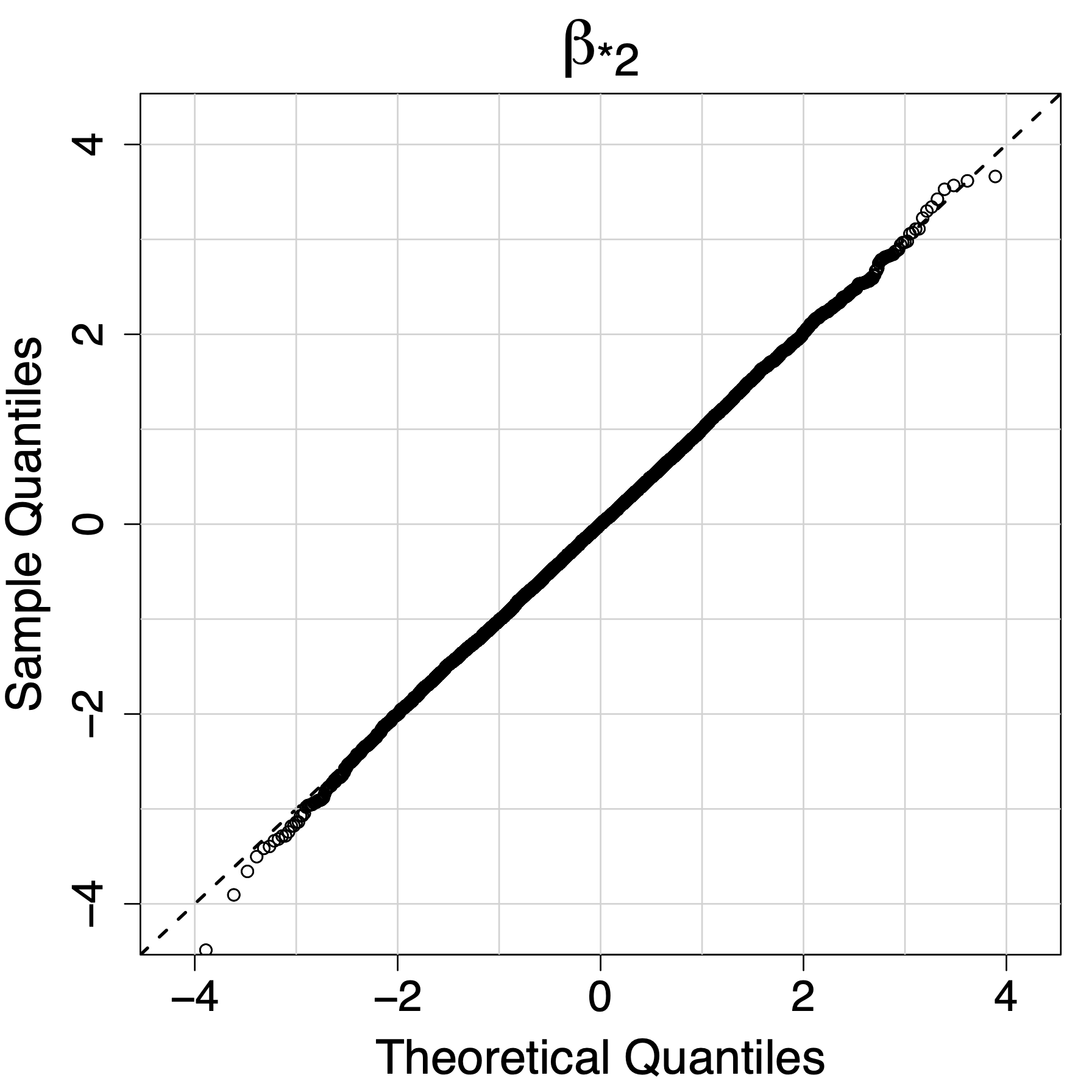} & \includegraphics[height = .32\textwidth]{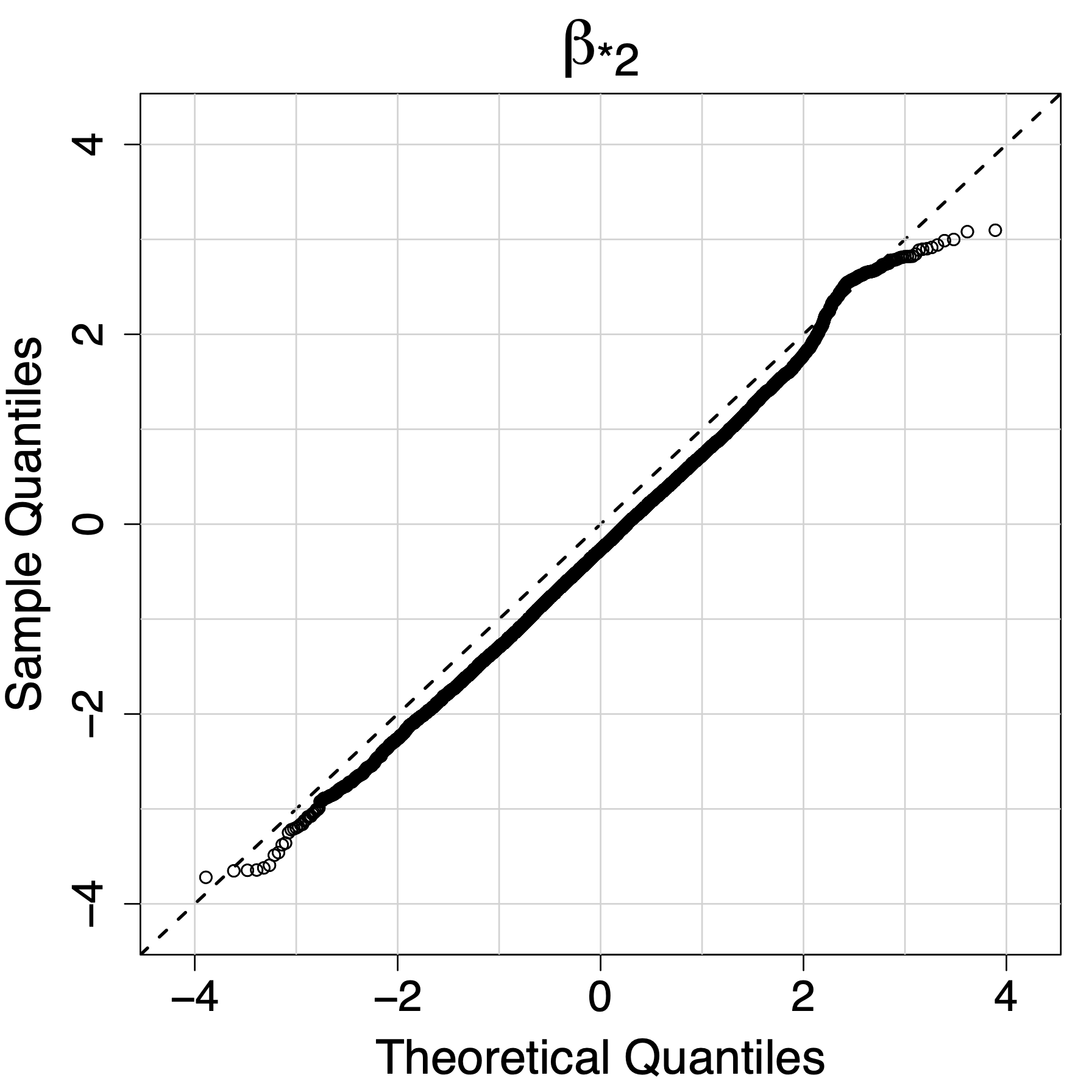}  & \includegraphics[height = .32\textwidth]{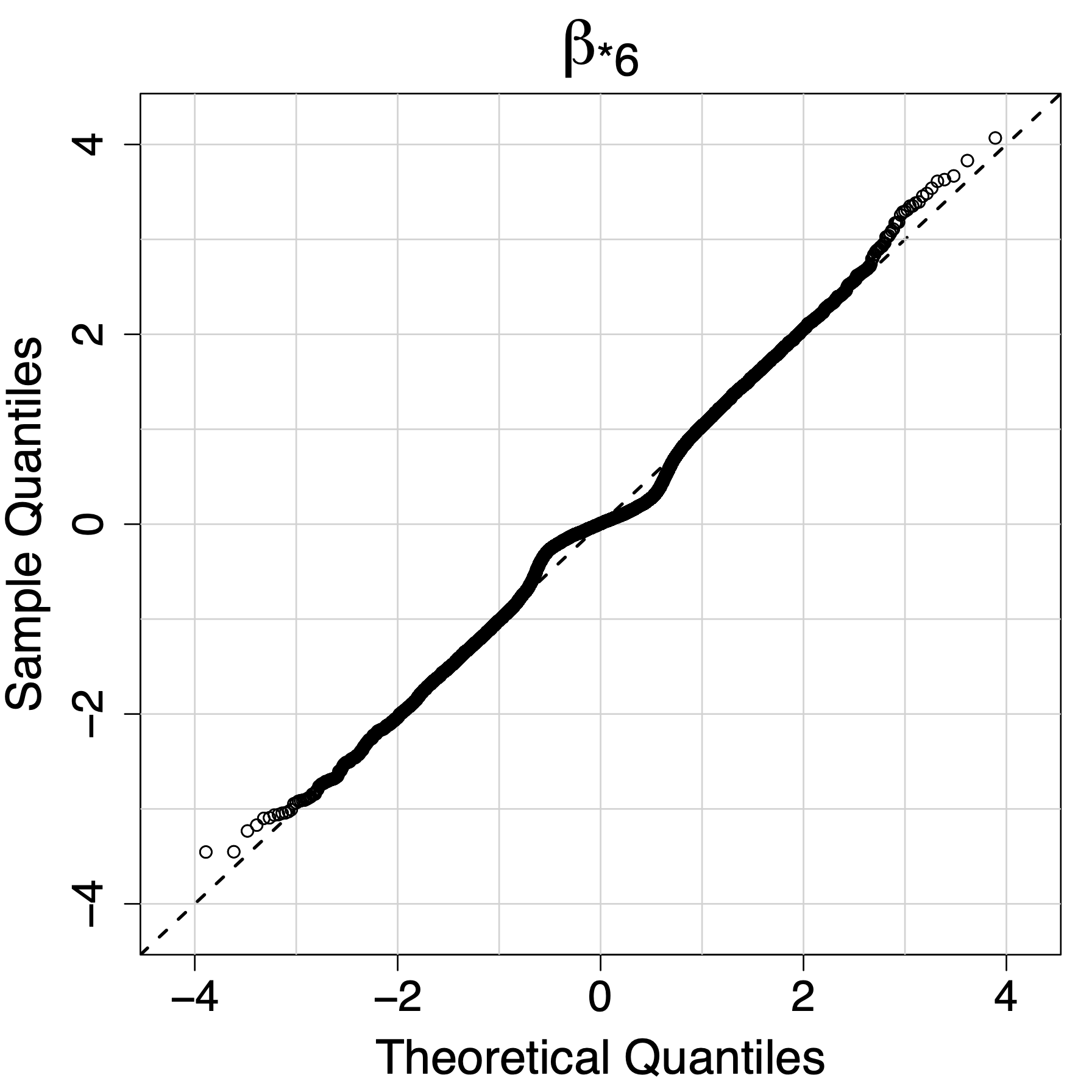} \\
\includegraphics[height = .32\textwidth]{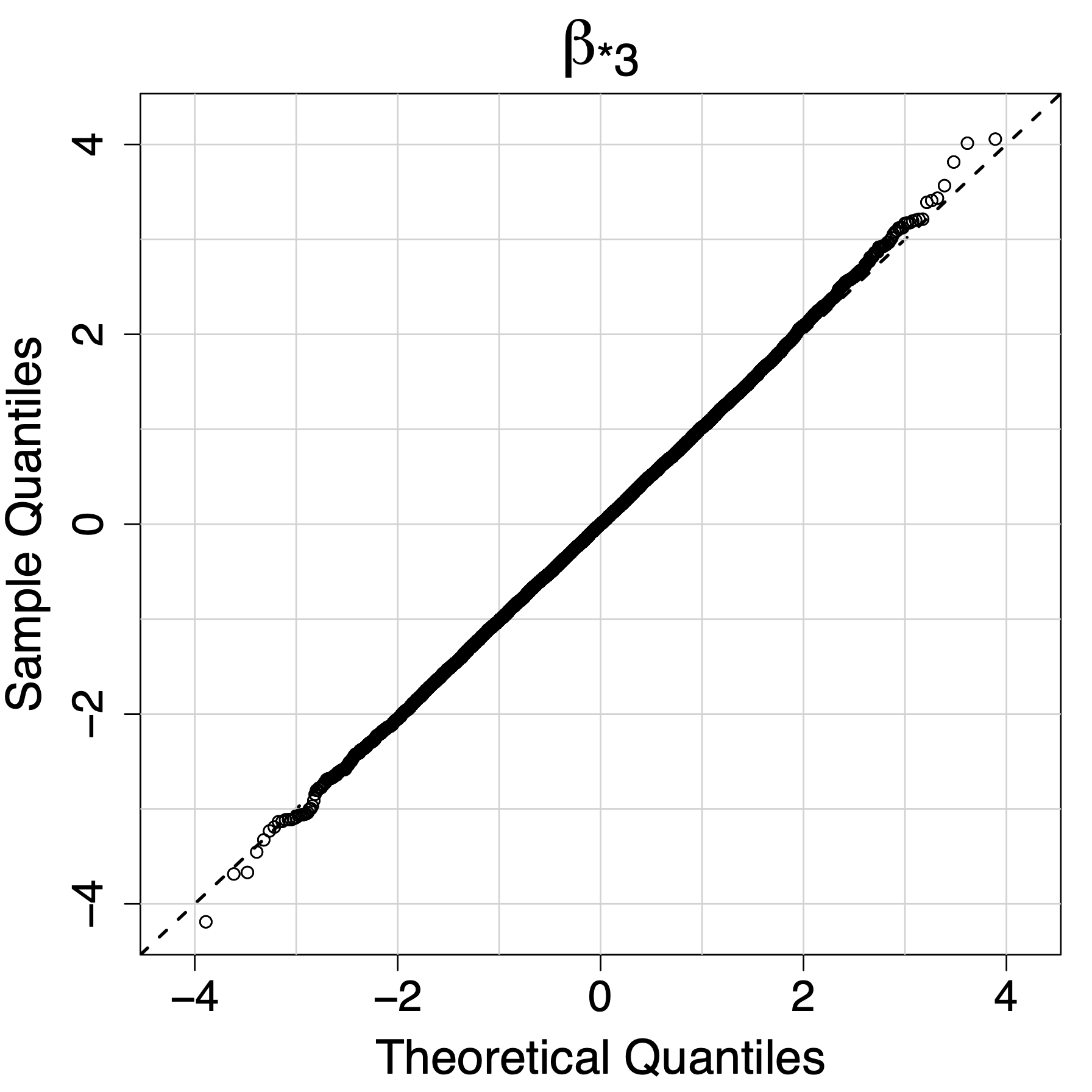} & \includegraphics[height = .32\textwidth]{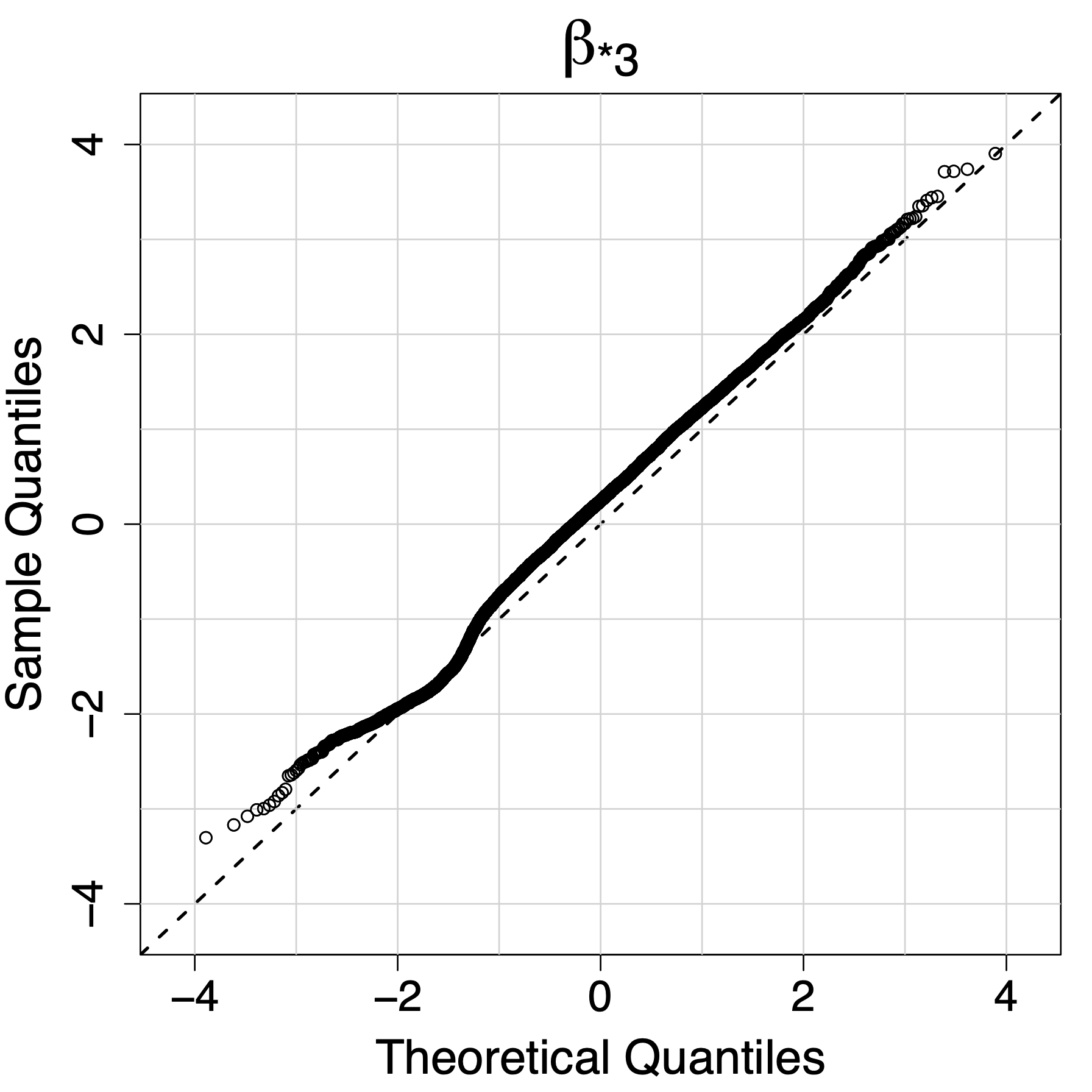}  & \includegraphics[height = .32\textwidth]{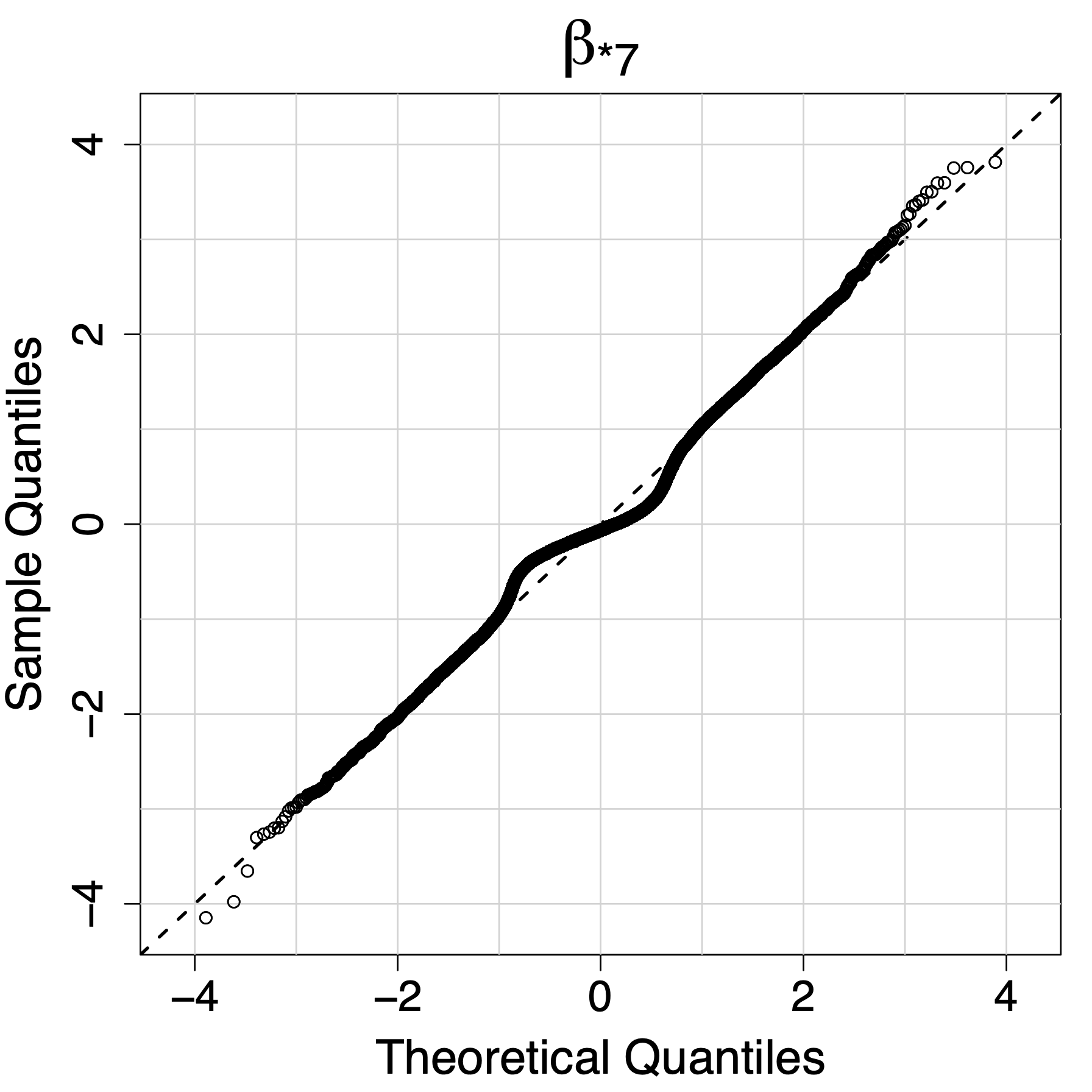}\\
\includegraphics[height = .32\textwidth]{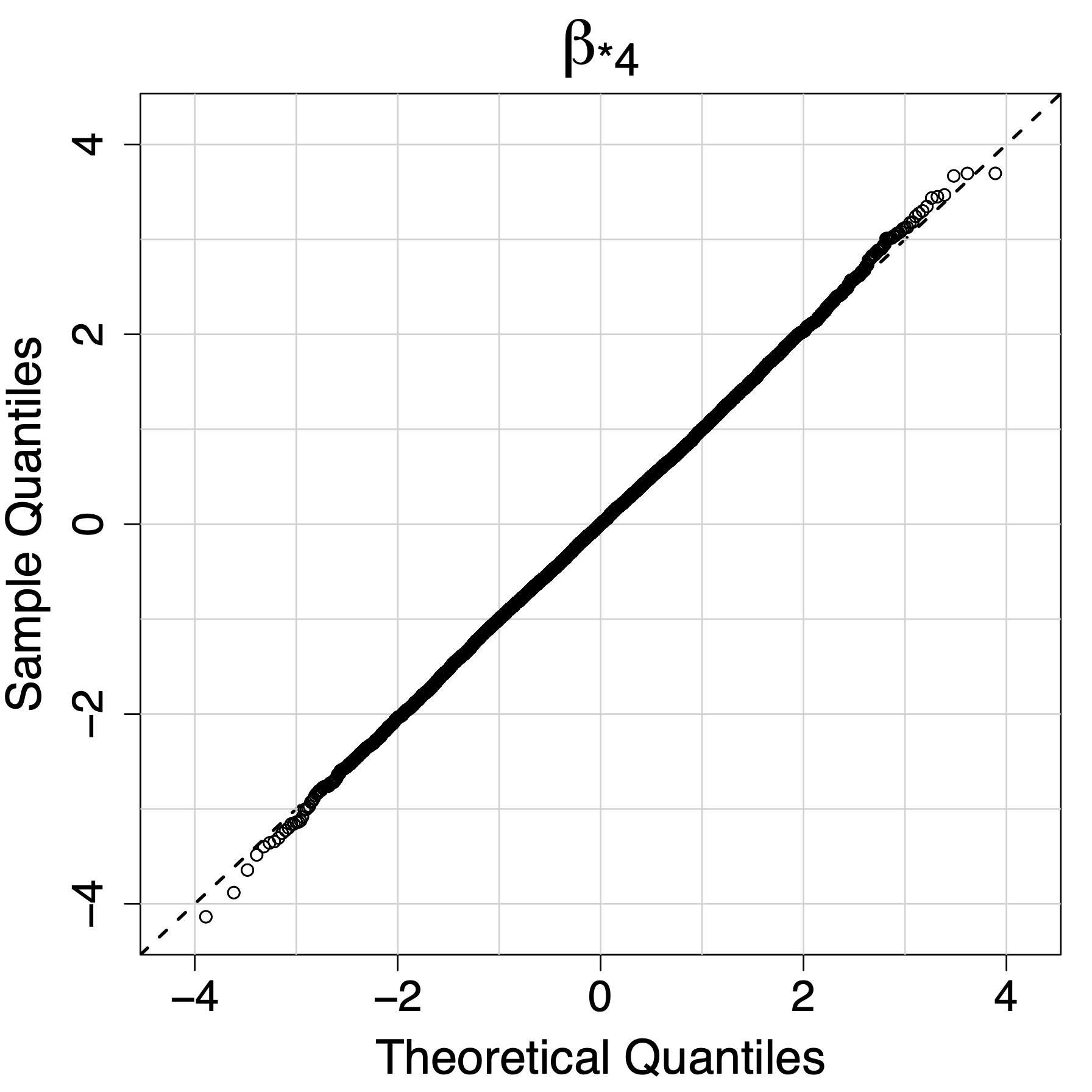}  & \includegraphics[height = .32\textwidth]{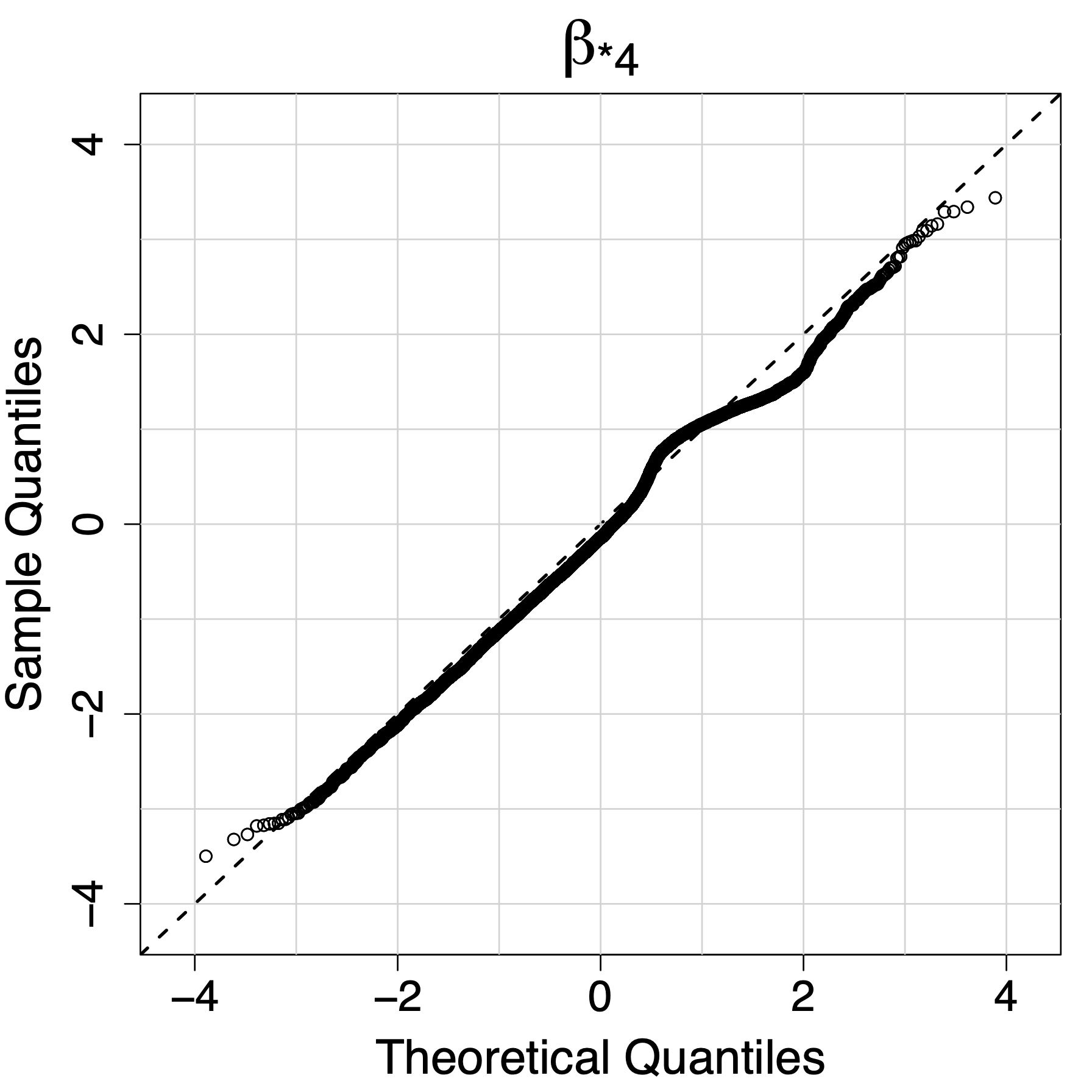}  &  \includegraphics[height = .32\textwidth]{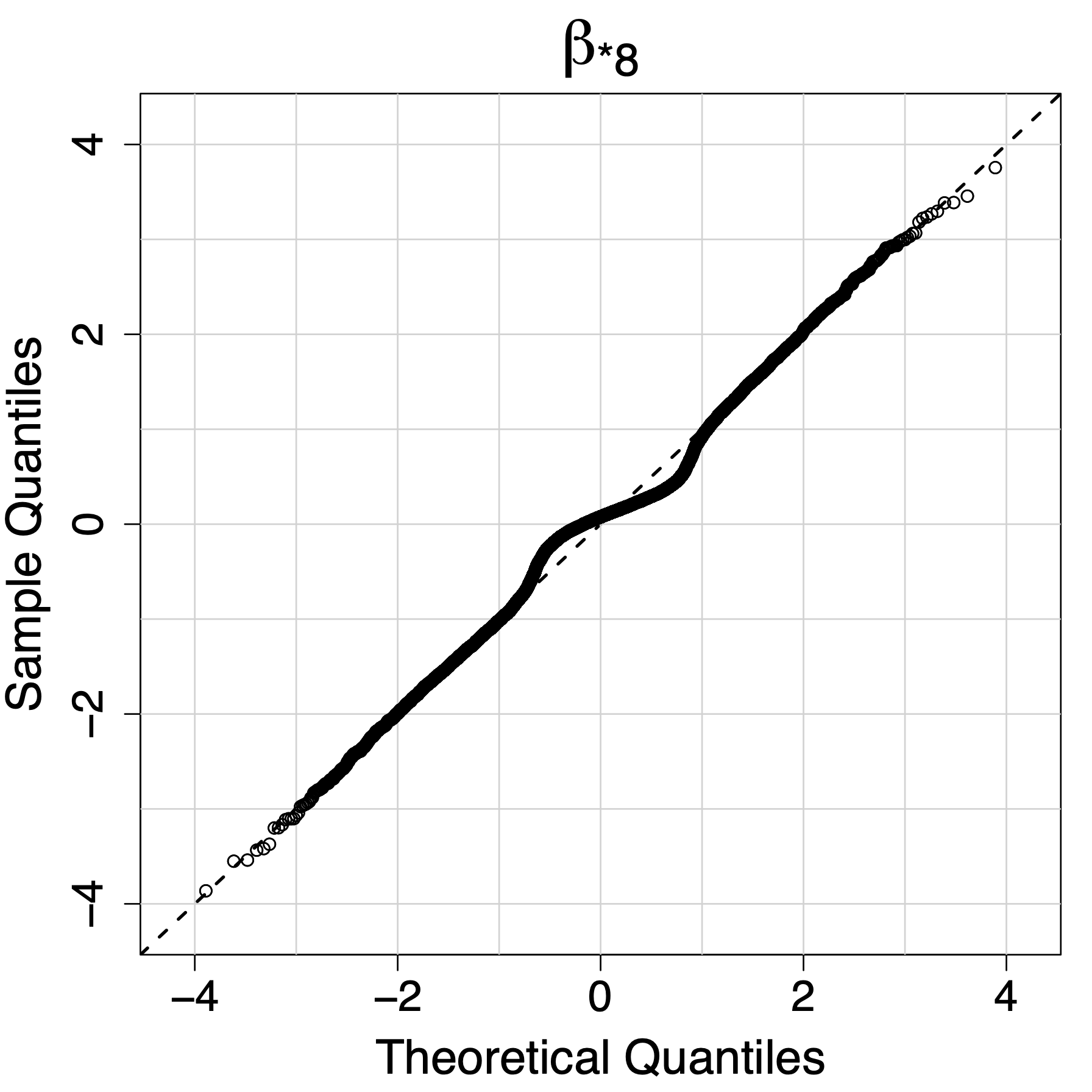} \\
\hline
\end{tabular}
\caption{Normal Q-Q plots of the centered and scaled regression coefficient estimates in the low-dimensional setting (left) and in the sparse, moderate dimensional setting when using $\ell_1$-penalization.}\label{fig:qqplot-low}
\end{figure}

In this subsection, we present two sets of simulations that investigate inference 
for the individual regression coefficients via asymptotic Normality in (i) the low-dimensional case, and (ii) the moderate dimensional and sparse case when using $\ell_1$-penalization, as asserted by Theorem \ref{theo:asymptotic} and Proposition \ref{prop:debias_lasso}, respectively.
\vskip1.5ex
\noindent {\em Low-dimensional setup}. We let $d = 4$ with $X_1 \sim 2 \cdot \mathrm{Beta}(1,1) - 1$, $X_2 \sim 2 \cdot \mathrm{Beta}(2,2) - 1$, $X_3 \sim 2\cdot\mathrm{Beta}(1,4) - 1$, 
$X_4 \sim 2\cdot\mathrm{Beta}(1,4) - 1$ with $\{ X_j \}_{j = 1}^4$ independent. Here, $\mathrm{Beta}(a,b)$ denotes the Beta distribution with parameters $a$ and $b$, respectively, whose probability density function $f_{a,b}$ takes the form $f_{a,b}(x) = \frac{\Gamma(a + b)}{\Gamma(a) \Gamma(b)} x^{a-1} (1 - x)^{b-1}, \; x \in (0,1)$, where $\Gamma$ denotes the Gamma function. The Beta family is a standard 
model for interval-valued data; we here scale by $2$ and translate by $-1$ to obtain the symmetric (around zero) interval $[-1,1]$ as the support. The response variable is generated as 
\begin{equation*}
Y = X^{\T} \beta_* + \sigma \eps, \quad \eps \sim N(0,1), \quad \sigma = \sqrt{1/2},     
\end{equation*}
with $\eps$ independent of $X$ and $\beta_* = (0.5, -\sqrt{3/4}, \sqrt{3/4}, -0.5)^{\T}$ so that 
$\nnorm{\beta_*}_2^2 = 2$. Each replication (out of 10k in total) is run with $n = 10^4$ samples, i.e., pairs $\{ (X_i, Y_i) \}$ each distributed as $(X,Y)$ with $X$ and $Y$ as described above.  With this setup, we choose $R = 1$ and $L = \sqrt{d} \nnorm{\beta_*} + \sigma \sqrt{2 \log n}$, which bounds $\max_{1 \leq i \leq n} |Y_i|$ with high probability. With these choices, the $\{ (X_i, Y_i) \}$ are quantized following the proposed protocol, and the quantized samples are used to obtain the estimator $\wh{\beta}$. We evaluate standard errors $\text{se}(\wh{\beta}_j) = 
n^{-1/2} e_j \Sigma^{-1}  \Gamma \Sigma^{-1} e_j$ according to Theorem \ref{theo:asymptotic} and construct confidence 95\% confidence intervals via $\wh{\beta}_j \pm 1.96 \cdot \text{se}(\wh{\beta}_j)$, $j=1,\ldots,d$, where $1.96$ approximately equals the $.975$-quantile of the $N(0,1)$-distribution. Deriving closed-form expressions for the standard errors is feasible, but tedious; we omit the details here. In the example, their estimated counterparts (cf.~\eqref{eq:est_standard_error}) are rather close to these closed-form standard errors. Figure \ref{fig:qqplot-low} show the Normal Q-Q plots in which the quantiles of the empirical distributions
of $(\wh{\beta}_j - \beta_{*j})/\text{se}(\wh{\beta}_j)$ over the 10k replications, $j=1,\ldots,d$, are compared to the quantiles of the standard Normal distribution. The alignment in the Q-Q plots is close to exact, and the confidence interval coverage agrees with the nominal level (cf.~Table \ref{tab:CI-coverage}). 

\vskip1.5ex
\noindent {\em Moderate-dimensional setup}. We investigate $\ell_1$-penalization based on formulation \eqref{eq:estimator_ell1}. The setup is that of the low-dimensional setup, with the following differences: (i) we add $36$ additional 
variables $X_5, \ldots, X_{40}$ with $X_{4k + 1}$ distributed as $X_1$, 
$X_{4k + 2}$ distributed as $X_2$,  $X_{4k + 3}$ distributed as $X_3$ and $X_{4k + 4}$ distributed as $X_4$, $k=1,\ldots,9$, and set $\beta_{*5} = \ldots = \beta_{*40} = 0$, (ii) estimates $(\wh{\beta}_j)_{j = 1}^d$ are computed by solving optimization problem \eqref{eq:estimator_ell1} without the $\ell_1$-ball constraint (i.e., $B = \infty$) and $\lambda = 2 \sqrt{\log(d)/n}$, (iii) subsequently, we obtain debiased regression coefficients $(\wh{\beta}_j^{\text{db}})_{j = 1}^d$ by adopting the approach outlined in Proposition \ref{prop:debias_lasso}, where 
$M_n$ is obtained as a matrix satisfying $\nnorm{M_n \wh{\Sigma} - I_d}_{\infty} \leq .1 \sqrt{\log(d)/n}$ according to the method \cite{Javanmard2014} and implementation \cite{sslasso}. We report bias, standard deviation, and confidence intervals of these debiased estimates in Table \ref{tab:CI-coverage} and Normal Q-Q plots in Figure \ref{fig:qqplot-low} (for space reasons, we confine ourselves to the first four coefficients outside the support of $\beta_*$). We observe that the debiased coefficients corresponding to the support still exhibit a noticeable bias, which is reflected by slight departures from the angle bisector in the Q-Q plots as well. However, this does not affect the confidence interval coverage; the latter is close to the nominal 95\%. We also note that the standard deviations of the first four coefficients are larger compared to the low-dimensional setup (which would correspond to knowing the support of $\beta_*$).    
\begin{table}
\centering
\begin{minipage}{.33\textwidth}
{\small \hspace*{3ex}$\text{Low-dimensional Setup}$ \\[1ex] \begin{tabular}{|l|c|c|c|c|}
\hline & & & & \\[-1ex]
& $\wh{\beta}_1$ & $\wh{\beta}_2$ & $\wh{\beta}_3$ & $\wh{\beta}_4$ \\[.5ex]
\hline & & & & \\[-2ex] 
Bias & $\sim0$ & $\sim0$ &  $\sim0$ & $\sim0$  \\[.5ex] 
\hline & & & & \\[-2ex] 
SD   & $.18$ & $.29$ & $.33$ & $.32$ \\[.5ex]
\hline & & & & \\[-2ex] 
COV  & $.95$ & $.95$ & $.94$ & $.94$ \\[.5ex]
\hline
\end{tabular}}
\end{minipage}
\hspace*{1ex}\begin{minipage}{.65\textwidth}
{\small \hspace*{15ex} $\text{Moderate-dimensional Setup}$ \\[1ex]\begin{tabular}{|l|c|c|c|c|c|c|c|c|c|}
\hline & & & & & & & & \\[-1ex]
& $\wh{\beta}_1^{\text{db}}$ & $\wh{\beta}_2^{\text{db}}$ & $\wh{\beta}_3^{\text{db}}$ & $\wh{\beta}_4^{\text{db}}$ & $\wh{\beta}_5^{\text{db}}$ & $\wh{\beta}_6^{\text{db}}$ & $\wh{\beta}_7^{\text{db}}$ & $\wh{\beta}_8^{\text{db}}$ \\[.5ex]
\hline & & & & & & & &\\[-2ex] 
Bias & $.03$ & $-.09$ & $.12$ & $-.05$ & $\sim0$ & $\sim0$ & $\sim0$ & $\sim0$  \\[.5ex] 
\hline & & & & & & & &\\[-2ex] 
SD   & $.21$ & $.36$ & $.57$ & $.57$ & $.20$ & $.35$ & $.55$ & $.54$ \\[.5ex]
\hline & & & & & & & & \\[-2ex] 
COV  & $.94$ &$.94$ & $.94$ &$.96$ & $.94$ & $.95$ & $.95$ & $.95$ \\[.5ex]
\hline
\end{tabular}}
\end{minipage}
\caption{Bias, standard deviations (SD), and confidence interval coverage (COV) of the regression coefficient estimates in the low-dimensional setup (left) and moderate-dimensional setup when using $\ell_1$-penalization followed by debiasing (right). For the latter, we only show the first four (out of $36$) coefficient estimates corresponding to the zeroes of $\beta_*$.}\label{tab:CI-coverage}
\end{table}

\section{Conclusion}
In this paper, we have studied the finite-sample and asymptotic performance of 
least squares estimation using exclusively one-bit information from the predictors and responses. Asymptotic inference and relative efficiency are key motivations for our study. A number of topics for future research arise from directions pursued in 
related work: these include (i) the use of the quantization scheme with 
triangular dither put forth in \cite{Chen2025cov, Chen2023heavy} and (ii) input
data with heavier tails than sub-Gaussian \cite{ Chen2023heavy}. Additional directions concern the use of general multi-bit dithered quantizers and quantization {\em without}
dither followed by likelihood-based inference and associated estimation techniques such as Expectation-Maximization (EM) as alluded to in $\S$\ref{sec:intro}. The latter direction is motivated by the sub-optimal efficiency arising from dithering, which generally introduces additional logarithmic factors in $n$ into the estimation error
rate; the discussion in $\S$\ref{subsec:lowerbound} suggests that these factors arise 
as a result of the quantization scheme rather than being a consequence of the specific 
estimator under consideration herein. 

While the present work is motivated by resource constraints, a potential additional benefit arising from data compression is enhanced data privacy. Note that quantization produces data with bounded range, which is beneficial for controlling sensitivity in differential privacy \cite[][Ch.~3]{Dwork2014}. The privacy aspect has been
studied for sketching \cite{Zhou2009, PilanciWainwright2015} as well as quantization; for the latter, formal privacy guarantees can be obtained by randomized label flipping, a mechanism commonly referred to as ``randomized response" \cite{Warner1965, Dwork2014}. Formal privacy guarantees for sketching followed by one-bit quantization and random flipping are developed in \cite{Li2023} albeit the utility in that work concerns information retrieval. Working out a similar analysis for the regression setup in the present paper thus remains an interesting avenue of future research as well.       

\bibliographystyle{plain}
\bibliography{quantreg}
\vskip8ex
\noindent {\Large {\bfseries \textsf{Appendix}}}

\renewcommand{\thesubsection}{\Alph{subsection}}

\subsection{Proof of Proposition \ref{prop:quantizer_bound}}
Regarding the event $\mc{R}$, the result follows from the union bound over the $\{X_{ij} \}$,
$1 \leq i \leq n$, $1 \leq j \leq d$, and invoking Lemma \ref{ref:tailbound}
with the choice $t = 2\sqrt{C_K} \cdot \sqrt{\log(d \cdot n)}$.

Regarding the event $\mc{L}$, we note that
\begin{equation*}
\nnorm{\nscp{X}{\beta_*}}_{\psi_2} = \nnorm{\nscp{\Sigma^{-1/2} X}{\Sigma^{1/2} \beta_*}}_{\psi_2} \leq \nnorm{X}_{\psi_2} \nnorm{\Sigma^{1/2} \beta_*}_2 
\end{equation*}
Using that for any fixed $i$, $|Y_i| \leq |\nscp{X_i}{\beta
^*}| + \sigma |\eps_i|$, a two-fold application of Lemma \ref{ref:tailbound}
with $t = 2\sqrt{C_{\overline{K}}} \nnorm{\Sigma^{1/2} \beta_*}_2 \cdot \sqrt{\log n}$ and $t = 2\sqrt{C_{K_{\eps}}} \cdot \sqrt{\log n}$, respectively, yields $|Y_i| \leq 2(\sqrt{C_{\overline{K}}} \nnorm{\Sigma^{1/2} \beta_*}_2 + \sqrt{C_{K_{\eps}}} ) \sqrt{\log n}$ with probability at most $4/n^2$. The assertion then follows from a union bound over $i=1,\ldots,n$. $\qed$ 

\subsection{Proof of Proposition \ref{prop:bias_control_truncation}}
We first note that  
\begin{equation*}
\E[\wh{\Sigma} | \mc{R} \cap \mc{L}] = \E[X X^{\T} | \mc{R} \cap \mc{L}], \qquad
\E[\wh{\Sigma}_{Xy} | \mc{R} \cap \mc{L}] = \E[X Y | \mc{R} \cap \mc{L}]
\end{equation*}
by construction of the quantizers $Q_X$, $Q_{X^2}, Q_{Y}$, the construction 
of the estimators $\wh{\Sigma}$ and $\wh{\Sigma}_{Xy}$, and the independence 
of the randomness in the quantization process and the events $\mc{R}$ and $\mc{L}$. It thus 
suffices to bound $\nnorm{\E[XX^{\T} | \mc{R} \cap \mc{L}] - \E[X X^{\T}]}_{\text{op}}$ and 
$\nnorm{\E[X Y| \mc{R} \cap \mc{L}] - \E[X Y]}_{\infty}$. 

Let $X_{\cdot j}$ denote the $j$-th component of $X$, $1 \leq j \leq d$. Moreover, let $(U,V)$ be a pair of random variables following the joint probability measure $\mu$ serving as placeholders for $(X_{\cdot j}, X_{\cdot k})$, 
$1 \leq j,k \leq d$, or $(X_{\cdot j}, Y)$, $1 \leq j \leq d$, and let $R_n$ and $L_n$ be bounds such that 
$\p(\mc{R}^{\textsf{c}}) \leq \epss$ and $\p(\mc{L}^{\textsf{c}}) \leq \delta$ with $\epss > 0$ and $\delta > 0$ to be chosen below. We then have 
\begin{align*}
\E[UV] &= \int_{\mc{R} \cap \mc{L}}  u\cdot v \,d\mu(u,v) + \int_{\mc{R}^{\textsf{c}} \cup \mc{L}^{\textsf{c}}}  u\cdot v \,d\mu(u,v) \\
&\leq \E[U V | \mc{R} \cap \mc{L}]  \p(\mc{R} \cap \mc{L}) + \sqrt{\E[U^2 V^2]} \sqrt{\p(\mc{R}^{\textsf{c}} \cup \mc{L}^{\textsf{c}})}
\end{align*}
by Cauchy-Schwarz. Thus, 
\begin{align*}
\E[U V | \mc{R} \cap \mc{L}]  - \E[UV] &\geq \E[UV] \, \frac{\p(\mc{R}^{\textsf{c}} \cup \mc{L}^{\textsf{c}})}{\p(\mc{R} \cap \mc{L})} - \frac{\sqrt{\E[U^2 V^2]}  \sqrt{\p(\mc{R}^{\textsf{c}} \cup \mc{L}^{\textsf{c}})}}{\p(\mc{R} \cap \mc{L})}\\
&\geq -\sqrt{\E[U^2 V^2]} \left( \frac{\epss + \delta + \sqrt{\epss} + \sqrt{\delta}}{1 - \epss - \delta}  \right) \\
&\geq -C\max\{ \nnorm{\Sigma^{1/2} \beta_*}_2 + \sigma, 1 \} \left( \frac{\epss + \delta + \sqrt{\epss} + \sqrt{\delta}}{1 - \epss - \delta} \right) \\
&\geq -C'\max\{ \nnorm{\Sigma^{1/2} \beta_*}_2 + \sigma, 1 \} (\sqrt{\epss} + \sqrt{\delta})
\end{align*}
where  $C = C_{K,\overline{K}, K_{\eps}} > 0$ depends only on $K = \max_{1 \leq j \leq d} \nnorm{X_{\cdot j}}_{\psi_2}$, $\overline{K} = \nnorm{\Sigma^{-1/2} X}_{\psi_2}$ and 
$K_{\eps} = \nnorm{\eps}_{\psi_2}$ and $C' = 4 C$ as long as $\max\{\delta,\epss \} < 1/4$.

Similarly, $\E[UV] \geq \E[UV | \mc{R} \cap \mc{L}] \p(\mc{R} \cap \mc{L}) - \sqrt{\E[U^2 V^2]} \sqrt{\p(\mc{R}^{\textsf{c}} \cup \mc{L}^{\textsf{c}})}$ and hence $\E[UV | \mc{R} \cap \mc{L}] - \E[UV]$ can be upper bounded analogously. Accordingly, we have
\begin{align*}
&\nnorm{\E[X Y | \mc{L} \cap \mc{R}] - \E[XY]}_{\infty} \leq C \max\{ \nnorm{\Sigma^{1/2} \beta_*}_2 + \sigma, 1 \} (\sqrt{\epss} + \sqrt{\delta}) \\
&\nnorm{\E[X X^{\T} | \mc{L} \cap \mc{R}] - \E[XX^{\T}]}_{\text{op}} \leq 
\nnorm{\E[X X^{\T} | \mc{L} \cap \mc{R}] - \E[XX^{\T}]}_{\textsf{F}} \\
&\qquad \qquad \qquad \qquad \qquad \qquad \quad \;\,\leq d \, C \max\{ \nnorm{\Sigma^{1/2} \beta_*}_2 + \sigma, 1 \} (\sqrt{\epss} + \sqrt{\delta}). 
\end{align*}
Now let $\epss = \delta =  2n^{-q}$ for $q > 0$. This choice is valid by letting $R_n = \sqrt{2 C_K (q+1) \log (n \cdot d)}$  and $L_n = \big(\sqrt{ 2(q+1) C_{\overline{K}}} \nnorm{\Sigma^{1/2} \beta_*}_2 + \sqrt{2 (q+1) C_{K_{\eps}}} \big) \sqrt{\log n}$, using a similar argument as in the proof of Proposition \ref{prop:quantizer_bound}. The assertion of the proposition follows. \qed

\subsection{Asymptotic Relative Efficiency compared to the quantization scheme in \cite{Dirksen2022covariance, Chen2023high}}\label{app:ARE_diagonalSigma}
The referenced work is based on data $\{(\wt{X}_i, \wt{\wt{X}}_i, \wt{Y}_i) \}_{i = 1}^n$, where
$\wt{\wt{X}}_i$ is obtained from $X_i$ through quantization in the same manner as and independently of $\wt{X}_i$, $1 \leq i \leq n$. For ease of demonstration, suppose that $d = 1$, $X \sim N(0,1)$, $\eps \sim N(0,1)$. Then the estimator 
proposed in the present paper and the estimator based on $\{(\wt{X}_i, \wt{\wt{X}}_i, \wt{Y}_i) \}_{i = 1}^n$, respectively, can be expressed as 
\begin{equation*}
\wh{\beta} = \frac{\frac{1}{n} \su \wt{X}_i \wt{Y}_i}{\frac{1}{n} \su \wt{X_i^2}}, \qquad \wt{\beta} = \frac{\frac{1}{2} \left( \frac{1}{n} \su \wt{X}_i \wt{Y}_i + \frac{1}{n} \su \wt{\wt{X}}_i \wt{Y}_i\right)}{\frac{1}{n} \su \wt{X}_i \wt{\wt{X}}_i}. 
\end{equation*}
The asymptotic variance of these two estimators as $n \rightarrow \infty$ can be determined via the 
Delta method \cite[][$\S3$]{vanderVaart1998}. This yields 
\begin{equation*}
n \var(\wh{\beta}) = \frac{\var(U)}{\E[V]^2} + \frac{\E[U]^2}{\E[V]^3} \var(V) - 2 \frac{\cov(U,V) \E[U]}{\E[V]^3} + o(1),   
\end{equation*}
as $n \rightarrow \infty$, where $U = \wt{X}_1 \wt{Y}_1$, $V = \wt{X_1^2}$. Using that 
$\E[V] = 1$, we obtain after collecting terms 
\begin{equation*}
n \var(\wh{\beta})  = \E[U^2] + \E[V^2] \beta_*^2 - 2 \E[X^4] \beta_*^2 = R^2 L^2 + R^2 \beta_*^2 - 6 \beta_*^2  + o(1)
\end{equation*}
Similarly, we find that with $U = \frac{1}{2} (\wt{X}_1 \wt{Y}_1 + \wt{\wt{X}}_1 \wt{Y}_1)$
\begin{equation*}
n \var(\wt{\beta}) = \E[U^2] + \E[V^2] \beta_*^2 - 2 \E[UV] \beta_* = \frac{1}{2} (R^2 L^2 + L^2) + 
R^{4} \beta_*^2 -2R^2 \beta_*^2 + o(1). 
\end{equation*}
We note that for $\beta_* / \sigma$ large, we have $L^2 \approx R^2 \beta_*^2$, thus 
for $R \asymp \sqrt{\log n}$, the term $R^4 \beta_*^2$ is of leading order in both expressions.
This yields $\var(\wt{\beta}) / \var(\wh{\beta}) \approx \frac{3}{2} R^4 \beta_*^2 / R^4 \beta_*^2 = \frac{3}{2}$, i.e., the asymptotic relative efficiency equals $3/2$.

\subsection{Proof of Theorem \ref{theo:asymptotic}}\label{app:theo:asymp}
We start with the following basic decomposition that we will refer to in all three cases discussed below. 
\begin{align}
\wh{\Sigma} \wh{\beta} - \wh{\Sigma}_{Xy} = 0 
  &\Leftrightarrow \; \wh{\Sigma} (\wh{\beta} - \beta_*) + \wh{\Sigma} \beta_* - \wh{\Sigma}_{Xy} = 0 \notag \\
  &\Leftrightarrow \; \wh{\Sigma} \sqrt{n}(\wh{\beta} - \beta_*) =    \sqrt{n}(\wh{\Sigma}_{Xy} -   \wh{\Sigma} \beta_*) \notag \\
  &\Leftrightarrow \;   \sqrt{n}(\wh{\beta} - \beta_*) = \wh{\Sigma}^{-1} \left( \sqrt{n}(\wh{\Sigma}_{Xy} -   \wh{\Sigma} \beta_*) \right) \label{eq:CLT_decomp}
\end{align}  
It holds $\wh{\Sigma} \rightarrow \Sigma$ and thus $\wh{\Sigma}^{-1} \rightarrow \Sigma^{-1} \invcoloneq \Omega$ in probability; this follows from an application of 
Matrix Bernstein's inequality (cf.~Lemma \ref{lem:matrix_bernstein} and Step I.~in the proof of Theorem \ref{theo:nonasymptotic}). Second, we show that
\begin{equation}\label{eq:CLT_term}
  \left( \sqrt{n}(\wh{\Sigma}_{Xy} -   \wh{\Sigma} \beta_*) \right) = \sqrt{n}  \left( \frac{1}{n} \su (\wt{X}_i \wt{Y}_i  - \{ \wt{X}_i \wt{X}_i^{\T} + \Delta_i \} \beta_*) \right)
\end{equation}
satisfies a central limit theorem (CLT). Combining this with the convergence for $\wh{\Sigma}$, the results to be established follow from Slutsky's theorem \cite[][Lemma 2.8]{vanderVaart1998}. 
\vskip1ex
\noindent \underline{Case 1: fixed range, random design}.\\
 The CLT for \eqref{eq:CLT_term} is immediate from the fact 
that the terms inside the sum are i.i.d.~zero mean and bounded (and thus have finite variance). It remains to calculate $\cov(\wt{X}_1 \wt{Y}_1 - (\wt{X}_1 \wt{X}_1^{\T} + \Delta_1)\beta_*)$, which is deferred to the end of this proof.
\vskip2ex
\noindent \underline{Case 2: unbounded range, random design}.\\
We first control the bias of the terms inside \eqref{eq:CLT_term}. We have 
\begin{align}
\E[\{\wt{X}_1 \wt{Y}_1 - (\wt{X}_1 \wt{X}_1^{\T} + \Delta_1) \beta_*\}] 
&= \E[\{\wt{X}_1 \wt{Y}_1 - (\wt{X}_1 \wt{X}_1^{\T} + \Delta_1) \beta_*\} \M{1}_{\{ |X_1| \leq \mc{R}_n \} \cap |Y_1| \leq \mc{L}_n \} }] + \notag \\
&\; + \E[\{\wt{X}_1 \wt{Y}_1 - (\wt{X}_1 \wt{X}_1^{\T} + \Delta_1) \beta_*\} \M{1}_{\{ |X_1| \leq \mc{R}_n \} \cap |Y_1| \leq \mc{L}_n \}^{\textsf{c}}}] \label{eq:case2_unbounded_base}
\end{align}
With $e_j$ denoting the $j$-th canonical basis vector, we obtain the following for the second term in \eqref{eq:case2_unbounded_base} ($1 \leq j \leq d$).
\begin{align}
&|\E[\{\wt{X}_{1j} \wt{Y}_1 - e_j^{\T} (\wt{X}_1 \wt{X}_1^{\T} + \Delta_1) \beta_{*}\} \M{1}_{\{ |X_1| \leq \mc{R}_n \} \cap \{ |Y_1| \leq \mc{L}_n \}^{\textsf{c}}}]| \notag \\
&\leq \E[|\{\wt{X}_{1j} \wt{Y}_1 - e_j^{\T} (\wt{X}_1 \wt{X}_1^{\T} + \Delta_1) \beta_{*}\}|  \M{1}_{\{ |X_1| \leq \mc{R}_n \} \cap \{ |Y_1| \leq \mc{L}_n \}^{\textsf{c}}}] \notag \\
&\leq \E[ \{ |\wt{X}_{1j} \wt{Y}_1| + |e_j^{\T} (\wt{X}_1 \wt{X}_1^{\T} + \Delta_1) \beta_{*}| \}  \M{1}_{\{ |X_1| \leq \mc{R}_n \} \cap \{ |Y_1| \leq \mc{L}_n \}^{\textsf{c}}}] \notag \\
&\leq (L_n R_n + R_n^2 \nnorm{\beta_*}_1)\p(\{ \mc{L}_n \cap \mc{R}_n \}^{{\textsf{c}}}) = O(\log n / n) = o(\tau_n/\sqrt{n}) \label{eq:case2_unbounded_base_loworder_1}
\end{align}
as $n \rightarrow \infty$. For the first term in \eqref{eq:case2_unbounded_base}, we have 
\begin{align}
 &\E[\{\wt{X}_1 \wt{Y}_1 - (\wt{X}_1 \wt{X}_1^{\T} + \Delta_1) \beta_*\} \M{1}_{\{ |X_1| \leq \mc{R}_n \} \cap \{ |Y_1| \leq \mc{L}_n \} }] \notag \\
 &=  \E[ \E[\{\wt{X}_1 \wt{Y}_1 - (\wt{X}_1 \wt{X}_1^{\T} + \Delta_1) \beta_*\} \M{1}_{\{ |X_1| \leq \mc{R}_n \} \cap \{ |Y_1| \leq \mc{L}_n \} } | X_1, Y_1]] \notag \\ 
 &= \E[(X_1 Y_1 - X_1 X_1^{\T} \beta_*) \M{1}_{\{ |X_1| \leq \mc{R}_n \} \cap \{ |Y_1| \leq \mc{L}_n \} }] \notag \\
 &= \{ \E[X_1 Y_1|\mc{L}_n \cap \mc{R}_n] - \E[X_1 X_1^{\T} | \mc{L}_n \cap \mc{R}_n] \beta_* \} \p(\mc{L}_n \cap \mc{R}_n) \notag \\
 &= O(d/n) = o(\tau_n/\sqrt{n}) \label{eq:case2_unbounded_base_loworder2}
\end{align}
according to Proposition \ref{prop:bias_control_truncation}. It remains to show that the variance 
of the terms inside \eqref{eq:CLT_term} is of the order of $\tau_n^2$. For any $j \in \{1,\ldots,d \}$, we have 
\begin{align}
&\var(\wt{X}_{1j} \wt{Y}_1 - (\wt{X}_{1j} \wt{X}_1^{\T} + e_j^{\T} \Delta_1) \beta_* | \mc{L}_n \cap \mc{R}_n) \notag \\
&\leq 4 (\E[\wt{X}_{1j}^2 \wt{Y}_1^2 | \mc{L}_n \cap \mc{R}_n] + \E[\wt{X}_{1j}^2 (\wt{X}_1^{\T} \beta_*)^2 |  \mc{L}_n \cap \mc{R}_n]) + 2 \E[ (e_j^{\T} \Delta_1 \beta_*)^2 | \mc{L}_n \cap \mc{R}_n] \notag \\
&\leq 4 ( R_n^2 L_n^2 + R_n^2  \E[(\wt{X}_1^{\T} \beta_*)^2 | \mc{L}_n \cap \mc{R}_n]) + 2 \beta_{*j}^2 \E[(\Delta_1)_{jj}^2 | \mc{L}_n \cap \mc{R}_n] \notag \\
&\leq  4 ( R_n^2 L_n^2 + R_n^4 \nnorm{\beta_*}_1^2) + 2 \beta_{*j}^2 R_n^4 = O(\tau_n^2). \label{eq:case2_varianceupper} 
\end{align}
To obtain a corresponding lower bound, let us first consider the case $\beta_* = 0$. We then have
\begin{align}
\var(\wt{X}_{1j} \wt{Y}_1 - (\wt{X}_{1j} \wt{X}_1^{\T} + e_j^{\T} \Delta_1) \beta_* | \mc{L}_n \cap \mc{R}_n)
&= \var(\wt{X}_{1j} \wt{Y}_1 |  \mc{L}_n \cap \mc{R}_n) \notag \\
& = \E[\wt{X}_{1j}^2 \wt{Y}_1^2 | \mc{L}_n \cap \mc{R}_n] - \{ \E[\wt{X}_{1j} \wt{\eps}_1 | \mc{L}_n \cap \mc{R}_n] \}^2 \notag \\
&= \E[\wt{X}_{1j}^2 \wt{Y}_1^2 | \mc{L}_n \cap \mc{R}_n] \notag \\ 
&= R_n^2 L_n^2 = \Omega(\tau_n^2), \; 1 \leq j \leq d, \label{eq:case2_variance_lower_1}
\end{align}
where we have used that $Y_1 = \eps_1$ if $\beta_* = 0$. In the opposite case, there exists an index $j_0$ such that $\beta_{*j_0} \neq 0$. We then have 
\begin{align}
&\var(\wt{X}_{1j_0} \wt{Y}_1 - (\wt{X}_{1j_0} \wt{X}_1^{\T} + e_{j_0}^{\T} \Delta_1) \beta_* | \mc{L}_n \cap \mc{R}_n) \notag \\
&\geq \E[\var(\wt{X}_{1j_0} \wt{Y}_1 - (\wt{X}_{1j_0} \wt{X}_1^{\T} + e_{j_0}^{\T} \Delta_1) \beta_*) | \wt{X}_{11}, \ldots, \wt{X}_{1d}, \wt{Y}_1, X_{1j_0}) | \mc{L}_n \cap \mc{R}_n] \notag \\
&=\E[\var(e_{j_0}^{\T} \Delta_1 \beta_* | \wt{X}_{11}, \ldots, \wt{X}_{1d}, \wt{Y}_1, X_{1j_0}) | \mc{L}_n \cap \mc{R}_n] \notag \\
&= \E[\var( \{\wt{X_{1j_0}^2} - \wt{X}_{1j_0}^2 \} \beta_{*j_0} | \wt{X}_{11}, \ldots, \wt{X}_{1d}, \wt{Y}_1, X_{1j_0}) | \mc{L}_n \cap \mc{R}_n] \notag \\
&=\E[\var( \{\wt{X_{1j_0}^2} - R_n^2 \} \beta_{*j_0} | X_{1j_0}) | \mc{L}_n \cap \mc{R}_n] \notag \\
&= \beta_{*j_0}^2 \E[\var( \wt{X_{1j_0}^2}   | X_{1j_0}) | \mc{L}_n \cap \mc{R}_n] = \Omega(R_n^4) = \Omega(\tau_n^2) \label{eq:case2_variance_lower_2}
\end{align}
since conditional on $X_{1j_0}$, the random variable $\wt{X_{1j_0}^2}$ follows a Bernoulli distribution
scaled by $R_n^2$. The variance of that Bernoulli distribution is bounded away from zero unless $\p( |X_{j_0}| \geq R_n)$ tends to one, which is in contradiction to the assumed (unbounded) sub-Gaussianity of $X_{j_0}$. Combining \eqref{eq:case2_unbounded_base}, \eqref{eq:case2_unbounded_base_loworder_1}, 
\eqref{eq:case2_unbounded_base_loworder2}, \eqref{eq:case2_varianceupper}, \eqref{eq:case2_variance_lower_1} \eqref{eq:case2_variance_lower_2} and the fact that 
as $n \rightarrow \infty$, we have $\p(\mc{L}_n \cap \mc{R}_n) \rightarrow 1$ per Proposition \ref{prop:bias_control_truncation}, it follows that 
\begin{equation*}
\frac{\sqrt{n}}{\tau_n}  \left( \frac{1}{n} \su (\wt{X}_i \wt{Y}_i  - \{ \wt{X}_i \wt{X}_i^{\T} + \Delta_i \} \beta_*) \right) \overset{\text{D}}{\rightarrow} N(0, \Gamma) 
\end{equation*}
as $n \rightarrow \infty$ with $\Gamma$ as defined in Theorem \ref{theo:asymptotic}, Part (II). The asymptotic distribution
of $\sqrt{n} (\wh{\beta} - \beta_*)$ is then obtained according to \eqref{eq:CLT_decomp} and the associated remarks.  
\vskip1ex
\noindent \underline{Case 3: fixed range, fixed design}.\\
From \eqref{eq:CLT_decomp}, we have for any unit vector $a \in \R^d$
\begin{align}\label{eq:clt_fixed_decomp}
a^{\T} \mathbb{V}_n^{-1/2} \overline{\Sigma} (\wh{\beta} - \beta_*) 
& = a^{\T} \mathbb{V}_n^{-1/2} \overline{\Sigma} \{ \wh{\Sigma}^{-1} (\wh{\Sigma}_{Xy} - \wh{\Sigma} \beta_*) \} \notag \\
&= a^{\T} \mathbb{V}_n^{-1/2} (\wh{\Sigma}_{Xy} - \wh{\Sigma} \beta_*) +a^{\T} \mathbb{V}_n^{-1/2} \overline{\Sigma} (\wh{\Sigma}^{-1} - \overline{\Sigma}^{-1}) (\wh{\Sigma}_{Xy} - \wh{\Sigma} \beta_*) \\
&=  a^{\T} \mathbb{V}_n^{-1/2} (\wh{\Sigma}_{Xy} - \wh{\Sigma} \beta_*) + o_{\p}(1), \notag
\end{align}
where we note that the second term in \eqref{eq:clt_fixed_decomp} tends to zero in probability since 
$\nnorm{\mathbb{V}_n^{-1/2}}_{\text{op}} = O(\sqrt{n})$, $\nnorm{\wh{\Sigma}_{Xy} - \wh{\Sigma} \beta_*}_2 = O_{\p}(n^{-1/2})$ and $\nnorm{\wh{\Sigma}^{-1} - \overline{\Sigma}^{-1}}_{\text{op}} = o_{\p}(1)$ since $\nnorm{\wh{\Sigma}^{-1} - \Sigma^{-1}}_{\text{op}} = o_{\p}(1)$ (as noted at the beginning of this proof) and $\nnorm{\Sigma^{-1} - \overline{\Sigma}^{-1}}_{\text{op}} = o(1)$ by assumption. Regarding the first term in \eqref{eq:clt_fixed_decomp}, it follows from the Lyapunov CLT \cite[][Thm.~27.3]{Billingsley1995} that  
$$ a^{\T} \mathbb{V}_n^{-1/2} (\wh{\Sigma}_{Xy} - \wh{\Sigma} \beta_*) \overset{\text{D}}{\rightarrow} N(0,1).$$  
By Slutky's Theorem, $a^{\T} \mathbb{V}_n^{-1/2} \overline{\Sigma} (\wh{\beta} - \beta_*) \overset{\text{D}}{\rightarrow} N(0,1)$ as well. The assertion then follows from the Cramer-Wold theorem \cite[][Thm.~29.4]{Billingsley1995}.

\vskip2ex
\noindent {\em Calculation of the covariance matrix} (under Case 1). Letting 
$\wh{\Sigma}_1 = \wt{X}_1 \wt{X}_1^{\T} + \Delta_1$, we have 
\begin{align*}
&\cov(\wt{X}_1 \wt{Y}_1 - (\wt{X}_1 \wt{X}_1^{\T} + \Delta_1)\beta_*) \\
&=\E[\wt{X}_1 \wt{X}_1^{\T} \wt{Y}_1^2] + \E[\wh{\Sigma}_1 \beta_* \beta_*^{\T} \wh{\Sigma}_1 ] - \E[\wh{\Sigma}_1 \beta_* \wt{X}_1^{\T} \wt{Y}_1] - \E[\wt{X}_1 \wt{Y}_1 \wh{\Sigma}_1 \beta_*] \\ 
&= L^2 \{ \Sigma - \text{diag}(\Sigma) + R^2 I_d \} + \Gamma_a - 
\Gamma_b,
\end{align*}
where $\Gamma_a =  \E[\wh{\Sigma}_1 \beta_* \beta_*^{\T} \wh{\Sigma}_1 ]$ and 
$\Gamma_b = \E[\wh{\Sigma}_1 \beta_* \wt{X}_1^{\T} \wt{Y}] + \E[\wt{X}_1 \wt{Y}_1 \wh{\Sigma}_1 \beta_*] $. We note that 
\begin{equation*}
\nnorm{\Gamma_a}_{\infty} = \max_{j,k}|\E[e_j^{\T} \wh{\Sigma}_1 \beta_* \beta_*^{\T} \wh{\Sigma}_1 e_k]| \leq \max_{j,k} \E[|e_j^{\T} \wh{\Sigma}_1 \beta_*| |\beta_*^{\T} \wh{\Sigma}_1 e_k|] \leq R^4 \nnorm{\beta_*}_1^2, 
\end{equation*}
using H\"older's inequality and the observation that $\nnorm{\wh{\Sigma}}_{\infty} \leq R^2$. Moreover, regarding $\Gamma_b$, we have 
\begin{align*}
 \E[\wh{\Sigma}_1 \beta_* \wt{X}_1^{\T} \wt{Y}] &= \E_{X_1, Y_1}[\E[\wh{\Sigma}_1 \beta_* \wt{X}_1^{\T} \wt{Y}_1 | X_1, Y_1]] \\
 &= \E_{X_1, Y_1}[\E[\wh{\Sigma}_1 \beta_* \wt{X}_1^{\T}|X_1] \, \E[\wt{Y}_1|Y_1]]  \\
 &= \E_{X_1, Y_1}[\E[\wh{\Sigma}_1 \beta_* \wt{X}_1^{\T}|X_1] \, Y_1] \\
 &= \E_{X_1}[\E[\wh{\Sigma}_1 \beta_* \wt{X}_1^{\T}|X_1] \, \E[Y_1 | X_1]] \\
 &= \E_{X_1}[\E[\wh{\Sigma}_1 \beta_* \wt{X}_1^{\T}|X_1] \, X_1^{\T} \beta_*]  \\
 &= \E[\wh{\Sigma}_1 \beta_* \wt{X}_1^{\T} \, X_1^{\T} \beta_*] 
\end{align*}
and thus $\max_{j,k}|\E[\wh{\Sigma}_1 \beta_* \wt{X}_1^{\T} \wt{Y}]| \leq 
\max_{j,k}\E[|e_j^{\T} \wh{\Sigma} \beta_*|\wt{X}_1^{\T} e_k||\wt{X}_1^{\T} \beta_*|] \leq R^4 \nnorm{\beta_*}_1^2$, which implies $\nnorm{\Gamma_b}_{\infty} \leq 2 R^4 \nnorm{\beta_*}_1^2$ and in turn $\nnorm{\Gamma_a + \Gamma_b}_{\infty} \leq 3 R^4 \nnorm{\beta_*}_1^2$. 

\subsection{Proof of Theorem \ref{theo:nonasymptotic}}
Consider $\beta_* = \argmin_{\beta \in \R^d} \{ \frac{1}{2} \beta^{\T} \Sigma \beta - \beta^{\T} \Sigma_{Xy} \}$ and $\wh{\beta} = \argmin_{\beta \in \R^d} \{ \frac{1}{2} \beta^{\T} \wh{\Sigma} \beta - \beta^{\T} \wh{\Sigma}_{Xy} \}$ with $\wh{\Sigma}$ and $\wh{\Sigma}_{Xy}$ as defined in 
\eqref{eq:plugin_estimators}. Since $\wh{\beta}$ is a minimizer of its associated optimization problem, we have the following basic inequality: 
\begin{align*}
\frac{1}{2} \wh{\beta}^{\T} \wh{\Sigma} \wh{\beta} \leq \frac{1}{2}
  \beta_*^{\T} \wh{\Sigma} \beta_* + \wh{\Sigma}_{Xy}^{\T}
  (\wh{\beta} - \beta_*)
\end{align*}
In turn, we have 
\begin{equation*}
\frac{1}{2} \wh{\delta}^{\T} \wh{\Sigma} \wh{\delta} \leq
\beta_*^{\T} \wh{\Sigma} \beta_* - \beta_*^{\T} \wh{\Sigma} \wh{\beta} + \wh{\Sigma}_{Xy}^{\T}
  (\wh{\beta} - \beta_*),
\end{equation*}
where $\wh{\delta} = \wh{\beta} - \beta_*$. The right hand side in the
previous display can be re-expressed as
\begin{align}\label{eq:basic_inequality}
\begin{split}
\frac{1}{2} \wh{\delta}^{\T} \wh{\Sigma} \wh{\delta} &\leq
\beta_*^{\T} \wh{\Sigma} (\beta_* - \wh{\beta}) - \wh{\Sigma}_{Xy}^{\T}
(\beta_* - \wh{\beta}) \\
                                                     &=
                                                       -\nscp{\Sigma\beta_*
                                                       -
                                                       \Sigma_{Xy}}{\wh{\beta}
                                                       - \beta_*}
                                                       - \beta_*^{\T}
                                                       (\wh{\Sigma} -
                                                       \Sigma)
                                                       (\wh{\beta} -
                                                       \beta_*) +
                                                       (\wh{\Sigma}_{Xy}
                                                       - \Sigma_{Xy})^{\T}
                                                       (\wh{\beta} -
                                                       \beta_*) \\
  &\leq -\beta_*^{\T} (\wh{\Sigma} - \Sigma) (\wh{\beta} - \beta_*)
    + (\wh{\Sigma}_{Xy}
                                                       - \Sigma_{Xy})^{\T}
                                                       (\wh{\beta} -
                                                       \beta_*) \\
&\leq \nnorm{\beta_*}_2 \nnorm{\wh{\Sigma} - \Sigma}_{\text{op}} \nnorm{\wh{\delta}}_2 + \nnorm{\wh{\Sigma}_{Xy} - \Sigma_{Xy}}_2 \nnorm{\wh{\delta}}_2.    
\end{split}
\end{align}
where in the second inequality we have used that $\Sigma \beta_* = \Sigma_{Xy}$. With the above inequality in place, the left hand side is lower bounded as 
\begin{equation}\label{eq:lambdamin_cond}
\wh{\delta}^{\T} \wh{\Sigma} \wh{\delta} \geq (\lambda_{\min}(\Sigma) - \nnorm{\wh{\Sigma} - \Sigma}_{\text{op}}) \nnorm{\wh{\delta}}_2^2.
\end{equation}
After (I) bounding $\nnorm{\wh{\Sigma} - \Sigma}_{\text{op}} \leq \lambda_{\min}(\Sigma)/2$ by invoking the requirement on the number of samples $n$, and (II) bounding 
$\nnorm{\wh{\Sigma}_{Xy} - \Sigma_{Xy}}_{\infty}$, the assertion of the theorem follows by diving both sides of \eqref{eq:basic_inequality} by $\nnorm{\wh{\delta}}_2$. In the sequel, we complete  
steps (I) and (II) in the situation that $X$ and $Y$ have fixed range $[-R,R]$ and $[-L,L]$, respectively; the variable range case is treated separately at the end of this proof. 
\vskip1ex
\noindent \underline{ {\em Step I}: Bounding $\nnorm{\wh{\Sigma} - \Sigma}_{\text{op}}$}.\\    
We apply the matrix Bernstein inequality as stated in Lemma \ref{lem:matrix_bernstein}. First, we note that by construction (cf.~Eq.~\eqref{eq:plugin_estimators}) 
\begin{equation*}
\wh{\Sigma} = \frac{1}{n} \su \wh{\Sigma}_i, \quad \wh{\Sigma}_i = \wt{X}_i \wt{X}_i^{\T} + \Delta_i, \quad \Delta_i \coloneq \text{diag}(\wt{X_{i1}^2} - \wt{X}_{i1}^2, \ldots, \wt{X_{id}^2} - \wt{X}_{id}^2), \quad 1 \leq i \leq n. 
\end{equation*}
{\em Random Design}: If the $\{X_i \}_{i = 1}^n$ are random, we have $\E[\wh{\Sigma}_i] = \Sigma$, $1 \leq i \leq n$.  We first bound 
\begin{align*}
\nnorm{\wh{\Sigma}_i - \Sigma}_{\text{op}} &\leq \nnorm{\wt{X}_i \wt{X}_i^{\T}}_{\text{op}} + \nnorm{\Delta_i}_{\text{op}} + \nnorm{\Sigma}_{\text{op}} \\
&\leq \nnorm{\wt{X}_i}_2^2 + \max_{1 \leq j \leq d} |\wt{X_{ij}^2} - \wt{X}_{ij}^2| + \nnorm{\Sigma}_{\text{op}} \\
&\leq d R^2 + R^2 + \nnorm{\Sigma}_{\text{op}} \leq (d+1) R^2 + \nnorm{\Sigma}_{\text{op}} \leq (2d + 1) R^2 \invcoloneq \eta, \quad 1 \leq i \leq n, 
\end{align*}
since $\nnorm{\Sigma}_{\text{op}} \leq \tr(\Sigma) \leq dR^2$. 
\vskip1ex
\noindent {\em Fixed Design}. For fixed $\{ X_i \}_{i = 1}^n$, $\Sigma$ in the previous display needs to be replaced 
by $\Sigma_i \coloneq X_i X_i^{\T}$ s.t.~$\nnorm{\Sigma_i}_{\text{op}} \leq \nnorm{X_i}_2^2 \leq d R^2$, yielding the matching bound $\nnorm{\wh{\Sigma}_i - \Sigma_i}_{\text{op}} \leq (2d + 1) R^2$. 
\vskip1ex
\noindent Next, we bound 
\begin{equation*}
\norm{ \frac{1}{n} \su \E[(\wh{\Sigma}_i - \E[\wh{\Sigma}_i])^2]}_{\text{op}} = \norm{ \frac{1}{n} \su \E[\wh{\Sigma}_i^2] - \frac{1}{n} \su \E[\wh{\Sigma}_i]^2}_{\text{op}} \leq \norm{\frac{1}{n} \su \E[\wh{\Sigma}_i^2]}_{\text{op}}  
\end{equation*}
where the inequality follows from the fact that for symmetric positive (semi-)definite matrices 
$A$ and $B$ such that $A - B$ is symmetric positive (semi-)definite, it holds that 
$\nnorm{A - B}_{\text{op}} = \lambda_{\max}(A - B) \leq \lambda_{\max}(A) = \nnorm{A}_{\text{op}}$. 
\noindent
We then have
\begin{align*}
\norm{\frac{1}{n} \E \left[\su \wh{\Sigma}_i^2 \right]}_{\text{op}} &\leq \norm{\E \left[ \frac{1}{n} \su \wt{X}_i  \wt{X}_i^{\T} \nnorm{\wt{X}_i}_2^2 \right]}_{\text{op}} + \norm{\frac{1}{n} \su \E[\Delta_i^2]}_{\text{op}} + \norm{\frac{1}{n} \su \E[\wt{X}_i \wt{X}_i^{\T} \Delta_i]}_{\text{op}} + \\
&\quad + \norm{\frac{1}{n} \su \E[\Delta_i \wt{X}_i \wt{X}_i^{\T} ]}_{\text{op}} \\
&\leq d R^2 \norm{\frac{1}{n} \su \E[\wt{X}_i  \wt{X}_i^{\T}]}_{\text{op}} + R^4 + \frac{2}{n} \su \E[\nnorm{\wt{X}_i}_2^2 \nnorm{\Delta_i}_{\text{op}}] \\
& \leq d R^2 \nnorm{\Sigma + \text{diag}(R^2 - \Sigma_{jj})_{j = 1}^d}_{\text{op}} + R^4 + 2dR^4 \\
 &\leq d R^2 (\nnorm{\Sigma}_{\text{op}} + R^2) + R^4 + 2dR^4 \leq 4 d R^4 + d R^2 \nnorm{\Sigma}_{\text{op}} \invcoloneq \tau^2, \quad 1 \leq i \leq n,  
\end{align*}
where in the second inequality, we have used that the entries of $\Delta_i$ are contained in 
$[-R^2, 0]$, the fact that $\nnorm{\E[A]}_{\text{op}} \leq \E[\nnorm{A}_{\text{op}}]$ for any matrix $A$ (which follows
from Jensen's inequality and the convexity of $\nnorm{\cdot}_{\text{op}}$), and the sub-multiplicativity of $\nnorm{\cdot}_{\text{op}}$. 
\vskip1ex
\noindent Invoking Lemma \ref{lem:matrix_bernstein} with $Q_i = \wh{\Sigma}_i - \Sigma_i$, $\sigma^2 = \tau^2$, $b = \eta$ and 
\begin{align}
\delta &= \sqrt{\frac{8 \tau^2 \log(n) \vee \log(2d)}{n}} + \frac{8\eta\log(n) \vee \log(2d)}{n} \label{eq:Bernstein_bound} \\
       &= \sqrt{\frac{\{ 32 d R^4 + 8dR^2 \nnorm{\Sigma}_{\text{op}} \}  \log(n) \vee \log(2d)}{n}} + \frac{8(2d+1) R^2 \{ \log(n) \vee \log(2d) \}}{n} \notag,
\end{align}
we obtain that $\nnorm{\wh{\Sigma} - \Sigma}_{\text{op}} \leq \delta$ with probability at least 
$1 - 1/n$. 

By requiring that both summands in \eqref{eq:Bernstein_bound} be less than $\lambda_{\min}(\Sigma)/4$, we obtain the following condition on the number of samples 
\begin{equation}\label{eq:cond_lambdamin_samples}
n \geq \left \{ \frac{512 d R^4 + 128 d R^2 \nnorm{\Sigma}_{\text{op}}}{\lambda_{\min}^2(\Sigma)} \vee  32 [(d+1)R^2 + \nnorm{\Sigma}_{\text{op}}] \right \} \{ \log(n) \vee \log(2d) \},    
\end{equation}
in which case \eqref{eq:lambdamin_cond} can be lower bounded by $\nnorm{\wh{\delta}}_2^2 \cdot \lambda_{\min}(\Sigma)/2$ as claimed.
\vskip1ex
\noindent \underline{ {\em Step II}: Bounding $\nnorm{\wh{\Sigma}_{Xy} - \Sigma_{Xy}}_{\infty}$}.\\
Pick $j \in \{1,\ldots,d\}$ arbitrary and consider 
\begin{equation*}
\wh{\Sigma}_{X_j Y} - \Sigma_{X_j Y} = \frac{1}{n} \sum_{i = 1}^n (\wt{X}_{ij} \wt{Y}_i - \E[X_j Y]) 
\end{equation*}
Since $|\wt{X}_{ij} \wt{Y}_i| \leq R L$ and using the fact that $\E[\wt{X}_{ij} \wt{Y}_i] = \E[X_j Y]$, $1 \leq i \leq n$, it follows from Hoeffding's inequality \cite[e.g.][Proposition 2.5]{Wainwright2019} that for any $t > 0$
\begin{equation*}
\p \left( \left|\frac{1}{n} \sum_{i = 1}^n (\wt{X}_{ij} \wt{Y}_i - \E[X_j Y]) \right| > t \right) \leq 2 \exp \left(-\frac{n t^2}{2 R^2 L^2} \right).  
\end{equation*}
Thus, choosing $t = \sqrt{4R^2L^2 \cdot \{ \log(d) \vee \log n \}/n}$ and using the union bound, we obtain that 
\begin{equation}\label{eq:bound_hoeffding}
\nnorm{\wh{\Sigma}_{XY} - \Sigma_{XY}}_{\infty} \leq \sqrt{\frac{4 R^2 L^2  \{ \log(2d) \vee \log n \}}{n}}
\end{equation}
with probability at least $1 - 1/n$. 
\vskip1ex
\noindent \underline{ {\em Step III}: Putting together the pieces}.\\
The result follows by combining \eqref{eq:basic_inequality}, \eqref{eq:lambdamin_cond}, \eqref{eq:Bernstein_bound}, \eqref{eq:cond_lambdamin_samples} and \eqref{eq:bound_hoeffding}. Accordingly, the left hand side of \eqref{eq:basic_inequality} can be lower bounded as 
$\frac{\lambda_{\min}(\Sigma)}{2} \nnorm{\wh{\delta}}_2^2$, and the right hand side of 
\eqref{eq:basic_inequality} can be upper bounded as $\nnorm{\wh{\delta}}_2 \times \text{Expression in \eqref{eq:Bernstein_bound}} \times \nnorm{\beta_*}_2 + \nnorm{\wh{\delta}}_2 \times \sqrt{d}\cdot\text{Expression in \eqref{eq:bound_hoeffding}}$. 
\vskip1ex 
\noindent \underline{ {\em Step IV}: Extension to unbounded $X$ and $Y$}.\\
To conclude, we discuss how the analysis above is extended to the setting in $\S$\ref{subsec:nofixed_range}. Starting from \eqref{eq:basic_inequality}, we bound 
\begin{align*}
&\nnorm{\wh{\Sigma} - \Sigma}_{\text{op}} \leq
\nnorm{\wh{\Sigma} - \E[\wh{\Sigma}] +  \E[\wh{\Sigma}] - \Sigma}_{\text{op}} \leq 
\nnorm{\wh{\Sigma} - \E[\wh{\Sigma}]}_{\text{op}} + \nnorm{\E[\wh{\Sigma}] - \Sigma}_{\text{op}} \\
&\nnorm{\wh{\Sigma}_{Xy} - \E[\wh{\Sigma}_{Xy}] + \E[\wh{\Sigma}_{Xy}]  - \Sigma_{Xy}}_{\infty} \leq \nnorm{\wh{\Sigma}_{Xy} - \E[\wh{\Sigma}_{Xy}]}_{\infty} + \nnorm{\E[\wh{\Sigma}_{Xy}]  - \Sigma_{Xy}}_{\infty}, 
\end{align*}
where the expectations are conditional on the truncation events $\mc{R}$ and $\mc{L}$ in Proposition 
\ref{prop:bias_control_truncation}, invoked with the choice $q = 3$. With this choice, the terms 
$\nnorm{\E[\wh{\Sigma}] - \Sigma}_{\text{op}} = O(d/n^{-3/2})$ and
$\nnorm{\E[\wh{\Sigma}_{Xy}]  - \Sigma_{Xy}}_{\infty} = O(n^{-3/2})$ are of lower order relative 
to $\nnorm{\wh{\Sigma} - \E[\wh{\Sigma}]}_{\text{op}}$ and \linebreak$\nnorm{\wh{\Sigma}_{Xy} - \E[\wh{\Sigma}_{Xy}]}_{\infty}$, respectively, which can be controlled according to Steps II and III above, setting $R = R_n$ and $L = L_n$ chosen per Proposition \ref{prop:bias_control_truncation}. Finally, the condition on the number of samples \eqref{eq:cond_lambdamin_samples} needs to be modified to account for the additional term $\nnorm{\E[\wh{\Sigma}] - \Sigma}_{\text{op}}$. Specifically, in addition to \eqref{eq:cond_lambdamin_samples}, we require that $n \geq C' \{ \lambda_{\min}^{-1}(\Sigma) d \nnorm{\Sigma^{1/2} \beta_*}_2 + \sigma, 1 \}\}^{2/3}$, where $C' = (4C)^{2/3}$ with 
$C$ as in Proposition \ref{prop:bias_control_truncation}. This yield the lower bound 
$\lambda_{\min}(\wh{\Sigma}) \geq \frac{1}{4} \lambda_{\min}(\Sigma)$ (as opposed to $\frac{1}{2} \lambda_{\min}(\Sigma)$), which increases the final error bound by a factor of $2$ compared to the fixed range case. Additionally, we incur lower order terms on the right hand side of \eqref{eq:basic_inequality} that reflect $\nnorm{\E[\wh{\Sigma}] - \Sigma}_{\text{op}}$ and $\nnorm{\E[\wh{\Sigma}_{Xy}]  - \Sigma_{Xy}}_{\infty}$, and the probability statement needs to be adjusted to account for the events $\mc{R}$ and $\mc{L}$. \qed 

\subsection*{Proof of Theorem \ref{theo:lowerbound}}
By direct computation, we obtain that 
\begin{align*}
\frac{1}{n} \left( -\frac{d^2}{d \beta^2} \log \textsf{L}(\beta) \dev{\beta}{\beta_*} \right) &= -\frac{\su c_i}{n} \left( \frac{\pi(\beta_*) \cdot \ddot{\pi}(\beta_*) - \{ \dot{\pi}(\beta_*) \}^2}{\{ \pi(\beta_*) \}^2} \right) \\
&\quad -\frac{n - \su c_i}{n} \left(\frac{(1 - \pi(\beta_*)) (- \ddot{\pi}(\beta_*)) - \{ \dot{\pi}(\beta_*)\}^2}{(1 - \pi(\beta_*))^2} \right)
\end{align*}
Taking expectations so that  $\E[\su c_i / n] = \pi(\beta_*)$, we obtain
\begin{equation}\label{eq:info_app}
\frac{1}{n} \E\left[-\frac{d^2}{d \beta^2} \log \textsf{L}(\beta) \dev{\beta}{\beta_*} \right] = \{ \dot{\pi}(\beta_*) \}^2 \left(\frac{1}{\pi(\beta_*)} + \frac{1}{1 - \pi(\beta_*)}\right),     
\end{equation}
i.e., the first equality in the theorem. Since $0 < \underline{c}(\sigma, \beta_*)< \pi(\beta_*) < \overline{c}(\sigma, \beta_*) < 1$ for constants $ \underline{c}$ and $\overline{c}$ depending
only on $\beta_*$ and $\sigma$, it suffices to upper bound $|\dot{\pi}(\beta_*)|$, which is done in the sequel. Letting $\Phi$ and $\phi$ denote the cumulative distribution and density function of the standard Normal distribution, we have the following for the inner double integral over $x$ and $y$ in \eqref{eq:collisionprob}: 
\begin{align}\label{eq:collisionprob_refactor}
& \frac{1}{2\pi}\int_{0}^{\infty}\int_{0}^{\infty}
  |V|^{-1/2} \exp\left(-\frac{1}{2}\begin{pmatrix}
      x - z_1 \\
      y - z_2 
      \end{pmatrix} V^{-1} \begin{pmatrix}
      x - z_1 \\
      y - z_2 
      \end{pmatrix} 
  \right) \;dx \, dy \notag \\
&=\frac{1}{2\pi} \int_{-z_2}^{\infty}\int_{-z_1}^{\infty}
  |V|^{-1/2} \exp\left(-\frac{1}{2}\begin{pmatrix}
      x  \\
      y  
      \end{pmatrix} V^{-1} \begin{pmatrix}
      x  \\
      y 
      \end{pmatrix} 
  \right) \;dx \, dy \notag \\
&= \int_{-z_1}^{\infty} \frac{1}{\sqrt{2\pi}} \exp \left(-\frac{x^2}{2}
                          \right)      \left( \int_{-z_2}^{\infty} \frac{1}{\sqrt{2\pi}}
                          \frac{1}{\sigma} \left(-\frac{(y - \beta
                          x)^2}{2 \sigma^2}
                          \right) \; dy \right) \;dx \notag \\
&= \int_{-z_1}^{\infty} \frac{1}{\sqrt{2\pi}} \exp \left(-\frac{x^2}{2}
                          \right)  \left(1 - \Phi\left(\frac{-z_2 -
                                                          \beta
                                                          x}{\sigma}
                                                          \right)
                                                          \right) \; dx, 
\end{align}
by factorizing the joint density of $(X,Y)$ as the density of $X$ times the density of $Y|X \sim N(X \beta, \sigma^2)$. Differentiating \eqref{eq:collisionprob_refactor} w.r.t.~$\beta$ under the integral sign and  evaluating the result at $\beta = \beta_*$, we obtain 
\begin{align}
&\left( \frac{d}{d\beta} \int_{-z_1}^{\infty} \frac{1}{\sqrt{2\pi}} \exp \left(-\frac{x^2}{2}
                          \right)  \left(1 - \Phi\left(\frac{-z_2 -
                                                          \beta
                                                          x}{\sigma}
                                                          \right)
                                                          \right) \; dx  \right) \dev{\beta}{\beta_*}\notag \\
&=\int_{-z_1}^{\infty} \frac{1}{\sqrt{2 \pi}}  \exp \left(-\frac{x^2}{2}
                          \right) \frac{x}{\sigma} \exp\left(-\frac{(z_2 + \beta_* x)^2}{2
  \sigma^2} \right) \; dx. \label{eq:collisionprob_refactor_diff}
\end{align}
By completing squares, \eqref{eq:collisionprob_refactor_diff} yields the integral
\begin{align*}
\frac{1}{2\pi \sigma} \exp\left(-\frac{z_2^2}{2 (\beta_*^2 + \sigma^2)}
\right) \int_{-z_1}^{\infty} x \, \exp \left(-\frac{1}{2}
  \frac{\left(x + \frac{z_2 \beta_*}{\sigma^2 + \beta_*^2}
  \right)^2}{\frac{\sigma^2}{\sigma^2 + \beta_*^2}} \right) \;dx.  
\end{align*}
Evaluating the latter integral and collecting terms, we obtain that
\begin{align}\label{eq:low_intm_1}
&\frac{1}{2\pi} \frac{1}{\sigma} \exp\left(-\frac{z_2^2}{2 (\beta_*^2 +
  \sigma^2)} \right)
  \Bigg\{ \frac{\sigma^2}{\sigma^2 + \beta_*^2} \exp
  \left(-\frac{(z_1 - z_2 \frac{\beta_*}{\sigma^2 +
  \beta_*^2})^2}{2\frac{\sigma^2}{\beta_*^2 + \sigma^2}} \right) \notag \\
&\qquad \qquad \qquad \qquad \qquad \qquad -z_2 \frac{\beta_*}{\sigma^2 + \beta_*^2}
                                                                  \sqrt{\frac{\sigma^2}{\beta_*^2
                                                                  +
                                                                  \sigma^2}}
                                                                  \sqrt{2\pi}       
                                                                  \Phi\left(\frac{z_1
                                                                  +
                                                                  z_2
                                                                  \frac{\beta_*}{\sigma^2
                                                                  +
                                                                  \beta_*^2}}{\sqrt{\frac{\sigma^2}{\sigma^2
                                                                  + \beta_*^2}}}
                                                                  \right)
                                                                  \Bigg \}.   
\end{align}
We now turn to the outer double integral in \eqref{eq:collisionprob}. For the first summand in \eqref{eq:low_intm_1}, we obtain that
\begin{align*}
&\int_{-L}^L \int_{-R}^{R} \frac{1}{2\pi} \frac{1}{\sigma} \exp\left(-\frac{z_2^2}{2 (\beta_*^2 +
  \sigma^2)} \right)  \exp
  \left(-\frac{(z_1 - z_2 \frac{\beta_*}{\sigma^2 +
  \beta_*^2})^2}{2\frac{\sigma^2}{\beta_*^2 + \sigma^2}} \right) dz_1 dz_2
  \times \frac{\sigma^2}{\beta_*^2 + \sigma^2} \\
&=  \frac{\sigma^2}{\beta_*^2 + \sigma^2}
                                               \text{Prob}_{\Sigma}(-L,L,-R,R)
                                               \leq 1, 
\end{align*}
where $\text{Prob}_{V}(-L,L,-R,R) = \p(G_1 \in [-R,R], G_2 \in
[-L,L])$ with $(G_1, G_2)$ being bivariate Gaussian with covariance
matrix $V$ as defined in \eqref{eq:collisionprob}. To see this, observe that
\begin{equation*}
G_2 \sim N(0, \beta_*^2 + \sigma^2), \qquad G_1|G_2 = g_2 \sim N \left(
\frac{\beta_*}{\sigma^2 + \beta_*^2} g_2,  \frac{\sigma^2}{\sigma^2 + \beta_*^2}\right), \;  \; g_2 \in \R.
\end{equation*}  
Turning to the second term in \eqref{eq:low_intm_1}, we have the following:
\begin{align}\label{eq:low_intm_2}
\int_{-L}^{L} \left(-\frac{z_2}{\sigma^2 + \beta_*^2} \right) \exp
  \left(-\frac{z_2^2}{2 (\sigma^2 + \beta_*^2)} \right) \Phi \left(\frac{z_1 + z_2
  \frac{\beta_*}{\sigma^2 + \beta_*^2}}{\sqrt{\frac{\sigma^2}{\sigma^2 + \beta_*^2}}}\right)
  \; dz_2 \times C(\sigma, \beta_*). 
\end{align}
for a constant $C(\sigma, \beta_*) > 0$. Using integration by parts to evaluate the integral in \eqref{eq:low_intm_2},  we obtain that
{\small \begin{align}\label{eq:low_intm_3}
 &\exp \left(-\frac{z_2^2}{2 (\beta_*^2 + \sigma^2)} \right)
  \Phi \left(\frac{z_1 + z_2
  \frac{\beta_*}{\sigma^2 + \beta_*^2}}{\sqrt{\frac{\sigma^2}{\sigma^2 +
      \beta_*^2}}}\right) \bou{-L}{L} - C'(\sigma, \beta_*) \int_{-L}^{L}  \exp
  \left(-\frac{z_2^2}{2 (\beta_*^2 + \sigma^2)} \right)
   \phi\left(\frac{z_1 + \frac{z_2 \beta_*}{\sigma^2 +
   \beta_*^2}}{\sqrt{\frac{\sigma^2}{\sigma^2 + \beta_*^2}}} \right) \; dz_2
  \notag \\
&= \exp\left(-\frac{L^2}{2(\beta_*^2 + \sigma^2)} \right)  \left\{  \Phi \left(\frac{z_1 + L
  \frac{\beta_*}{\sigma^2 + \beta_*^2}}{\sqrt{\frac{\sigma^2}{\sigma^2 +
      \beta_*^2}}}\right) - \Phi \left(\frac{z_1 - L
  \frac{\beta_*}{\sigma^2 + \beta_*^2}}{\sqrt{\frac{\sigma^2}{\sigma^2 +
       \beta_*^2}}}\right) \right\} \notag \\
&\quad-C''(\sigma, \beta_*) \int_{-L}^L  \exp
  \left(-\frac{z_2^2}{2 (\beta_*^2 + \sigma^2)} \right) \,
                                       \exp\left(-\frac{\left( z_1 +
                                       \frac{\beta_*}{\sigma^2 +
                                       \beta_*^2} z_2 \right)^2}{2
                                       \frac{\sigma^2}{\sigma^2 + \beta_*^2}} \right)\;dz_2  
\end{align}} 
for constants $C'(\sigma, \beta_*), C''(\sigma, \beta_*) > 0$.
For the second term in \eqref{eq:low_intm_3}, we complete squares to obtain
\begin{align*}
&  C''(\sigma, \beta_*) \exp\left(-\frac{z_1^2 \frac{\sigma^2}{\sigma^2 +
      \beta_*^2}}{2
      \frac{\sigma^2}{\sigma^2 + \beta_*^2}} \right) \int_{-L}^L
  \exp \left(-\frac{\left( z_1 \sqrt{\frac{\beta_*^2}{\sigma^2 + \beta_*^2}} +
      z_2 \frac{1}{\sqrt{\beta_*^2 + \sigma^2}} \right)^2}{2
                 \frac{\sigma^2}{\sigma^2 + \beta_*^2}}\right) \; dz_2 \\
&= C'''(\sigma, \beta_*)  \exp\left(-\frac{z_1^2}{2} \right) \left(\Phi\left(\frac{z_1
                                                                          \sqrt{\frac{\beta_*^2}{\sigma^2
                                                                          +
                                                                          \beta_*^2}}
                                                                          +
                                                                          \frac{L}{
                                                                          \sqrt{\beta_*^2
                                                                          +
                                                                          \sigma^2}}}{\sqrt{\frac{\sigma^2}{\sigma^2
                                                                          +
                                                                          \beta_*^2}}}
                                                                          \right)
                                                                          -
                                                                          \Phi\left(\frac{z_1
                                                                          \sqrt{\frac{\beta_*^2}{\sigma^2
                                                                          +
                                                                          \beta_*^2}}
                                                                          -
                                                                          \frac{L}{
                                                                          \sqrt{\beta_*^2
                                                                          +
                                                                          \sigma^2}}}{\sqrt{\frac{\sigma^2}{\sigma^2
                                                                          +
                                                                          \beta_*^2}}}
                                                                          \right)
                                                                          \right)   
\end{align*}
for a constant $C'''(\sigma, \beta_*) > 0$. Integrating the first term in \eqref{eq:low_intm_3} over $z_1$, we obtain that
\begin{align}\label{eq:low_intm_4}
&\exp \left(-\frac{L^2}{2(\beta_*^2 + \sigma^2)} \right) \left| \int_{-R}^R  \left\{  \Phi \left(\frac{z_1 + L
  \frac{\beta_*}{\sigma^2 + \beta_*^2}}{\sqrt{\frac{\sigma^2}{\sigma^2 +
      \beta_*^2}}}\right) - \Phi \left(\frac{z_1 - L
  \frac{\beta_*}{\sigma^2 + \beta_*^2}}{\sqrt{\frac{\sigma^2}{\sigma^2 +
                 \beta_*^2}}}\right) \right\}  \; dz_1 \right| \notag \\
&\leq 2 R / \log n \leq 2C, 
\end{align}  
where we have used that the integrand is uniformly bounded by $1$ and the fact that $L \geq \sqrt{2(\beta_*^2 + \sigma^2) (\log \log n)^{1/2}}$ and $R \leq C \sqrt{\log
  n}$ by assumption. Integrating over the second term in \eqref{eq:low_intm_3}, we obtain that
\begin{align}\label{eq:low_intm_5}
&\int_{-R}^R \exp\left(-\frac{z_1^2}{2} \right)  \, \left(\Phi\left(\frac{z_1
                                                                          \sqrt{\frac{\beta_*^2}{\sigma^2
                                                                          +
                                                                          \beta_*^2}}
                                                                          +
                                                                          \frac{L}{
                                                                          \sqrt{\beta_*^2
                                                                          +
                                                                          \sigma^2}}}{\sqrt{\frac{\sigma^2}{\sigma^2
                                                                          +
                                                                          \beta_*^2}}}
                                                                          \right)
                                                                          -
                                                                          \Phi\left(\frac{z_1
                                                                          \sqrt{\frac{\beta_*^2}{\sigma^2
                                                                          +
                                                                          \beta_*^2}}
                                                                          -
                                                                          \frac{L}{
                                                                          \sqrt{\beta_*^2
                                                                          +
                                                                          \sigma^2}}}{\sqrt{\frac{\sigma^2}{\sigma^2
                                                                          +
                                                                          \beta_*^2}}}
                                                                          \right)
                                                                          \right)
                                                                          \;dz_1 \notag \\          
                                                                          &\leq  \int_{-R}^R \exp\left( -\frac{z_1^2}{2} \right) \leq \sqrt{2 \pi}                                                                         
\end{align}
Combining \eqref{eq:collisionprob} and \eqref{eq:collisionprob_refactor} through \eqref{eq:low_intm_5}, we have established that $|\dot{\pi}(\beta_*)| \leq \frac{1}{L \cdot R} C(\beta_*, \sigma)$. Inserting this bound into \eqref{eq:info_app} yields the assertion. \qed

\subsection{Proof of Theorem \ref{theo:lasso}}
Note that according to the first-order optimality condition of the optimization
problem defining the estimator \eqref{eq:estimator_ell1}, we have that 
\begin{equation*}
\nscp{\wh{\Sigma} \wh{\beta} - \wh{\Sigma}_{Xy} + \lambda \nabla
\nnorm{\wh{\beta}}_1}{\beta_* - \wh{\beta}} \geq 0. 
\end{equation*}
where $\nabla \nnorm{\wh{\beta}}_1$ is an element of the
sub-differential of the map $\beta \mapsto \nnorm{\beta}_{1}$ at $\beta =
\wh{\beta}$ \cite[cf.][Proposition B.24 (f)]{Bertsekas1999}.
Let $\wh{\delta} = \wh{\beta} - \beta_*$. Consequently, we have that
\begin{align}
0 &\geq \nscp{\wh{\Sigma} \wh{\beta} - \wh{\Sigma}_{Xy} + \lambda \nabla
    \nnorm{\wh{\beta}}_1}{\wh{\delta}} \notag \\
  &= \wh{\delta}^{\T} \wh{\Sigma} \wh{\delta} + \nscp{\wh{\Sigma}
    \beta_* - \wh{\Sigma}_{Xy}}{\wh{\delta}} + \lambda \nscp{\nabla
    \nnorm{\wh{\beta}}_1}{\wh{\delta}} \notag \\
    &\geq  \wh{\delta}^{\T} \wh{\Sigma} \wh{\delta} - \lambda_0 \nnorm{\wh{\delta}}_1 + \lambda \nscp{\nabla
    \nnorm{\wh{\beta}}_1}{\wh{\delta}} \label{eq:firstorder},
\end{align}
where the last inequality is conditional on the event in Lemma \ref{lem:grad:sup} (cf.~Appendix \ref{app:add_tech}). Below, it will shown that for all $v \in \mc{C}$ with 
\begin{align}\label{eq:tangentcone}
\begin{split}
\mc{C} \coloneq &\mc{C}_0 \cap \left\{v \in \R^d: \nnorm{v}_2 > \wt{C} B  \frac{\nnorm{\Sigma}_{\text{op}}^{1/2} + R}{\{ \lambda_{\min}(\Sigma) \}^{1/2}} \sqrt{\frac{s \log d}{n}}  \right\} \\
& \mc{C}_0 \coloneq \{v \in \R^d: v = u - \beta_*, \; u \in \R^d, \; \nnorm{u}_1 \leq \sqrt{s} B \},
\end{split}
\end{align}
where $\wt{C} > 0$ is a constant depending only on $\overline{K}$ as specified in \eqref{eq:normv_lower}, it holds that 
\begin{equation}\label{eq:RSC}
v^{\T} \wh{\Sigma} v \geq \underbrace{\frac{1}{4} \lambda_{\min}(\Sigma)}_{\invcoloneq \kappa} \nnorm{v}_2^2  - \underbrace{\wt{C}^2 (\nnorm{\Sigma}_{\text{op}}^{1/2} + R) \alpha_0}_{\invcoloneq \tau} \frac{\log d}{n} \nnorm{v}_1^2, \quad \alpha_0 \coloneq \frac{2 R^2}{\lambda_{\min}(\Sigma)},
\end{equation}
Note that $\wh{\delta} \in \mc{C}_0$. If $\wh{\delta} \notin \mc{C}$, the bound (i) follows immediately. Otherwise, \eqref{eq:firstorder} and \eqref{eq:RSC} imply that 
\begin{align}
0 &\geq  \kappa \nnorm{\wh{\delta}}_2^2 - \tau \frac{\log d}{n} \nnorm{\wh{\delta}}_1^2 - \lambda_0 \nnorm{\wh{\delta}}_1 + \lambda \nscp{\nabla
    \nnorm{\wh{\beta}}_1}{\wh{\delta}} \notag \\
     &\overset{\text{(i)}}{\geq}   \kappa \nnorm{\wh{\delta}}_2^2 - \tau  \frac{\log d}{n}
    \nnorm{\wh{\delta}}_1 (2 \sqrt{s} B) - \lambda_0 \nnorm{\wh{\delta}}_1 + \lambda \nscp{\nabla
      \nnorm{\wh{\beta}}_1}{\wh{\delta}} \notag \\
    &\overset{\text{(ii)}}{\geq} \kappa \nnorm{\wh{\delta}}_2^2 - \tau  \frac{\log d}{n}
    \nnorm{\wh{\delta}}_1 (2 \sqrt{s} B) - \lambda_0
      \nnorm{\wh{\delta}}_1 + \lambda (\nnorm{\wh{\beta}}_1 -
      \nnorm{\beta_*}_1) \notag \\
     &\overset{\text{(iii)}}{\geq} \kappa \nnorm{\wh{\delta}}_2^2 - \tau  \frac{\log d}{n}
    \nnorm{\wh{\delta}}_1 (2 \sqrt{s} B) - \frac{\lambda}{\gamma}
      \nnorm{\wh{\delta}}_1 + \lambda (\nnorm{\wh{\beta}}_1 -
       \nnorm{\beta_*}_1) \notag \\
       &\overset{\text{(iv)}}{\geq} \kappa \nnorm{\wh{\delta}}_2^2 -\frac{\lambda}{2}
         \nnorm{\wh{\delta}}_1 + \lambda (\nnorm{\wh{\beta}}_1 -
       \nnorm{\beta_*}_1) \label{eq:firstorder:RSC},
\end{align}
where (i) follows from the triangle inequality and the fact that 
$\nnorm{\wh{\beta}}_1 \leq \sqrt{s} B$ and $\nnorm{\beta_*}_1 \leq \sqrt{s} B$, 
(ii) is a consequence of the subgradient inequality of the form 
$f(y) \geq f(x) + \nscp{g_x}{y - x}$ for convex $f$ with
subgradient $g_x$ at $x$ \cite[cf.][Eq.~(B.21)]{Bertsekas1999}, in (iii) we have used that $\lambda \geq \gamma \lambda_0$ for $\gamma >
2$, and (iv) follows from 
\begin{equation}\label{eq:condn_tau}
\frac{\tau \log d \, (2 \sqrt{s} B)}{n} \leq \left(\frac{1}{2} - \frac{1}{\gamma} \right)
\lambda, 
\end{equation}
which is implied by the lower bound on $n$ in the statement of the theorem. Inequality \eqref{eq:firstorder:RSC} yields
\begin{equation*}
  \frac{\lambda}{2} \nnorm{\wh{\delta}}_1 \geq \lambda
  (\nnorm{\wh{\beta}}_1 - \nnorm{\beta_*}_1).
\end{equation*}
Let $S$ denote the support of $\beta_*$, and let the subscript $S$
denote the corresponding sub-vector. The preceding inequality 
implies that $\frac{3}{2} \nnorm{\wh{\delta}_S}_1 \geq \frac{1}{2}
\nnorm{\wh{\delta}_{S^c}}_1$
and thus $\nnorm{\wh{\delta}_{S^c}}_1 \leq 3
\nnorm{\wh{\delta}_S}_1$. Back-substitution into \eqref{eq:firstorder:RSC} then yields
\begin{equation*}
\nnorm{\wh{\delta}}_2^2 \leq \kappa^{-1} \lambda \frac{3}{2}
\nnorm{\wh{\delta}}_1 \leq \kappa^{-1} 6 \lambda \nnorm{\wh{\delta}_S}_1 \leq \kappa^{-1} 6
\sqrt{s} \lambda \nnorm{\wh{\delta}}_2.
\end{equation*}
The final bound is then obtained by dividing both sides by
$\nnorm{\wh{\delta}}_2$.

It remains to establish property \eqref{eq:RSC}. We decompose
$\wh{\Sigma} = \wh{\Sigma}_0 - \wh{D}$ with 
$\wh{\Sigma}_0 = \frac{1}{n} \su \wt{X}_i \wt{X}_i^{\T}$
and $\wh{D} = \frac{1}{n} \su \text{diag}(R^2 - \wt{X_{i1}^2}, \ldots, R^2 - \wt{X_{id}^2})$. Accordingly, we let $\Sigma_0 = \E[\wh{\Sigma}_0]$ and $D = \E[\wh{D}]$.  Consequently, we have for any $v \in \mc{C}$ with $\mc{C}$ in \eqref{eq:tangentcone} 
\begin{align}\label{eq:RE_basic_1}
  v^{\T} \wh{\Sigma} v &= v^{\T} \wh{\Sigma}_0 v - v^{\T} \wh{D} v \notag \\
                       &\geq \nnorm{ \wh{\Sigma}_0^{1/2} v}_2^2 -
                         v^{\T} D v -  \big|\{ v^{\T} (D -
                         \wh{D})  v \} \big|,
\end{align}
where, with some abuse of notation, we write $\wh{\Sigma}_0^{1/2} = \wt{\M{X}}/\sqrt{n} \coloneq (\wt{X}_{ij}/\sqrt{n})_{1 \leq i \leq n, \, 1 \leq j \leq d}$ for the scaled matrix
of quantized $X$s.  Similarly, we have 
\begin{equation}\label{eq:RE_basic_2}
\nnorm{\wh{\Sigma}_0^{1/2} v}_2 \geq \nnorm{\Sigma_0^{1/2} v}_2 - \underbrace{\left|\nnorm{\wh{\Sigma}_0^{1/2} v}_2 - \nnorm{\Sigma_0^{1/2}
    v}_2 \right|}_{\invcoloneq \Gamma(v)}. 
\end{equation}
We first bound $\Gamma(v)$. We have 
\begin{equation}\label{eq:Gamma_decompose}
\Gamma(v) \leq  \nnorm{v}_1  \cdot \left\{ \sup_{v \in\mathbb{B}_1(1)}  \left|\nnorm{\wh{\Sigma}_0^{1/2} v}_2 - \nnorm{\Sigma_0^{1/2}
    v}_2 \right| \right\},
\end{equation}
where $\mathbb{B}_1(1) = \{ v \in \R^d: \; \nnorm{v}_1 \leq 1\}$. We have 
\begin{align*}
\sup_{v \in\mathbb{B}_1(1)}  \left|\nnorm{\wh{\Sigma}_0^{1/2} v}_2 - \nnorm{\Sigma_0^{1/2}
    v}_2 \right|  &= \sup_{v \in\mathbb{B}_1(1)}  \left|\norm{\frac{\M{X}}{\sqrt{n}} \Sigma_0^{-1/2} \Sigma_0^{1/2} v}_2 - \nnorm{\Sigma_0^{1/2}
    v}_2 \right| \\
    &= \sup_{v \in\mathbb{B}_1(1)}  \left|\norm{\frac{\wt{\M{X}}}{\sqrt{n}}  \Sigma_0^{1/2} v}_2 - \nnorm{\Sigma_0^{1/2}
    v}_2 \right| \\
    &= \sup_{v \in \Sigma_0^{1/2} \mathbb{B}_1(1)}  \left|\norm{\frac{\wt{\M{X}}}{\sqrt{n}} v}_2 - \nnorm{v}_2 \right| \invcoloneq \mc{T}
\end{align*}
Note that the matrix $\wt{\M{X}} = \M{X} \Sigma_0^{-1/2}$ has i.i.d.~isotropic 
and sub-Gaussian rows. It hence follows from Lemma \ref{lem:matrix_deviation} 
with $T = \mathbb{B}_1(1)$ that 
\begin{equation}\label{eq:matrix_deviation}
\p \left(\mc{T} > C \overline{K}^2 \nnorm{\Sigma_0}_{\text{op}}^{1/2} \sqrt{\frac{\log d}{n}}  \right) \leq \frac{2}{d},
\end{equation}
where we have used that the Gaussian width of $\Sigma_0^{1/2} \mathbb{B}_1(1)$ is bounded by
$C \nnorm{\Sigma_0}_{\text{op}}^{1/2} \sqrt{\log d}$ \cite[cf.][$\S$7.5]{Vershynin2018}. Conditional on the event in \eqref{eq:matrix_deviation}, we have 
\begin{equation}\label{eq:Gammabar}
\Gamma(v) \leq \nnorm{v}_1 \cdot C \overline{K}^2 \{\nnorm{\Sigma}_{\text{op}}^{1/2} + R \}  \sqrt{\log(d)/n} \invcoloneq \overline{\Gamma}(v)
\end{equation}
since $\nnorm{\Sigma_0}_{\text{op}}^{1/2} \leq \nnorm{\Sigma}_{\text{op}}^{1/2} + \nnorm{D}_{\text{op}}^{1/2} \leq \nnorm{\Sigma}_{\text{op}}^{1/2}  + R$.
Inserting this bound for the second term on the right hand side of \eqref{eq:RE_basic_2} yields 
\begin{equation}\label{eq:Sigma0v_lower}
\nnorm{\wh{\Sigma}_0^{1/2} v}_2 \geq \nnorm{\Sigma_0^{1/2} v}_2  - \overline{\Gamma}(v)
\end{equation}
Note that since $\Sigma_0 = \Sigma + D$ with $D$ positive semidefinite, we have 
$\nnorm{\Sigma_0^{1/2} v}_2 \geq \{ \lambda_{\min}(\Sigma) \}^{1/2} \nnorm{v}_2$. Now observe that since $v \in \mc{C}$ defined in \eqref{eq:tangentcone}, we have 
\begin{equation}\label{eq:normv_lower}
\nnorm{v}_2 > \wt{C} B \frac{\nnorm{\Sigma}_{\text{op}}^{1/2} + R}{\{ \lambda_{\min}(\Sigma) \}^{1/2}} \sqrt{s \log(d) /n}, \quad \wt{C} \coloneq 2 C \overline{K}^2. 
\end{equation}
At the same time, $\nnorm{v}_1 \leq 2B \sqrt{s}$ for all $v \in \mc{C} \subseteq \mc{C}_0$ . Combining this fact, \eqref{eq:Gammabar} and \eqref{eq:normv_lower}, it follows that 
the right hand side of \eqref{eq:Sigma0v_lower} is positive, and the inequality is 
preserved when taking squares. Using the elementary inequality 
$(x - y)^2 \geq (1 - \frac{1}{\alpha}) x^2 - (\alpha - 1) y^2$ for any $x,y \geq 0$ and $\alpha > 0$, we thus obtain that 
\begin{equation}\label{eq:Sigma0sqv_lower}
\nnorm{\wh{\Sigma}_0^{1/2} v}_2^2 \geq \left(1 - \frac{1}{\alpha} \right)  v^{\T} \Sigma_0 v  - (\alpha - 1) \overline{\Gamma}^2(v). 
\end{equation}
Substituting this lower bound into \eqref{eq:RE_basic_1}, we obtain that for any $\alpha > 1$
\begin{align}\label{eq:RE_lower_final}
v^{\T} \wh{\Sigma} v &\geq \left(1 - \frac{1}{\alpha} \right)  v^{\T} \Sigma_0 v
- v^{\T} D v - \big| \{ v^{\T} (D -
                         \wh{D})  v \} \big|  - (\alpha - 1)
                         \overline{\Gamma}^2(v)  \notag \\
&\geq  \left(1 - \frac{1}{\alpha} \right) v^{\T} \Sigma v  - \frac{1}{\alpha} v^{\T} D v   -  \big| \{ v^{\T} (D -
                         \wh{D})  v \} \big| - (\alpha - 1)
                                        \overline{\Gamma}^2(v) \notag \\
&\geq \left[ \left(1 - \frac{1}{\alpha} \right) \lambda_{\min}(\Sigma)  - \frac{1}{\alpha} \max_{1
                                                      \leq j \leq d}
                                                      \{ R^2 -
                                                      \Sigma_{jj} \} \right] \nnorm{v}_2^2
                                                      -  \big| \{ v^{\T} (D -
                         \wh{D})  v \} \big| - (\alpha - 1)
                                                      \overline{\Gamma}^2(v)  \notag \\
&\geq \left[ \left(1 - \frac{1}{\alpha} \right) \lambda_{\min}(\Sigma)  - \frac{1}{\alpha} 
                                                      (R^2 -
                                                      \lambda_{\min}(\Sigma)) \right] \nnorm{v}_2^2
                                                      -   \big| \{ v^{\T} (D -
                         \wh{D})  v \} \big| - (\alpha - 1)
                                                      \overline{\Gamma}^2(v)  \notag \\       
&\geq \frac{1}{2} \lambda_{\min}(\Sigma) \nnorm{v}_2^2   -  \big| \{ v^{\T} (D -
                         \wh{D})  v \} \big|  - (\alpha - 1)
                                                                     \overline{\Gamma}^2(v)
\end{align}
by choosing $\alpha = \alpha_0 \coloneq \frac{2 R^2}{\lambda_{\min}(\Sigma)} \geq 2$ and noting that $R^2 \geq \max_{1 \leq j \leq d} \Sigma_{jj} \geq \lambda_{\min}(\Sigma)$. Furthermore, we have
\begin{equation*}
\big|\{ v^{\T} (D -
                         \wh{D})  v \} \big| \leq \nnorm{v}_2^2 \max_{1
                           \leq j \leq d} |D_{jj} - \wh{D}_{jj}|. 
                       \end{equation*}
Note that $D_{jj} - \wh{D}_{jj}= (R^2 - \Sigma_{jj}) - \left(R^2 -
\frac{1}{n} \su \wt{X_{ij}^2} \right)$, $1 \leq j \leq d$. Since
$\E[\wt{X_{ij}^2}] = \Sigma_{jj}$, $1 \leq j \leq d$, an application of Hoeffdings's inequality combined with the union bound implies that the event 
\begin{equation}\label{eq:event_Dfluctuation}
\left\{ \max_{1 \leq j \leq d} |D_{jj} - \wh{D}_{jj}| \leq R^2 \sqrt{\log(d)/n} \right\}
\end{equation}
occurs with probability at least $1 - 2/d$. Note that if 
\begin{equation}\label{eq:ncond_Ddev}
n \geq \frac{16 R^4}{\lambda_{\min}^2(\Sigma)} \log(d),
\end{equation}
the right hand side in \eqref{eq:event_Dfluctuation} is bounded by $\frac{1}{4} \lambda_{\min}(\Sigma)$. Combining this with \eqref{eq:RE_lower_final} yields \eqref{eq:RSC}. 

The proof of the bound (i) is concluded by collecting the conditions on $n$ per \eqref{eq:condn_tau} and \eqref{eq:ncond_Ddev}, and the probabilities associated with the use of Lemma \ref{lem:grad:sup} and the events in \eqref{eq:matrix_deviation}, \eqref{eq:event_Dfluctuation}. 
\vskip1.5ex
We now turn to bound (ii). If $\wh{\delta} \in \mc{C}$, we obtain $\nnorm{\wh{\delta}_{S^c}}_1 \leq 3 \nnorm{\wh{\delta}_S}_1$ as demonstrated above, which yields $\nnorm{\wh{\delta}}_1 \leq 4 \sqrt{s} \nnorm{\wh{\delta}}_2$, and thus the assertion. By contrast, if $\wh{\delta} \notin \mc{C}$, let 
\begin{equation*}
\Delta = \wt{C} B  \frac{\nnorm{\Sigma}_{\text{op}}^{1/2} + R}{\{ \lambda_{\min}(\Sigma) \}^{1/2}} \sqrt{\frac{s \log d}{n}}
\end{equation*}
denote the bound on the right hand side of \eqref{eq:tangentcone}. We hence obtain the following from \eqref{eq:firstorder}: 
\begin{align*}
0 &\geq  \wh{\delta}^{\T} \wh{\Sigma} \wh{\delta} - \lambda_0 \nnorm{\wh{\delta}}_1 + \lambda \nscp{\nabla
    \nnorm{\wh{\beta}}_1}{\wh{\delta}} \\
    &\geq \wh{\delta}^{\T} (\wh{\Sigma}_0 - \wh{D}) \wh{\delta} - \lambda_0 \nnorm{\wh{\delta}}_1 + \lambda (\nnorm{\wh{\beta}}_1 - \nnorm{\beta_*}_1) \\
    & \geq -\wh{\delta}^{\T} \wh{D} \wh{\delta} - \lambda_0 \nnorm{\wh{\delta}}_1 + \lambda (\nnorm{\wh{\beta}}_1 - \nnorm{\beta_*}_1) \\
    &\geq -R^2 \nnorm{\wh{\delta}}_2^2 - \lambda_0 \nnorm{\wh{\delta}}_1 + \lambda (\nnorm{\wh{\beta}}_1 - \nnorm{\beta_*}_1). 
\end{align*}
The last inequality implies that 
\begin{align*}(\lambda - \lambda_0) \nnorm{\wh{\delta}_{S^c}}_1 \leq (\lambda + \lambda_0) \nnorm{\wh{\delta}_{S}}_1 + R^2 \Delta^2 \; \; \iff  \;\,
\nnorm{\wh{\delta}_{S^c}}_1 &\leq 3 \nnorm{\wh{\delta}_S}_1 + \frac{R^2 \Delta^2}{\lambda - \lambda_0} \\
&\leq 3 \nnorm{\wh{\delta}_S}_1 + \frac{R^2 \Delta^2}{\lambda_0} \\
&\leq  3 \nnorm{\wh{\delta}_S}_1 + \Delta \frac{\wt{C} B (\nnorm{\Sigma}_{\text{op}}^{1/2} + R) \sqrt{s}}{\{ \lambda_{\min}(\Sigma ) \}^{1/2} \cdot C},
\end{align*}
from the definitions of $\Delta, \lambda$ and $\lambda_0$, where $C > 0$ is the 
absolute constant in the definition of $\lambda_0$. Setting $\wt{C}' = \wt{C}/C$, bound (ii) in the case
$\wh{\delta} \notin \mc{C}$ then follows from 
\begin{equation*}
\nnorm{\wh{\delta}}_1  \leq 4 \nnorm{\wh{\delta}_S}_1 + \Delta \frac{\wt{C} B (\nnorm{\Sigma}_{\text{op}}^{1/2} + R)\sqrt{s}}{\{ \lambda_{\min}(\Sigma ) \}^{1/2} \cdot C} 
\leq \left(4 +  \frac{\wt{C}' B (\nnorm{\Sigma}_{\text{op}}^{1/2} + R)}{\{ \lambda_{\min}(\Sigma ) \}^{1/2}}  \right) \sqrt{s} \Delta.  
\end{equation*}

\subsection*{Proof of Proposition \ref{prop:debias_lasso}}
From the definition of $\wh{\beta}^{\text{{db}}}$ and re-arranging terms, we obtain that 
\begin{equation*}
\wh{\beta}^{\text{db}} - \beta_* = (M_n\wh{\Sigma}  - I)(\beta_* - \wh{\beta}) + M_n(\wh{\Sigma}_{Xy} - \wh{\Sigma} \beta_*). 
\end{equation*}
Note that under the assumptions made and by invoking Theorem \ref{theo:lasso}
\begin{align*}
\nnorm{(M_n\wh{\Sigma}  - I)(\beta_* - \wh{\beta})}_{\infty} \leq 
\nnorm{M_n\wh{\Sigma}  - I}_{\infty} \nnorm{\wh{\beta} - \beta_*}_1 = O_{\p}\left( \sqrt{\frac{\log p}{n}} \right) O_{\p} \left(s \sqrt{\frac{\log p}{n}} \right)
\end{align*}
Thus, $\sqrt{n} \nnorm{(M_n\wh{\Sigma}  - I)(\beta_* - \wh{\beta})}_{\infty} = o_{\p}(1)$ since $s = o(\sqrt{n}/\log d)$. It follows that 
\begin{equation*}
\sqrt{n} (\wh{\beta}_j^{\text{db}} - \beta_j^*) = \sqrt{n} (m_{j,n}^{\T} (\wh{\Sigma}_{Xy} -\wh{\Sigma} \beta_*))+ o_{\p}(1), \quad m_{j,n} = M_n^{\T} e_j, \;\;\;\; j=1,\ldots,d. 
\end{equation*}
Observe that the first term on the right hand side can be expressed as 
\begin{align*}
    \sqrt{n} (m_{j,n}^{\T} (\wh{\Sigma}_{Xy} -\wh{\Sigma} \beta_*)) =
    \frac{1}{\sqrt{n}} \su \zeta_i, \;\;\; \zeta_i \coloneq &\sum_{k=1}^{d} (m_{j,n})_k (\wt{X}_{ik} \wt{Y}_i - e_k^{\T} \{\wt{X}_i \wt{X}_i^{\T} + \Delta_i\} \beta_*),\\ 
    &1 \leq i \leq n. 
\end{align*}
Note that by assumption for all $1 \leq j \leq d$, we have $\var(\zeta_i) \rightarrow \vartheta_j$, $1 \leq i \leq n$, as $n \rightarrow \infty$. The assertion hence follows from the central limit theorem and Slutsky's theorem.  


\subsection{Additional technical Lemmas}\label{app:add_tech}

\begin{lemma}\label{lem:grad:sup} In the situation of Theorem \ref{theo:lasso}, there exists a universal constant $C > 0$ such that 
\begin{equation*}
\nnorm{\wh{\Sigma} \beta_* - \wh{\Sigma}_{Xy}}_{\infty} \leq \lambda_0 \invcoloneq  C  (L R + R^2 B) \left( \sqrt{\frac{\log d}{n}} + \frac{\log d}{n} \right)
\end{equation*}
with probability at least $1 - 2/d$. 
\end{lemma}
\begin{bew} Pick $j \in \{1,\ldots,d\}$ arbitrary and consider the $j$-th coordinate of $\wh{\Sigma} \beta_* - \wh{\Sigma}_{Xy}$. We have
\begin{align}\label{eq1:lem:grad_sup}
(\wh{\Sigma} \beta_* - \wh{\Sigma}_{Xy})_j &= \frac{1}{n} \su \left( \wt{X}_{ij} (\wt{X}_i^{\T} \beta_*)  + (\wt{X^2_{ij}}  - R^2)\beta_{*j} - \wt{X}_{ij} \wt{Y}_i \right) \\
 &=\frac{1}{n} \su \Big( \{ \wt{X}_{ij} (\wt{X}_i^{\T} \beta_*) - \E[\wt{X}_{ij} (\wt{X}_i^{\T} \beta_*) ] \} \notag \\
 &\qquad \qquad + \{ (\wt{X^2_{ij}}  - R^2)\beta_{*j} - \E[(\wt{X^2_{ij}}  - R^2)\beta_{*j}] \}  \notag \\
  &\qquad \qquad -\{ \wt{X}_{ij} \wt{Y}_i - \E[\wt{X}_{ij} \wt{Y}_i ] \} \Big) \notag \\
  &= \frac{1}{n} \sum_{i = 1}^n (T_{1i} + T_{2i}  + T_{3i}), \notag
\end{align}
where the second identity follows from the fact that each summand in \eqref{eq1:lem:grad_sup} has mean zero. Note that 
$T_{2i} = (\wt{X^2_{ij}}  - R^2)\beta_{*j} - \E[(\wt{X^2_{ij}}  - R^2)\beta_{*j}]$ and $T_{3i} =  \wt{X}_{ij} \wt{Y}_i - \E[\wt{X}_{ij} \wt{Y}_i ]$, $1 \leq i \leq n$, are zero-mean sub-Gaussian random variables whose sub-Gaussian norms are proportional to $R^2 |\beta_{*j}| \leq R^2 B$ and $L \cdot R$, respectively.  
Moreover, the $T_{1i} = \wt{X}_{ij} (\wt{X}_i^{\T} \beta_*) = T_{11i} \cdot (T_{12i})$ are products 
of zero-mean sub-Gaussian random variables with sub-Gaussian norms proportional
to $R$ and $R \nnorm{\beta_*}_2$, respectively, thus the $T_{1i}$ are zero-mean sub-exponential
random variables with sub-Exponential norm proportional to $R^2 \nnorm{\beta_*}_2 \leq R^2 B$ as follows from Lemma \ref{ref:product_subG}, $1 \leq i \leq n$. It follows that $T_{1i} + T_{2i} + T_{3i}$ are also zero-mean sub-Exponential random variables with sub-Exponential norm bounded by $ L \cdot R + 2 R^2 B$, $1 \leq i \leq n$. By Bernstein's inequality for sub-Exponential random variables (Lemma \ref{lem:Bernstein_subexp}), we have for any $t > 0$ 
\begin{equation}\label{eq2:lem:grad_sup}
\p\left(\left|\frac{1}{n} \su (T_{1i} + T_{2i} + T_{3i})  \right| \geq t \right)
\leq 2 \exp\left(-c \left(\frac{n t^2}{L^2 R^2 + R^4 B^2} \wedge \frac{n t}{L \cdot R + R^2 B} \right) \right),
\end{equation}
where $c > 0$ is a universal constant. Choosing 
\begin{equation*}
t = \frac{2}{\sqrt{c} \wedge c} (L R + R^2 B) \left( \sqrt{\frac{\log d}{n}} + \frac{\log d}{n} \right) = C  (L R + R^2 B) \left( \sqrt{\frac{\log d}{n}} + \frac{\log d}{n} \right),
\end{equation*}
where $C = \frac{2}{\sqrt{c} \wedge c}$, the right hand side of \eqref{eq1:lem:grad_sup} is upper bounded by $2/d^2$, and the assertion of the lemma follows from a union bound over $\{1,\ldots,d\}$.

\end{bew}

\subsubsection*{Properties of sub-Gaussian and sub-Exponential random variables} 
The following statements concern basic properties of sub-Gaussian random variables that 
can be found in standard literature such as \cite{Vershynin2018}. 
\renewcommand{\thelemma}{A.\arabic{lemma}}
\setcounter{lemma}{0}
\begin{lemma}\label{ref:tailbound} Let $X$ be zero-mean, unit variance sub-Gaussian random variable with $K = \nnorm{X}_{\psi_2}$. Then for any $t > 0$ 
\begin{equation*}
\p(|X| \geq t) \leq 2 \exp(-t^2/(2 C_K)), 
\end{equation*}
where $C_K = C K^2$ for some constant $C > 0$. 
\end{lemma}

\begin{lemma}\label{ref:product_subG}\cite[][Lemma 2.7.7]{Vershynin2018}
Let $X$ and $X'$ be sub-Gaussian random variables. Then 
$X \cdot X'$ is sub-Exponential with $\nnorm{X \cdot X'}_{\psi_1} \leq 
\nnorm{X}_{\psi_2} \cdot \nnorm{X'}_{\psi_2}$
\end{lemma}

\begin{lemma}\label{ref:linearcomb_subG}
Let $(X_i)_{i = 1}^m$ be independent zero-mean sub-Gaussian random variables, and 
let $(a_i)_{i = 1}^m \subset \R$ arbitrary. Then $\nnorm{\sum_{i = 1}^m a_i X_i}_{\psi_2}  \leq \max_{1 \leq i \leq m} C \nnorm{X_i}_{\psi_2}$ for some constant $C > 0$. 
\end{lemma}

\begin{lemma}\label{lem:Bernstein_subexp} (Bernstein's inequality for sub-Exponential random variables) \cite[][Theorem 2.8.1]{Vershynin2018}
Let $X_1, \ldots, X_n$ be independent, zero mean, sub-Exponential random variables such that $\max_{1 \leq i \leq n} \nnorm{X_i}_{\psi_1} \leq K$. Then, for every $t \geq 0$, we have 
\begin{equation*}
\p\left(\left| \frac{1}{n} \su X_i \right| \geq t \right) \leq 2 \exp \left(-c  \left\{ \frac{n t^2}{K^2} \wedge \frac{n t}{K} \right\} \right),
\end{equation*}
where $c > 0$ is an absolute constant. 
\end{lemma}

\subsubsection*{Matrix Bernstein inequality} 

The following Bernstein bound for random matrices is taken from \cite{Wainwright2019}
\begin{lemma}\label{lem:matrix_bernstein} Let $\{ Q_i \}_{i = 1}^n$ be a sequence of independent, zero-mean, symmetric random 
matrices of dimension $d$ so that $\max_{1 \leq i \leq n} \nnorm{Q_i}_{\text{{\em op}}} \leq b < \infty$ almost surely. Moreover, let $\sigma^2 = \nnorm{\frac{1}{n} \su \E[Q_i^2]}$. Then for all $\delta > 0$, it holds that  
\begin{equation*}
\p \left(\norm{\frac{1}{n} \su Q_i}_{\text{{\em op}}} \geq \delta \right) \leq 2d \exp \left(-\frac{n\delta^2}{2(\sigma^2  + b \delta)} \right).
\end{equation*}
\end{lemma}

\subsubsection*{Matrix deviation inequality} 

The lemma below is called ``matrix deviation inequality" in
\cite[][cf.~Theorem 9.1.1 and Exercise 9.1.8]{Vershynin2018}. 
We note that in the version we consider only sets $T$ that
are symmetric around the origin, hence the Gaussian complexity
of $T$ in the original statement may be replaced by the Gaussian width
of $T$. 

\begin{lemma}\label{lem:matrix_deviation}
Let $T = -T$ be a compact subset of $\R^d$, and let 
$w(T) = \E_{g \sim N(0, I_d)}[\sup_{t \in T} \nscp{w}{g}]$ be
the Gaussian width of $T$, and let further $\text{{\em rad}}(T) = \sup_{t \in T} \nnorm{t}_2$. Let $A$ be an $n \times d$ random matrix whose rows
$(A_i)$ are independent, isotropic, and sub-Gaussian random vectors in 
$\R^d$ with $\max_{1 \leq i \leq n} \nnorm{A_i}_{\psi_2} \leq K$. We then
have for any $u > 0$
\begin{equation*}
\p \left(\sup_{t \in T} \left|\norm{\frac{A t}{\sqrt{n}}} - \nnorm{t}_2 \right| > C K^2 \frac{\left[w(T) + u \cdot \text{\em{rad}$(T)$} \right]}{\sqrt{n}} \right) \leq 2\exp\left(-u^2 \right).  
\end{equation*}

\end{lemma}


\end{document}

%% file: preamble_report.tex
\usepackage{amsfonts}
\usepackage{amssymb}
\usepackage{amsmath}
\usepackage{amsthm}
\usepackage{latexsym}
\usepackage{verbatim}
\usepackage{epsfig}
\usepackage{amsbsy}
\usepackage{amscd}
\usepackage{bm}

\usepackage{graphicx}

\usepackage{ifthen} 
\usepackage{xcolor}
\usepackage[colorlinks=true, linkcolor=teal, anchorcolor=teal,
            citecolor=blue,  filecolor=blue, menucolor=blue, pagecolor=teal,											urlcolor=blue, plainpages=false, pdfpagelabels]{hyperref}

\newcommand{\mymail}[1]{\href{mailto:#1}{\texttt{#1}}}

\usepackage{fancyhdr}
\usepackage[realmainfile]{currfile}
\fancypagestyle{firststyle}{

\cfoot{{\footnotesize $^\ast$: Supported by a scholarship of the Air Force Office for Scientific Research}} 
}

\newcommand{\setauthA}[1]{\def\authA{#1}}
\newcommand{\setauthB}[1]{\def\authB{#1}}

\def\printA{\begin{tabular}{l} \authA \end{tabular}}
\def\printB{\begin{tabular}{l} \authB \end{tabular}}

\newcommand{\makemytitle}[1]{\begin{center}{\textsf{\LARGE #1}}
  \end{center}
}

\usepackage[sf,bf]{titlesec}

\usepackage{algorithmic}
\usepackage{algorithm}

\usepackage[nonamebreak,square,numbers,sort]{natbib}

\usepackage[letterpaper, textheight = 676pt]{geometry}
\providecommand{\M}[1]{\mathbf#1}
\providecommand{\mc}[1]{\mathcal#1}
\providecommand{\mc}[1]{\mathcal#1}

\newcommand{\R}{{\mathbb R}}

\DeclareMathOperator{\E}{\mathbf{E}}
\DeclareMathOperator{\p}{\mathbf{P}}

\DeclareMathOperator{\var}{Var}
\DeclareMathOperator{\cov}{Cov}

\DeclareMathOperator{\tr}{tr}

\providecommand{\T}{\top} 

\DeclareMathOperator*{\argmin}{argmin}

\providecommand{\wt}[1]{\widetilde{#1}}
\providecommand{\wh}[1]{\widehat{#1}}

\providecommand{\norm}[1]{\left \lVert#1 \right \rVert}
\providecommand{\nnorm}[1]{ \lVert#1 \rVert}

\newcommand{\nscp}[2]{\langle#1, #2\rangle}

\newcommand{\dev}[2]{\Big\lvert _{{#1}={#2}}}

\newcommand{\bou}[2]{\Big\lvert _{{#1}}^{#2}}

\newcommand{\blanco}[1]{  }

\newcommand{\deriv}[3]{%
\ifthenelse{#1 = 1}{\frac{d\,#2}{d\,#3}}{\frac{d^{{#1}} #2}{d{#3}^{{#1}}}}
}

\newcommand{\partials}[3]{%
\ifthenelse{#1 = 1}{\frac{\partial\,#2}{\partial\,#3}}{\frac{\partial^{#1}
    #2}{\partial#3^{#1}}}
} 

\def\su{\sum_{i=1}^n}

\def \coloneq{\mathrel{\mathop:}=}
\def \invcoloneq{=\mathrel{\mathop:}}
\def \eps{\varepsilon}

\newtheorem{theo}{Theorem}
\newtheorem{propo}{Theorem}

\newtheorem{prop}[propo]{Proposition}

 \newtheorem{lemma}{Lemma}

\newenvironment{bew}{\begin{proof}[Proof]}{\end{proof}}
